\documentclass[11pt,a4paper]{article}
\usepackage{amsmath}
\usepackage{amsfonts}
\usepackage{amssymb}
\usepackage{amsthm}
\usepackage{verbatim}
\numberwithin{equation}{section}

\newtheorem{Thm}{Theorem}[section]
\newtheorem{Prop}{Proposition}[section]
\newtheorem{Lem}{Lemma}[section]
\newtheorem{Cor}{Corollary}[section]
\theoremstyle{definition}
\newtheorem{Dfn}{Definition}[section]
\newtheorem{Exmp}{Example}[section]
\newtheorem{Rmk}{Remark}[section]

\def\Bosonnormalord{\,\lower.8ex \hbox{$\circ$} \llap{\raise.8ex\hbox{$\circ$}} \,}
\def\normalord{\,\lower.8ex \hbox{$\cdot$} \llap{\raise.8ex\hbox{$\cdot$}} \,}
\def\Hom{\mathop{\rm Hom}\nolimits}
\def\Aut{\mathop{\rm Aut}\nolimits}
\def\charge{\mathop{\rm charge}\nolimits}

\def\cL{{\mathcal L}}
\def\cV{{\mathcal V}}
\def\ab{\mathop{\rm ab}\nolimits}
\def\cW{{\mathcal W}}
\def\cVd{{\mathcal V}^{\dagger}}
\def\cWd{{\mathcal W}^{\dagger}}
\def\ad{\mathop{\rm ad}\nolimits}
\def\charge{\mathop{\rm charge}\nolimits}
\def\Res{\mathop{\rm Res}\limits}
\def\res{\mathop{\rm Res}\nolimits}

\def\Ker{\mathop{\rm Ker}\nolimits}
\def\Coker{\mathop{\rm Coker}\nolimits}

\def\bC{{\mathbf C}}
\def\bP{{\mathbf P}}
\def\bR{{\mathbf R}}
\def\bZ{{\mathbf Z}}
\def\cB{{\mathcal B}}
\def\cC{{\mathcal C}}
\def\cD{{\mathcal D}}
\def\cF{{\mathcal F}}
\def\cS{{\mathcal S}}
\def\cU{{{\mathcal U}}}
\def\cV{{{\mathcal V}}}
\def\cO{{{\mathcal O}}}
\def\cVd{{{\mathcal V}^\dagger}}
\def\cFd{{\mathcal F}^\dagger}
\def\gF{\mathfrak{F}}
\def\gX{\mathfrak{X}}
\def\hZ{{{\mathbf Z}_h}}
\def\oe{{\overline{e}}}
\def\ovpsi{{\overline{\psi}}}

\def\Id{\mathop{\rm Id}\nolimits}
\def\UGM{\mathop{\rm UGM}\nolimits}
\def\sumQj{\sum_{j=1}^NQ_j}
\def\QED{{{} \hfill QED }\par \medskip}

\begin{document}
\title{Abelian Conformal Field Theory and Determinant Bundles\thanks{
This research was conducted
 partly by the first author for the Clay Mathematics Institute at University of California,
 Berkeley and for MaPhySto --
Centre for Mathematical Physics and Stochastics, funded by The
Danish National Research Foundation. The second author is partially supported by Grant in Aid
for Scientific Research {\sc no}. 14102001 of JSPS.}
}
\date{ }
\author{J{\o}rgen Ellegaard Andersen \& Kenji Ueno}

\maketitle

\begin{abstract}
Following \cite{KNTY} we study a so-called
$bc$-ghost system of zero conformal dimension from the viewpoint
of \cite{TUY} and \cite{U2}. We show that the ghost vacua construction results
in holomorphic line bundles with connections over holomorphic families of curves.
We prove that the curvature of these connections are up to a scale the
same as the curvature of the connections constructed in \cite{TUY} and \cite{U2}.
We study the sewing construction for nodal curves and its explicit relation to the
constructed connections. Finally we construct preferred holomorphic sections
of these line bundles and analyze their behaviour near nodal
curves. These results are used in \cite{AU2} to construct modular
functors
form the conformal field theories given in \cite{TUY} and \cite{U2} by twisting with an
appropriate factional power of this Abelian theory.
\end{abstract}

\tableofcontents

\section*{Introduction}

The present paper is the first in a series of three papers
(\cite{AU2} and \cite{AU3}), in which
we shall construct modular functors and the
Reshetikhin-Turaev Topological Quantum Field Theories
from the conformal field theories developed in \cite{TUY} and \cite{U2}.

The basic idea behind the construction of a modular functor
for a simple Lie algebra is the following.
The sheaf of vacua construction for a simple Lie algebra gives a vector
bundle with connection over Teichm\"uller space of any oriented pointed surface.
The vector  space that the modular functor associates to the oriented pointed surface should be
the covariant constant sections of the bundle.
However the connection on the bundle is only projectively flat, so
we need to find a suitable line bundle with a connection, such that
the tensor product of the two has a flat connection.

We shall construct a line bundle with a connection on any family of
$N$-pointed curves with formal coordinates. By computing the curvature of this
line bundle, we conclude that we actually need a fractional power of
this line bundle so as to obtain a flat connection after
tensoring . In order to functorially extract this fractional power, we need to
construct a preferred section of the line bundle.

We shall construct the line bundle by the use of the so-called
$bc$-ghost systems (Faddeev-Popov ghosts) first introduced in
covariant quantization \cite{FP}. The $bc$ system have two
anticommuting fields $b(z)$, $c(z)$ of conformal dimension $j$,
$1-j$, respectively, where $j$ is an integer or half integer. In
the case $j =1/2$ a mathematically rigorous treatment was given in
the paper \cite{KNTY}. The case $j=1/2$ corresponds to the study
of the determinant bundle of half-canonical line bundles on smooth
curves, i.e.. on compact Riemann surfaces. Since we cannot define
the half-canonical line bundles for curves with node, whose
normalization has at least two components, the boundary behavior
of the sheaf of vacua is complicated \cite{KSUU}. Therefore, in
the present paper we shall consider the case $j=0$, following the
ideas of \cite{KNTY}, but describing it from the viewpoint of
\cite{TUY}. We shall also use the terminology of \cite{TUY} and
\cite{U2}. In particular we shall define the sheaf of ghost vacua
on a family of $N$-pointed Riemann surfaces with formal
coordinates and we shall
 introduce a  connection on it.  All the necessary properties
 which we need to construct for our modular functor construction in \cite{AU2} will be proved
 in this paper.
 The sheaf of the ghost vacua is isomorphic to the invertible sheaf associated
to the determinant bundle of the relative canonical sheaf of the family.

Let us explain briefly the contents of the present paper. In section 1
we shall introduce the fermion operators and the fermion Fock space
and fix the notation which will be used in the present paper.

In  section 2 the universal Grassmann manifold due to M. Sato will be defined.
The universal Grassmann  manifold is an infinite dimensional manifold which can be
embedded into an infinite dimensional projective space. The pull-back of
the hyperplane line bundle of the projective space is the determinant
bundle of the universal subbundle of the universal Grassmann manifold.

In section 3 we shall develop the theory of  $j=0$ ghost systems. The main purpose of
this section is to define the space of ghost vacua for an $N$-pointed
curve with formal coordinates.
Also basic properties of the space of ghost vacua will be discussed and we proved
important theorems such as propagation of vacua.  The space of ghost vacua is
a one-dimensional vector space. This will be proved in section 4 and section 5.

In section 4 we shall consider a family of $N$-pointed curves with formal coordinates
and define the sheaf of ghost vacua  attached to the family.  The projectively flat
connection will be defined on the sheaf of ghost vacua.

In section 5 we shall consider  smoothings  of nodal curves. In particular we shall
construct a section of
the sheaf of ghost vacua starting from an element of the ghost vacua of the nodal curve.
This construction is called sewing and it is the key to proving that
the space ghost vacua is a one-dimensional vector space.

In section 6 we shall construct the preferred element of an $N$-pointed curve
with formal coordinates and study its basic properties.

\section{Fermion Fock space}

Let $\hZ$ be the set of all half integers. Namely
$$
\hZ = \{ n + 1/2\, | \, n \in \bZ \,\}.
$$
Let $\cWd$ be an infinite-dimensional vector space over $\bC$
with a filtration $\{F^m\cWd\}_{m \in \bZ}$ which satisfies the following conditions.
\begin{enumerate}
\item  The filtration $\{F^m\cWd\}$ is decreasing;
\item  $\bigcup_{m \in \bZ} F^m\cWd = \cWd$, \quad
$\bigcap_{m \in \bZ} F^m\cWd = \{0\}$;
\item $\dim_\bC F^m\cWd/ F^{m+1}\cWd = 1$;
\item The vector space $\cWd$ is complete with respect to the uniform topology such that
$\{F^m\cWd\}$ is a basis of open neighbourhoods of 0.
\end{enumerate}
We introduce a  basis $\{e^\nu \}_{\nu \in \hZ}$ of $\cWd$ in such a way that
$$
e^{m+1/2} \in F^m\cWd \setminus F^{m+1}\cWd.
$$
Then, each element $ u \in \cWd$ can uniquely be expressed in the form
$$
u= \sum_{\nu> n_0, \nu \in \hZ}^\infty a_\nu  e^\nu
$$
for some $n_0$ and with respect to this basis the filtration is given  by
$$
F^m\cWd = \left\{ u \in \cWd \, \left|
\, u=  \sum_{\nu> m, \nu \in \hZ}^\infty a_\nu  e^\nu\, \right. \right\}.
$$
We fix the basis $\{e^\nu\}_{\nu \in \hZ}$ throughout the present
paper.

Let $\bC((\xi))$ be a field of formal Laurent series over the complex number field.
Then the basis gives us a filtration preserving linear isomorphism
\begin{eqnarray*}
\bC((\xi)) & \cong & \cWd \\
\xi^n  &\mapsto &e^{n+1/2}.
\end{eqnarray*}
By mapping $\xi^n d\xi$ to $e^{n+1/2}$ we of course also get a filtration
preserving linear isomorphism between $\bC((\xi))d\xi$ and $\cWd$.

We let $\{\oe_\nu\}_{\nu \in \hZ}$ be the dual basis of $\{e^\nu\}_{\nu \in \hZ}$.
Then,  put
$$
\cW = \bigoplus_{\nu \in \hZ} \bC\oe_\nu
$$
Then $\cW$ is the topological dual of the vector space $\cWd$.
There is a natural pairing $(\phantom{X}|\phantom{X}) : \cWd\times \cW \rightarrow \bC$
defined by
$$
(e^\nu|\oe_\mu) = \delta^\nu_\mu.
$$
In other word we have
$$
(u|v) = v(u).
$$
The complete topological vector space $\cWd$ will be used in the next section to
define the universal Grassmann manifold due to M. Sato ([SA]).

Here let us introduce the semi-infinite exterior product of the vector spaces
$\cW$ and $\cWd$.  For that purpose we first introduce the notion of
a Maya diagram.

\begin{Dfn}{\rm
A Maya diagram $M$ of the charge $p$, $p \in \bZ$ is a set
$$
M = \left\{ \mu(p-1/2), \mu(p-3/2), \mu(p-5/2), \ldots\right\},
$$
where $\mu$ is an increasing function
$$
\mu : \hZ_{<p} = \{ \nu \in \hZ\, | \, \nu < p\,\} \rightarrow \hZ
$$
such that there exists an integer $n_0$ such that
$$
\mu(\nu) = \nu
$$
for all $\nu < n_0$.

The function $\mu$ is called the {\it characteristic function} of the
Maya diagram $M$.
The set of Maya diagrams of charge $p$ is written as ${\mathcal M}_p$.
}
\end{Dfn}

A Maya diagram $M$ of charge $p$ with its characteristic function $\mu$ can
be expressed by a diagram as in Figure 1.

\setlength{\unitlength}{1mm}
\begin{figure}
\begin{picture}(100,50)(-10,-10)
\put(15,-10){$\mu(\nu)=\nu$ for $\nu\le -\frac52$, $\mu(-\frac32) = -\frac12$,
$\mu(-\frac12) = \frac23$, $\mu(\frac12) = \frac52$}
\put(5,0){\line(1,0){99}}
\put(5,10){\line(1,0){99}}
\put(10,0){\line(0,1){10}}
\put(20,0){\line(0,1){10}}
\put(30,0){\line(0,1){10}}
\put(40,0){\line(0,1){10}}
\put(50,0){\line(0,1){10}}
\put(60,0){\line(0,1){10}}
\put(70,0){\line(0,1){10}}
\put(80,0){\line(0,1){10}}
\put(90,0){\line(0,1){10}}
\put(100,0){\line(0,1){10}}
\put(5,5){\line(1,1){5}}
\put(5,0){\line(1,1){5}}
\put(10,0){\line(1,1){10}}
\put(10,5){\line(1,1){5}}
\put(10,8){\line(1,1){2}}
\put(15,0){\line(1,1){5}}
\put(18,0){\line(1,1){2}}
\put(20,0){\line(1,1){10}}
\put(20,5){\line(1,1){5}}
\put(20,8){\line(1,1){2}}
\put(25,0){\line(1,1){5}}
\put(28,0){\line(1,1){2}}
\put(30,0){\line(1,1){10}}
\put(30,5){\line(1,1){5}}
\put(30,8){\line(1,1){2}}
\put(35,0){\line(1,1){5}}
\put(38,0){\line(1,1){2}}
\put(50,0){\line(1,1){10}}
\put(50,5){\line(1,1){5}}
\put(50,8){\line(1,1){2}}
\put(55,0){\line(1,1){5}}
\put(58,0){\line(1,1){2}}
\put(50,0){\line(1,1){10}}
\put(50,5){\line(1,1){5}}
\put(50,8){\line(1,1){2}}
\put(55,0){\line(1,1){5}}
\put(58,0){\line(1,1){2}}
\put(70,0){\line(1,1){10}}
\put(70,5){\line(1,1){5}}
\put(70,8){\line(1,1){2}}
\put(75,0){\line(1,1){5}}
\put(78,0){\line(1,1){2}}
\put(80,0){\line(1,1){10}}
\put(80,5){\line(1,1){5}}
\put(80,8){\line(1,1){2}}
\put(85,0){\line(1,1){5}}
\put(88,0){\line(1,1){2}}
\put(12, 13){$-\frac{9}{2}$}
\put(22, 13){$-\frac{7}{2}$}
\put(32, 13){$-\frac{5}{2}$}
\put(42, 13){$-\frac{3}{2}$}
\put(52, 13){$-\frac{1}{2}$}
\put(63, 13){$\frac{1}{2}$}
\put(73, 13){$\frac{3}{2}$}
\put(83, 13){$\frac{5}{2}$}
\put(93, 13){$\frac{7}{2}$}
\end{picture}
\caption{Maya diagram of charge 1 and degree 5}
\end{figure}

For a Maya diagram M we have $\mu(\nu) = \nu$ for almost all $\nu$. Therefore
the set
$$
\{ \mu(\nu) - \nu \,|\, \nu \in \hZ, \, \mu(\nu) - \nu >0\,\}
$$
is finite and the number
$$
d(M) = \sum_{\nu \in \hZ} (\mu(\nu) - \nu)
$$
is finite. The number $d(M)$ is also written as $d(\mu)$ and it is
called the {\it degree} of the Maya diagram $M$ with characteristic function  $\mu$.
The finite set of Maya diagrams of degree $d$ and change $p$ is
denoted ${\mathcal M}_p^d$. Clearly ${\mathcal M}_p = \coprod_{d} {\mathcal
M}_p^d$.

For a Maya diagram $M$ of charge $p$ we define two semi-infinite
products
\begin{eqnarray*}
|M\rangle &= & \oe_{\mu(p-1/2)}\wedge \oe_{\mu(p-3/2)} \wedge \oe_{\mu(p-5/2)}
\wedge \cdots \\
\langle M| &=& \cdots \wedge
 e^{\mu(p-5/2)}\wedge  e^{\mu(p-3/2)} \wedge  e^{\mu(p-1/2)}
\end{eqnarray*}
Formally, these semi-infinite products is just another notation
for the corresponding Maya diagram. This notation is
particular convenient for the following discussion. However, by using
the basis $e_\nu$, we clear indicate the relation to the vector
spaces $\cW$ and $\cWd$.

For any integer $p$ put
\begin{eqnarray*}
|p \rangle &= & \oe_{p-1/2}\wedge \oe_{p-3/2} \wedge \oe_{p-5/2}
\wedge \cdots \\
\langle  p| &=& \cdots \wedge
 e^{p-5/2}\wedge  e^{p-3/2} \wedge  e^{p-1/2}
\end{eqnarray*}

Now the {\it fermion Fock space} $\cFd(p)$
of {\it charge} $p$ and the {\it dual fermion Fock space} $\cF(p)$ of
{\it charge} $p$ are defined by
\begin{eqnarray*}
\cF(p) &=&  \bigoplus_{M \in {\mathcal M}_p} \bC |M\rangle \\
\cFd(p) &=& \prod_{M \in {\mathcal M}_p} \bC \langle M|
\end{eqnarray*}
We observe that
\[\cF(p) = \bigoplus_{d\geq 0} \cF_d(p),\]
where
\[\cF_d(p) =  \bigoplus_{M\in {\mathcal M}_p^d} \bC |M\rangle. \]

The dual pairing
$$
\langle \cdot | \cdot \rangle : \cFd(p) \times \cF(p) \rightarrow \bC
$$
is given by
$$
\langle M| N \rangle = \delta_{M,N}, \quad M,  N \in {\mathcal M}_p
$$
Put also
\begin{eqnarray*}
\cF      &=&  \bigoplus_{p \in \bZ} \cF(p) \\
\cFd   &=& \bigoplus_{p \in \bZ} \cFd(p)
\end{eqnarray*}
The vector space $\cFd$ is called the {\it fermion Fock space} and $\cF$ is
called the {\it dual fermion Fock space}.
These are the
semi-infinite exterior products of the vector spaces $\cWd$ and
$\cW$ respectively, which we shall be interested in. We only
define the fermion Fock space by using the basis $e_\nu$,
since we are fixing this basis throughout.

The above pairing can be extended to the one on $\cFd \times \cF$ by
assuming that the paring is zero on $\cFd(p) \times \cF(p')$ if $p \ne p'$.

Let us introduce the {\it fermion operators} $\psi_\nu$ and $\ovpsi_\nu$  for
all half integers $\nu \in \hZ$ which act on $\cF$ from the left  and on
$\cFd$ from the right.
\begin{eqnarray}
\hbox{\rm Left action on $\cF$} &\quad & \psi_\nu = i(\oe_\nu), \quad
\ovpsi_\nu = \oe_{-\nu} \wedge \\
\hbox{\rm Right action on $\cFd$} &\quad & \psi_\nu = \wedge  e^\nu, \quad
\ovpsi_\nu = i(e^{-\nu})
\end{eqnarray}
where $i(\cdot)$ is the interior product. For example we have
\begin{eqnarray*}
\psi_{-3/2}|0\rangle &=& i(\oe_{-3/2}) \oe_{-1/2}\wedge \oe_{-3/2}\wedge \cdots =
- \oe_{-1/2}\wedge\oe_{-5/2}\wedge \oe_{-7/2}\wedge \cdots , \\
\langle 0 | \ovpsi_{5/2} &=& \cdots \wedge e^{-5/2}\wedge e^{-3/2} \wedge e^{-1/2} i(e^{-5/2})
=\cdots \wedge e^{-7/2}\wedge e^{-3/2} \wedge e^{-1/2}
\end{eqnarray*}
Note that $\psi_\nu$ maps $\cF(p)$ to $\cF(p-1)$, hence decreases the charge by one,
and $\ovpsi_\nu$ maps $\cF(p)$ to $\cF(p+1)$, hence increase the charge by one.
Similarly the right action of $\psi_\nu$ maps $\cFd(p)$ to $\cFd(p+1)$ and
$\ovpsi_\nu$ maps  $\cFd(p)$ to $\cFd(p-1)$. It is easy to show that
for any $\langle u| \in \cF$ and $|v\rangle \in \cFd$ we have
$$
\langle u | \psi_\nu v\rangle = \langle u  \psi_\nu | v\rangle,
\quad \langle u | \ovpsi_\nu v\rangle = \langle u  \ovpsi_\nu | v\rangle .
$$
Also it is easy to show that
\begin{eqnarray*}
\psi_\nu |0 \rangle = 0  && \hbox{\rm if and only if $\nu >0$,} \\
\ovpsi_\nu |0 \rangle = 0  && \hbox{\rm if and only if $\nu >0$.}
\end{eqnarray*}
Similarly we have
\begin{eqnarray*}
\langle 0 | \psi_\nu = 0  && \hbox{\rm if and only if $\nu <0$,} \\
\langle 0 | \ovpsi_\nu = 0  && \hbox{\rm if and only if $\nu <0$.}
\end{eqnarray*}

The fermion operators have the following anti-commutation relations as operators
on $\cF$ and $\cFd$.
\begin{eqnarray}
{[\psi_\nu, \psi_\mu ]}_+  &=&  0, \label{anticomm1}\\
{[\ovpsi_{\nu}, \ovpsi_{\mu} ]}_+ &=& 0,  \label{anticomm2}\\
{[\psi_\nu, \ovpsi_{\mu} ]}_+  &=& \delta_{\nu +\mu, 0}, \label{anticomm3}
\end{eqnarray}
where we define
$$
{[A,B]}_+ = AB+BA.
$$

Note that for each Maya diagram $M$ of charge $p$ we can find
non-negative half integers
$$
\mu_1<\mu_2<\cdots<\mu_r<0, \quad
 \nu_1<\nu_2<\cdots<\nu_s<0, \quad r\ge 0, \,s\ge 0
 $$
with $r- s=p$ and $\mu_i \neq \nu_j$ such that
\begin{equation}
\label{maya}
 |M\rangle = (-1)^{\sum_{i=1}^s\nu_i  +s/2}
  \ovpsi_{\mu_1}\ovpsi_{\mu_2}\cdots\ovpsi_{\mu_r}
 \psi_{\nu_s}\psi_{\nu_{s-1}}\cdots\psi_{\nu_1}|0\rangle.
\end{equation}
The negative half integers $\mu_i$'s and $\nu_j$'s are uniquely determined by
the Maya diagram $M$.
The {\it normal ordering} $\normalord\phantom{X}\normalord$ of the
fermion operators are defined as follows.
$$
\normalord A_\nu B_\mu \normalord = \left\{
\begin{array}{ll}
- B_\mu A_\nu & \hbox{\rm if $\mu <0$ and $\nu>0$}, \\
A_\nu B_\mu & \hbox{\rm otherwise,}
\end{array}
\right.
$$
where $A$ and $B$ is $\psi$ or $\ovpsi$. By \eqref{anticomm1}, \eqref{anticomm2} and
\eqref{anticomm3} the normal ordering is non-trivial if and only if $\mu <0$ and
$A_\nu ={\bar B}_{-\mu} $.  For example we have
$$
\normalord \ovpsi_{1/2} \psi_{-1/2}\normalord = - \psi_{-1/2}\ovpsi_{1/2},
$$
although as operators on $\cF$ or $\cFd$ we have
$$
\ovpsi_{1/2} \psi_{-1/2} =  - \psi_{-1/2}\ovpsi_{1/2} + \Id.
$$
Thus the normal ordering has the effect of subtracting the identity operator,
whenever it is non-trivial.

The field operators $\psi(z)$ and $\ovpsi(z)$ are defined by
\begin{eqnarray}
\psi(z)&=& \sum_\mu \psi_\mu z^{-\mu -1/2} , \\
\overline{\psi}(z) &=& \sum_\mu \overline{\psi}_\mu z^{-\mu -1/2} .
\end{eqnarray}
The current operator $J(z)$ is defined by
\begin{equation}
J(z) = \normalord \overline{\psi}(z)\psi(z) \normalord = \sum_{n \in \bZ} J_n z^{-n-1}
\end{equation}

Note that thanks to the normal ordering, the operator $J_n$ can operate on $\cF$ and $\cFd$
even though $J_n$ is an infinite sum of  operators.

For any integer or half integer $j$ the energy-momentum tensor $T^{(j)}(z)$ is defined
by
\begin{equation}
T^{(j)}(z) = \normalord (1-j)\frac{d \psi(z)}{dz} \ovpsi(z) -
 j \psi(z)\frac{d \ovpsi(z)}{dz}\normalord = \sum_{n \in \bZ} L^{(j)}_n z^{-n-2}.
\end{equation}
Again due to the normal ordering, the coefficients $L_n^{(j)}$ operates on $\cF$ and $\cFd$.

These operators satisfy the following commutation relations.
\begin{eqnarray}
[J_n, \, J_m] &=& n \delta_{n+m, 0},  \\
\left[ L^{(j)}_n,\,  L^{(j)}_m \right] & = & (n-m)L_{n+m}^{(j)} -
\frac16(6 j^2 - 6j +1)  (n^3-n)\delta_{m+n.0}, \\
\left[ L^{(j)}_n,\, J_m \right] &=& -mJ_{n+m} -
\frac12(2j -1) (n^2+n) \delta_{n+m,0}, \\
\left[ J_n, \, \psi(z)\right] &=& -z^n \psi(z), \\
\left[ J_n, \, \ovpsi(z) \right] &=& z^n \ovpsi(z), \\
\left[ L^{(j)}_n,\, \psi(z) \right] & = & z^n\left(
z\frac{d}{dz} + j(n+1) \right) \psi(z), \\
\left[ L^{(j)}_n,\, \overline{\psi}(z) \right] & = & z^n \left (z\frac{d}{dz}
+ (1-j)(n+1)\right) \ovpsi(z).
\end{eqnarray}

Thus the set $\{L^{(j)}_n\}_{n \in \bZ}$ forms an infinite-dimensional Lie algebra
called the {\it Virasoro} algebra with central charge $c=-2(6 j^2 - 6j +1) $.

The field operators $\psi(z)$ and $\ovpsi(z)$, the current operator $J(z)$
and the energy-momentum tensor $T^{(j)}(z)$ form the so-called
spin $j$ $bc$-system or ghost system in the physics literature.
For $j=1/2$  the $bc$-system is usually called abelian conformal field theory.
In the present paper we shall discuss the $j=0$  ghost system and we shall also
refer to this case as abelian conformal field theory.

\section{Universal Grassmann Manifold}
In this section we shall briefly recall the theory of the universal Grassmann manifold
due to M. Sato (\cite{SA}). Let $\cWd$ be the vector space introduced in the
previous section.

\begin{Dfn}
The {\it universal Grassmann manifold} $\UGM^p$ of charge $p$, $p \in
\bZ$,
is the set of closed
subspace $U\subset \cWd$ such that
\begin{enumerate}
\item The kernel and the cokernel of the natural linear map $f : U \rightarrow
\cWd/ F^0\cWd$ are of finite dimension;
\item $\dim \Ker f - \dim \Coker f = p$.
\end{enumerate}
Put
$$
\UGM =  \bigsqcup_{p \in \bZ} \UGM^p
$$
\end{Dfn}

Also, we can introduce the induced filtration $F^mU$ on $U$ by
$$
F^mU = U \cap F^m\cWd.
$$
Then we have
$$
\dim F^mU/F^{m+1}U \le 1.
$$
Put
$$
M(U) = \left\{ m + 1/2 \, |\, \dim F^mU/F^{m+1}U = 1\,\right\}.
$$
It is easy to show that $M(U)$ is a Maya diagram of charge $p$.
\begin{Dfn}
For $U \in UGM^p$ a {\it frame}  $\Xi$ of $U$ is
a basis
$$
\Xi=\{\cdots, \zeta^{p-5/2}, \zeta^{p-3/2}, \zeta^{p-1/2}\}
$$
of $U$ such that there exists a half integer $\nu_0$ such that for any
$\nu \le \nu_0$
\begin{eqnarray*}
&&\zeta^\nu \in F^{\nu -1/2}U \setminus F^{\nu +1/2}U \\
&&\zeta^\nu \equiv e^\nu \pmod {F^{\nu - 1/2}U}
\end{eqnarray*}
\end{Dfn}
For a frame
$$
\Xi=\{\cdots, \zeta^{p-5/2}, \zeta^{p-3/2}, \zeta^{p-1/2}\}
$$
of $U \in \UGM^p$ we can define the semi-infinite wedge product
\begin{equation}
\cdots \wedge \zeta^{p-5/2}\wedge  \zeta^{p-3/2}\wedge \zeta^{p-1/2}
\end{equation}
as an element of the fermion Fock space $\cFd(p)$ by the following procedure.
Put $n_0 = \nu_0 -1/2$.
For $n \le n_0 $  the element $\zeta^{n+1/2} $ can be written as
$$
\zeta^{n+1/2}  = e^{n+1/2}  + \sum_{k=n+1}^\infty a^{(n)}_{k +1/2}e^{k+1/2}.
$$
For $n \le n_0$ we define the wedge product
$$
\cdots \wedge e^{n-3/2}\wedge e^{n-1/2} \wedge \zeta^{n+1/2}\wedge
\cdots \wedge \zeta^{p-1/2}
= \langle n | \wedge \zeta^{n+1/2}\wedge
\cdots \wedge \zeta^{p-1/2}
$$
as the limit of the wedge products
$$
 \langle n | \wedge \zeta^{n+1/2}_{m_{n+1}}\wedge \cdots
 \zeta^{p-1/2}_{m_p}
$$
where $\zeta^\nu_m$ is defined by
$$
\zeta^\nu_m = \sum_{k<m}b^{(\nu)}_{k+1/2}e^{k+1/2} \equiv \zeta^\nu \mod F^m\cFd.
$$
Then the wedge product
$$
\langle n-1 | \wedge \zeta^{n-1/2}\wedge
\cdots \wedge \zeta^{p-1/2}
$$
contains all the terms appearing in
$$
\langle n | \wedge \zeta^{n+1/2}\wedge
\cdots \wedge \zeta^{p-1/2}.
$$
Hence taking the limit $n \rightarrow -\infty$ we can define the semi-infinite
product
$$
\cdots \wedge \zeta^{n-3/2}\wedge \zeta^{n-1/2} \wedge \zeta^{n+1/2}\wedge
\cdots \wedge \zeta^{p-1/2}.
$$

If we use another frame
$\Xi'$ of $U$ the resulting wedge product is a non-zero constant multiple of
the wedge product defined by the frame $\Xi$.
Therefore, to each element
$U \in \UGM^p$ we can associate a one-dimensional subspace of $\cFd(p)$, hence
associated a point of the projective space $\mathbf{P}(\cFd(p)) =
\cFd(p) \setminus\{0\}/\bC^*$. We denote this 1-dimensional
subspace of $\cFd(p)$ by $\det U$.

\begin{Exmp}
\label{exmp2.1}
Let $\gX=(C;Q;\xi)$ be a one-pointed projective curve with
formal coordinate, that is we assume that
$Q$ is a smooth point of the curve $C$ and $\xi$ is  a formal coordinate
of the curve $C$ with center $Q$.  Then to each meromorphic
function $f \in H^0(C,\cO_C(*Q))$
we can associate its Laurent expansion $f(\xi)$ at $Q$ with respect to the coordinate $\xi$.
This give an injective linear map from $H^0(C,\cO_C(*Q))$ to $\bC((\xi))$. We identify
$H^0(C,\cO_C(*Q))$ with the image and denote it by $U(\gX)$.
Identifying $\bC((\xi))$ with $\cWd$ as above, we can show that $U(\gX)$ is closed subspace of
$\bC((\xi))$. Under this identification $\cWd/F^0\cWd$ becomes identified with
  $\bC[\xi^{-1}]\xi^{-1}$.  The
natural map  $f : U(\gX) \rightarrow \bC[\xi^{-1}]\xi^{-1}$ is nothing but the map given by
taking the principal part of the Laurent expansion of each meromorphic function. Hence
$\Ker f = H^0(C,\cO_C)$. On the other hand by the exact sequence of sheaves
$$
0 \rightarrow \cO_C \rightarrow \cO_C(*Q) \rightarrow
\bC[\xi^{-1}] \xi^{-1} \rightarrow  0
$$
we obtain an  exact sequence
\begin{equation}
\label{longex}
0 \rightarrow H^0(C,\cO_C)
\rightarrow H^0(C, \cO_C(*Q)) \rightarrow \bC[\xi^{-1}]\xi^{-1}
\rightarrow H^1(C, \cO_C)
\rightarrow 0.
\end{equation}
Therefore, we have
$$
\Coker f \cong H^1(C, \cO_C) .
$$
Put
$$
g = \dim H^1(C, \cO_C) ,
$$
and call it the genus of the curve $C$. If the curve $C$ is non-singular this number
$g$ is the usual genus of the curve $C$. Then, $U(\gX)$ is a point of $\UGM^{1-g}$.
A frame $\Xi$ is given by meromorphic functions $f_j$ whose Laurent expansion
is of the form
$$
f_{1-g -j+1/2} = \xi^{-n_j} + \sum_{n= -n_j+1}^\infty a_n^{(j)}\xi^n, \quad j=1,2, \ldots
$$
where we may choose $f_{1-g-1/2}=1$, hence $n_1=0$.
If our curve $C$ is non-singular or has only nodes, then by the Riemann-Roch theorem we
have that
$$
n_j= j+g-1 ,
$$
for $j \ge g+1$.
Put
$$
e(f_{1-g -j+1/2}) = e^{-n_j +1/2} + \sum_{n= -n_j+1}^\infty a_n^{(j)}e^{n+1/2}.
$$
For example we have $e(f_{1-g-1/2}) = e^{1/2}$.
Then,  the wedge product
\begin{equation}
\label{eq2.3}
\cdots \wedge e(f_{1-g -5/2}) \wedge  e(f_{1-g -3/2})\wedge  e(f_{1-g -1/2})
\end{equation}
gives an element of $\cFd(1-g)$ spanning $\det U(\gX)$ and defines the
 point ${\mathbf P}(\cFd(1-g))$ associated
to $U(\gX)$.  Since $H^0(C,\cO_C(*Q)))$ is isomorphic to $U(\gX)$ we may
define $\det H^0(C,\cO_C(*Q)) = \det U(\gX)$ and so we can regard the wedge product \eqref{eq2.3} as an
element of $ \det H^0(C, \cO_C(*Q))$.
Moreover, by \eqref{longex} we get a natural isomorphism.
$$
 \det H^0(C,\cO_C)    \otimes
(\det H^0(C, \cO_C(*Q))) ^{-1} \otimes
\det \bC[\xi^{-1}]\xi^{-1}
\otimes  (\det   H^1(C, \cO_C))^{-1} \cong {\mathbf C}.
$$
Thus we may regard $\det H^0(C, \cO_C(*Q))$ as
the determinant of the structure sheaf of the curve $C$, since there is a canonical
isomorphism $\det \bC[\xi^{-1}]\xi \cong {\mathbf C}$.
\end{Exmp}

The map
$$
\Phi : \UGM^p \rightarrow {\mathbf P}(\cFd(p)) ,
$$
which associates to each element $U \in \UGM^p$ the point
$\det U \in{\mathbf P}(\cFd(p))$ is called the Pl\"ucker embedding.
This is a generalization of
the usual Pl\"ucker embedding for the usual Grassmann manifold.
The following theorem is due to M. Sato and it plays an important role in
soliton theory.
\begin{Thm}
\label{thm2.1}
The Pl\"ucker embedding
$\Phi : \UGM^p \rightarrow {\mathbf P}(\cFd(p))$  is a holomorphic embedding.
The image is a closed submanifold of ${\mathbf P}(\cFd(p))$ which is defined by
Pl\"ucker's relations.
\end{Thm}

\section{Ghost Vacua}

Now we shall develop the theory of the ghost system for the case $j=0$ in section 1.
In this case the Virasoro algebra have the following commutation relations
\begin{equation}
\label{3.1}
\left[ L^{(0)}_n,\,  L^{(0)}_m \right]  =  (n-m)L_{n+m}^{(0)} -
\frac16 (n^3-n)\delta_{m+n.0}.
\end{equation}
Also we have the following important commutation relations
\begin{eqnarray}
\label{3.2}
\left[ L^{(0)}_n,\, \psi(z) \right] & = & z^{n+1}\frac{d}{dz}  \psi(z), \\
\label{3.3}
\left[ L^{(0)}_n,\, \overline{\psi}(z) \right] & = & z^n \left (z\frac{d}{dz}
+ (n+1)\right) \ovpsi(z).
\end{eqnarray}
These commutation relations suggest that the field operator $\psi(z)$ behaves
like a meromorphic function and $\ovpsi(z)$ behaves like a meromorphic
one-form. This fact will be used to define the ghost vacua of the $j=0$ ghost
system on a pointed curves with formal coordinates.

By an $N$-pointed curve
with formal coordinates  $\gX=(C; Q_1, \ldots,Q_N;\xi_1,\ldots,\xi_N)$  we mean that
the curve $C$ is reduced and projective but not necessarily connected and that
the points  $Q_j$  are  non-singular points of the curve $C$ and $\xi_j$ is a formal
coordinate of the curve $C$ with center $Q_j$. See \cite{U2} for
further details regarding curves with formal coordinates. We will
always assume that each of the connected components of $C$
contains at least one of the $Q_j's$.
Put
\begin{eqnarray*}
\cF_N &=& \bigoplus_{p_1, \ldots,p_N \in\bZ} \cF(p_1) \otimes \cdots \otimes \cF(p_N) ,\\
\cFd_N &=& \bigoplus_{p_1, \ldots, p_N \in \bZ}
{\cFd}(p_1)  \hat{\otimes}  \cdots \hat{\otimes} {\cFd}(p_N),
\end{eqnarray*}
where $\hat{\otimes}$ means the complete tensor product.
\begin{Dfn}
\label{dfn3.1}
The ghost vacua  $\cVd_{\ab}(\gX)$
of the spin $j=0$ ghost system is the linear subspace of
$\cFd_N$ consisting of elements $\langle \Phi|$ satisfying the following conditions:
\begin{enumerate}
\item  For all $|v\rangle \in \cF_N$, there exists a meromorphic function $f \in H^0(C,\cO_C(*\sumQj))$
such that
$\langle \Phi| \rho_j(\psi(\xi_j))|v\rangle$ is the Laurent expansion of $f$
 at the point $Q_j$ with respect to the formal coordinate $\xi_j$;
\item  For all $|v\rangle \in \cF_N$, there exists  a meromorphic one-form
$\omega \in H^0(C,\omega_C(*\sumQj))$ such that
$\langle \Phi| \rho_j(\ovpsi(\xi_j))|v\rangle d\xi_j$ is the Laurent expansion
of $\omega$
at the point $Q_j$ with respect to the coordinates $\xi_j$,
\end{enumerate}
where $\rho_j(A)$ means that the operator $A$ acts on the $j$-th component of
$\cF_N$  as
$$
\rho_j(A)|u_1 \otimes u_2 \otimes \cdots \otimes u_N\rangle =
(-1)^{p_1+\cdots+p_{j-1}} |u_1 \otimes \cdots \otimes u_{j-1}\otimes
Au_j \otimes u_{j+1} \otimes \cdots \otimes u_N\rangle.
$$
\end{Dfn}
We will reformulate the above two conditions into gauge conditions. For that purpose we
introduce the following notation.

For a meromorphic one-form $\omega \in H^0(C,\omega_C(*\sumQj))$ we let
$$
\omega_j = (\sum_{n=-n_0}^\infty a_n \xi_j^n )d\xi_j
$$
be the Laurent expansion at $Q_j$ with respect to the coordinate $\xi_j$.
Then, for the field operator $\psi(z)$ let us define $\psi[\omega_j]$ by
$$
\psi[\omega_j]= \Res_{\xi_j=0}(\psi(\xi_j)\omega_j) =
\sum_{n= -n_0}^\infty a_n\psi_{n+1/2}.
$$
Similarly we can define $\ovpsi[\omega_j]$.
For  a meromorphic function $f \in H^0(C,\cO_C(*\sumQj ))$ we let $f_j(\xi_j)$
be
the Laurent expansion of $f$ at $Q_j$ with respect to the coordinate $\xi_j$.
For the field  operator $\psi(z)$ define $\psi[f_j]$ by
$$
\psi[f_j] = \Res_{\xi_j=0}(\psi(\xi_j)f_j(\xi_j)d\xi_j)
$$
Put
\begin{eqnarray*}
\psi[\omega] = (\psi[\omega_1], \ldots, \psi[\omega_N]) , &&
\ovpsi[\omega] = (\ovpsi[\omega_1], \ldots, \ovpsi[\omega_N]) \\
\psi[f]= (\psi[f_1], \ldots, \psi[f_N]), &&
\ovpsi[f]= (\ovpsi[f_1], \ldots, \ovpsi[f_N]).
\end{eqnarray*}
Then, these operate  on $\cF_N$ from the left
and on $\cFd_N$ from the right. For example,
$\ovpsi[f]$ operates on $\cF_N$ from the left by
\begin{eqnarray*}
\ovpsi[f] |u_1\otimes \cdots \otimes u_N \rangle  & = &
\sum_{j=1}^N\rho_j(\ovpsi[f_j])|u_1\otimes \cdots \otimes u_N \rangle \\
&= & \sum_{j=1}^N(-1)^{p_1+\cdots+p_{j-1}}|u_1 \otimes \cdots \otimes u_{j-1}
\otimes \ovpsi[f_j]u_j \otimes u_{j+1}\otimes \cdots \otimes u_N\rangle
\end{eqnarray*}
for $|u_j\rangle \in \cFd(p_j)$ and  operates on $\cFd_N$ from the right by
\begin{eqnarray*}
\langle v_N \otimes \cdots v_1|\ovpsi[f]   & = &
\sum_{j=1}^N \langle v_N \otimes \cdots v_1|\rho_j(\ovpsi[f_j]) \\
&= & \sum_{j=1}^N(-1)^{p_1+\cdots+p_{j-1}} \langle v_N \otimes
\cdots \otimes v_{j+1}
\otimes v_j \ovpsi[f_j] \otimes v_{j-1}\otimes \cdots \otimes v_1 |
\end{eqnarray*}
for $\langle v_j| \in \cFd(p_j)$.
\begin{Thm}
\label{thm3.1}
The element $\langle \Phi| \in \cFd_N$ belongs to
the space of ghost vacua  $\cVd_{\ab}(\gX)$
of the $j=0$ ghost system if and only if $\langle \Phi|$ satisfies
the following two conditions.
\begin{enumerate}
\item $\langle\Phi|\psi[\omega] =0$ for any meromorphic one-form $\omega \in
H^0(C, \omega_C(*\sumQj))$.
\item  $\langle\Phi|\ovpsi[f] =0$
for any meromorphic function $f \in H^0(C,\cO_C(*\sumQj))$.
\end{enumerate}
\end{Thm}

The first (resp. second) condition in the above theorem is called the
first (resp. second) gauge condition. The first and second gauge conditions can be
rewritten in the following form:
\begin{enumerate}
\item  $\sum_{j=1}^N(-1)^{p_1+\cdots+p_{j-1}}
\langle \Phi|   u_1 \otimes \cdots
\otimes \cdots \otimes u_{j-1}
\otimes \psi[\omega_j]u_j \otimes u_{j+1}\otimes \cdots \otimes u_N\rangle =0$
for any $\omega \in H^0(C, \omega_C(*\sumQj))$ and $|u_j\rangle \in \cF(p_j)$,
$j=1,2,\ldots,N$.
\item  $\sum_{j=1}^N(-1)^{p_1+\cdots+p_{j-1}}
\langle \Phi|  u_1 \otimes \cdots \otimes u_{j-1}
\otimes \ovpsi[f_j]u_j \otimes
u_{j+1}\otimes \cdots \otimes u_N\rangle =0$
for any $f \in H^0(C, \cO_C(*\sumQj))$ and $|u_j\rangle \in \cF(p_j)$,
$j=1,2,\ldots,N$.
\end{enumerate}

It is easy to show that the ghost vacua  $\cVd_{\ab}(\gX)$ is a
finite dimensional vector space. More strongly we can prove the
following theorem.
\begin{Thm}
\label{thm3.2}
 For any $N$-pointed curve $\gX= (C;Q_1,\ldots,Q_N;\xi_1,\ldots
\xi_N)$ with formal coordinates we have
$$
\dim_\bC \cVd_{\ab}(\gX) =1.
$$
\end{Thm}

A proof is given in section 5.

In the later application we need to consider a disconnected curve.
The following proposition is an immediate consequence of the definition.
\begin{Prop}
\label{prop3.1}
Let
$$\gX_1=(C_1;Q_1, \ldots,Q_M;\xi_1, \ldots, \xi_M)$$
and
$$\gX_2=(C_2;Q_{M+1}, \ldots,Q_M;\xi_{M+1}, \ldots, \xi_N)$$
be pointed  curves with formal coordinates. Let $C$ be the disjoint union
$C_1 \sqcup C_2$ of the curves $C_1$, $C_2$.  Put
$$
\gX=(C; Q_1, \ldots, Q_N; \xi_1, \ldots,\xi_N).
$$
Then we have
$$
\cVd_{\ab}(\gX) = \cVd_{\ab}(\gX_1) \otimes
\cVd_{\ab}(\gX_2).
$$
\end{Prop}

Now we can introduce the dual  ghost vacua.
\begin{Dfn}
\label{dfn3.2}
Let $\cF_{\ab}(\gX)$ be the subspace of $\cF_N$ spanned by
$\psi[\omega] \cF_N$,
$\omega \in H^0(C,\omega_C(*\sumQj))$
and $\ovpsi[f] \cF_N$,  $f \in H^0(C,\cO(*\sumQj))$.
Put
$$
\cV_{\ab}(\gX) = \cF_N/\cF_{\ab}(\gX).
$$
The quotients space $\cV_{\ab}(\gX)$ is called
 the space of {\it dual ghost vacua} of the $j=0$ ghost system.
\end{Dfn}

Since
$\cVd_{\ab}(\gX)$ is finite dimensional,  $\cV_{\ab}(\gX)$ is
dual to $\cVd_{\ab}(\gX)$.

\begin{Exmp}
\label{exmp3.1}
Let us consider the one-dimensional complex projective space
${\mathbf P}^1=\bC \cup \{\infty\}$. Let $z$ be the coordinate of $\bC$.
We shall show that
$$
\dim_\bC \cVd_{\ab}(({\mathbf P}^1;0 ; z))=1.
$$
A basis of $H^0({\mathbf P}^1, \omega_{{\mathbf P}^1}(*0))$ is given by
$$
\frac{dz}{z^{m+2}}, \quad m=0,1,2, \ldots .
$$
A basis of $H^0({\mathbf P}^1, \cO_{{\mathbf P}^1}(*0))$ is given by
$$
\frac{1}{z^m}, \quad m=0,1,2,\ldots .
$$
Then, we have that
$$
\psi[\frac{dz}{z^{m+2}}] = \psi_{-m-3/2}, \quad \ovpsi[\frac{1}{z^m}] = \ovpsi_{-m+1/2},
\quad  m=0,1,2,\ldots .
$$
Hence, an element $\langle \Phi | \in \cVd_{\ab}(({\mathbf P}^1;0 ; z))$
satisfies the equations
\begin{eqnarray*}
\langle \Phi | \psi_{-m-3/2}&=&
 \langle \Phi | \wedge e^{-m-3/2}=0, \quad m=0,1,2,\ldots , \\
\langle \Phi | \ovpsi_{-m+1/2}&=&
 \langle \Phi |  i(e^{m-1/2})=0, \quad m=0,1,2,\ldots .
 \end{eqnarray*}
 Note that any element of $\cF$ is a finite linear combination of elements of
 the form
 \begin{eqnarray}
 \label{Basis2}
 &&\ovpsi_{\mu_1}\ovpsi_{\mu_2}\cdots\ovpsi_{\mu_r}
 \psi_{\nu_s}\psi_{\nu_{s-1}}\cdots\psi_{\nu_1}|0\rangle,  \\
 &&\mu_1<\mu_2<\cdots<\mu_r<0, \quad
 \nu_1<\nu_2<\cdots<\nu_s<0, \quad r\ge 0, \,s\ge 0 . \nonumber
 \end{eqnarray}
 By the first gauge condition, if $r>0$ we have
$$
 \langle \Phi | \ovpsi_{\mu_1}\ovpsi_{\mu_2}\cdots\ovpsi_{\mu_r}
 \psi_{\nu_s}\psi_{\nu_{s-1}}\cdots\psi_{\nu_1}|0\rangle=0.
$$
 Moreover, by the second gauge condition if $s>0$ and
 $\nu_s < -1/2$ we have
$$
 \langle \Phi | \psi_{\nu_s}\psi_{\nu_{s-1}}\cdots\psi_{\nu_1}|0\rangle=0.
$$
Thus we conclude that $\langle \Phi|$ is zero on all the element of the form
\eqref{Basis2}  except $\ovpsi_{-1/2}|0\rangle = |-1\rangle$. Therefore, $\langle \Phi |$
is a constant multiple of $\langle -  1|$ and we conclude
$$
\cVd_{\ab}(({\mathbf P}^1;0 ; z)) = \bC\langle -1|.
$$
\end{Exmp}

Let us assume that  $C$ is a non-singular curve of genus $g\ge 1$. Let us consider
a one-pointed curve $\gX = (C;Q;\xi)$ with a formal coordinate.  Choose a basis
$\{\omega_1,\ldots,\omega_g\}$ of holomorphic one-forms on $C$. We let
$$
\omega_i = \bigl(\sum_{n=0}^\infty a^{(i)}_n \xi^n\bigr)d\xi , \quad i=1,2, \ldots,g
$$
be the Taylor expansions of $\omega_i$'s at the point $Q$.
For any positive integer $j$ choose a meromorphic
one-form $\omega_{g+j} \in H^0(C,\omega_C(*Q))$ in such a way that
it has the Laurent expansion
$$
\omega_{g+j}= \bigl(\xi^{-(j+1)} +\sum_{n= -j}^\infty
a^{(g+j)}_n \xi^n  \bigr) d \xi
$$
at $Q$.
Put
\begin{eqnarray*}
e(\omega_i) &=& \sum_{n=0}^\infty a^{(i)}_n e^{n+1/2},  \\
e(\omega_{g+j}) &=& e^{-j-1/2} + \sum_{n= -j}^\infty a^{(g+j)}_n e^{n+1/2}.
\end{eqnarray*}
Then, the infinite sums $e(\omega_i)$ and $e(\omega_{g+j})$ are regarded as
elements of $\cWd$ (see section 2).  Put
\begin{equation}
\langle \omega(\gX) | = \cdots \wedge e(\omega_{g+2})\wedge e(\omega_{g+1})
\wedge e(\omega_g)\wedge \cdots \wedge e(\omega_1) .
\end{equation}
 The set
$$
\{\ldots, \omega_{g+2}, \omega_{g+1},\omega_g,
\ldots, \omega_1\}
$$
is a frame of $H^0(C,\omega_C(*Q)) \subset
\bC((\xi))d\xi$. Since $ H^0(C,\omega_C(*Q))  \in \UGM^{g-1}$, we
see that
$\langle \omega(\gX) | \in \cFd(g-1)$ is non-zero.
By a similar arguments as in Example \ref{exmp2.1} we may also regard $\langle \omega(\gX) | $
as an element of the determinant of the canonical sheaf
$\omega_C$. Note of course that $\langle \omega(\gX) | $ depends
on the choice of the basis $(\omega_i)$.
\begin{Lem}
\label{lem3.1}
The semi-infinite wedge product $\langle \omega(\gX) | $ is a non-zero
element of $\cFd(g-1)$ which satisfies the first and second gauge conditions.
Hence
$$
\langle \omega(\gX) |  \in \cVd_{\ab}(\gX).
$$
\end{Lem}

{\it Proof}. \quad For any element $\omega \in H^0(C,\omega_C(*Q))$ we let
$(\sum_{n=-n_0}^\infty a_n\xi^n)d\xi$ be its Laurent expansion. Then we have
$$
\psi[\omega]= \sum_{n=-n_0}^\infty a_n \psi_{n+1/2}
$$
whose right action on $\cFd$ is given by
$$
\wedge\left( \sum_{n=-n_0}^\infty a_n e^{n+1/2} \right) = \wedge e(\omega).
$$
Hence we have that
$$
\langle \omega(\gX)| \wedge  e(\omega)=0.
$$
Thus $\langle \omega(\gX)| $  satisfies the first gauge condition.

Now let $\sum_{m=-m_0}^\infty b_m\xi^m$ be the Laurent expansion of
a meromorphic function $f \in H^0(C,\cO_C(*Q))$ at $Q$. Then we have
$$
\ovpsi[f]= \sum_{m=-m_0}^\infty b_m \psi_{m+1/2}
$$
and its right action is given by
$$
\sum_{m=-m_0}^\infty b_m i(e^{-m-1/2}).
$$
Then we have that
\begin{eqnarray*}
e(\omega)\ovpsi[f] &=& \left( \sum_{n=-n_0}a_ne^{n+1/2} \right)
\left( \sum_{m=-m_0}^\infty b_m i(e^{-m-1/2} )\right) \\
&=& \sum_{m=-m_0}^\infty a_{-m-1} b_m  \\
&=& \Res_{\xi=0}(f(\xi)\omega) =0.
\end{eqnarray*}
Thus the second gauge condition is also satisfied. \QED
\begin{Cor}
\label{cor3.7}
Let $1 = l_1<l_2<\cdots <l_g\le 2g-1$ be the Weierstrass gap values of
the curve $C$ at the point $Q$. Then,
$$
\ovpsi[\xi^{-l_g}]\ovpsi[\xi^{-l_{g-1}}]\cdots\ovpsi[\xi^{-l_1}] |-1\rangle
= \oe_{l_g-1/2}\wedge\oe_{l_{g-1} -1/2}\wedge \cdots \wedge \oe_{l_1-1/2}
\wedge \oe_{-3/2}\wedge\oe_{-5/2} \wedge \cdots
$$
defines a non-zero element of $\cV_{\ab}(\gX)=\cF/\cF(\gX)$.
\end{Cor}

{\it Proof}. \quad The Weierstrass gap values have the property (for details
see Lemma \ref{lem6.2}):
$$
\langle e(\omega_g)\wedge\cdots\wedge e(\omega_1)|
\oe_{l_g-1/2}\wedge\oe_{l_{g-1} -1/2}\wedge \cdots \wedge\oe_{l_1-1/2}
\rangle \ne 0.
$$
On the other hand in the infinite wedge product
$$
\cdots \wedge e(\omega_{g+n} )\wedge \cdots \wedge e(\omega_{g+2} )
\wedge e(\omega_{g+1})
$$
the term $\langle -1|$ only appears when we choose the term $e^{-j-1/2}$ of
$e(\omega_{g+j})$ for each $j \ge 1$.  Other terms do not have a term $e^{-m - 1/2}$
for a certain positive integer $m$.
Since $\langle e(\omega_g)\wedge\cdots\wedge e(\omega_1)|$ does not contain
a term $e^{-m-1/2}$ for any $m \ge 1$ we have that
$$
\langle \omega(\gX)|
\ovpsi[\xi^{-l_g}]\ovpsi[\xi^{-l_{g-1}}]\cdots\ovpsi[\xi^{-l_1}] |-1\rangle
= \langle e(\omega_g)\wedge\cdots\wedge e(\omega_1)|
\oe_{l_g-1/2}\wedge\oe_{l_{g-1} -1/2}\wedge \cdots \wedge\oe_{l_1-1/2}
\rangle \ne 0.
$$
{} \QED

\begin{Thm}
\label{thm3.3}
The space of ghost vacua
$\cVd_{\ab}(\gX)$ is isomorphic to the determinant of the canonical bundle
$\omega_C$.
\end{Thm}

Notice that this theorem follows directly from Lemma \ref{lem3.1},
Theorem \ref{thm3.2}
and the discussion in section 2. The isomorphism of course depends
on the choice of a basis for $H^0(C,\omega_C(*Q))$.

Let $\gX=(C;Q_1,\ldots,Q_N;\xi_1,\ldots,\xi_N)$  be an $N$-pointed curve with
formal coordinates. Let $Q_{N+1}$ be a non-singular point and choose a
formal coordinate $\xi_{N+1}$ of $C$ with center $Q_{N+1}$.  Put
$$
\widetilde{\gX} = (C;Q_1,\ldots,Q_N,Q_{N+1} ;\xi_1,\ldots,\xi_N, \xi_{N+1}).
$$
Then the canonical linear mapping
\begin{eqnarray*}
\iota : \cF_N & \rightarrow & \cF_{N+1} \\
|v\rangle & \mapsto & |v\rangle \otimes |0\rangle
\end{eqnarray*}
induces  the canonical mapping
$$
\iota^* :\cFd_{N+1}  \rightarrow \cFd_N.
$$
\begin{Thm}
\label{thm3.4}
The canonical mapping $\iota^*$ induces
an isomorphism
$$
\cVd_{\ab}(\widetilde{\gX}) \cong \cVd_{\ab}(\gX) .
$$
\end{Thm}

This isomorphism is denoted the ``Propagation of vacua''
isomorphism.

{\it Proof}. \quad  Since we have injective maps
\begin{eqnarray*}
H^0(C,\omega_C(*\sumQj)) &\hookrightarrow &
H^0(C,\omega_C(*\sum_{j=1}^{N+1} Q_j)),\\
H^0(C,\cO_C(*\sumQj)) &\hookrightarrow &
H^0(C,\cO_C(*\sum_{j=1}^{N+1} Q_j)),
\end{eqnarray*}
and for $f(\xi) \in \bC[[\xi]]$ we have
$$
\psi[f(\xi)d\xi]|0\rangle =0 , \quad \ovpsi[f(\xi)| 0\rangle =0,
$$
the image $\iota^*(\cVd_{\ab}(\hat{{\gX}}))$ is contained in
$\cVd_{\ab}(\gX)$.

Therefore,  it is enough to show that any element $\langle \phi|
\in \cVd_{\ab}(\gX)$ uniquely determines an element $\langle\Phi|
\in \cVd_{\ab}(\hat{{\gX}})$  such that  $\iota^*(\langle \Phi|) =
\langle \phi |$.
Note that any element of $\cF$ is a finite linear combination of
elements of the form
 \begin{eqnarray*}
 &&\ovpsi_{\mu_1}\ovpsi_{\mu_2}\cdots\ovpsi_{\mu_r}
 \psi_{\nu_s}\psi_{\nu_{s-1}}\cdots\psi_{\nu_1}|0\rangle,  \\
 &&\mu_1<\mu_2<\cdots<\mu_r<0, \quad
 \nu_1<\nu_2<\cdots<\nu_s<0, \quad r\ge 0, \,s\ge 0 .
 \end{eqnarray*}
 By a double induction on $r$ and $s$ we shall show that $\langle \phi |$
 uniquely determines the value
$$
 \langle \Phi | \ovpsi_{\mu_1}\ovpsi_{\mu_2}\cdots\ovpsi_{\mu_r}
 \psi_{\nu_s}\psi_{\nu_{s-1}}\cdots\psi_{\nu_1}|0\rangle
$$
in such a way that $\langle \Phi |$ satisfies the first and second gauge conditions.

For $(r,s)=(0,0)$ put
$$
\langle \Phi | u \otimes 0\rangle = \langle \phi |  u \rangle
$$
for any $|u\rangle \in \cF_N$.

Choose $|u_j \rangle \in \cF(p_j)$ and put
$$
|u \rangle = |u_1\rangle \otimes \cdots \otimes |u_N\rangle.
$$
Choose  $\omega \in H^0(C, \omega_C(*\sumQj+(n+1)Q_{N+1}))$
in such a way that its Laurent expansion at the point $Q_{N+1}$ has the form
\begin{equation}
\label{Lafunc}
\omega_{N+1} =  \left(\xi^{-n-1}_{N+1} +
\sum_{n=0}^\infty a_n \xi^n_{N+1}\right) d \xi
\end{equation}
Then we have
$$
\psi[\omega_{N+1}] |0\rangle = \psi_{-n-1/2}|0\rangle =
i(\oe_{-n-1/2})|0\rangle.
$$
This is no-zero if and only if  $n$ is a non-negative integer.
For a non-negative integer $n$ define
\begin{equation}
\label{defpsi1}
\langle \Phi | u \otimes  \psi_{-n-1/2}|0\rangle
= (-1)^{p_1+\cdots p_N+1} \sum_{j=1}^N
\langle \Phi | \rho_j(\psi[\omega_j])u \otimes 0\rangle.
\end{equation}
This is independent of the choice of $\omega$ satisfying \eqref{Lafunc},
since if $\omega' \in H^0(C,\omega_C(*\sumQj+(n+1)Q_{N+1}))$ satisfies
\eqref{Lafunc} $\omega - \omega' $ is
holomorphic at $Q_{N+1}$.
Thus $\langle \Phi | u \otimes  \psi_{-n-1/2}|0\rangle $ is well-defined
for any integer $n$.  Then for any meromorphic form $\tau
\in H^0(C,\omega(*\sum_{n=1}^{N+1}Q_j))$ we have
\begin{equation}
\label{gauge0}
\sum_{j=1}^{N+1} \langle \Phi | \rho_j(\psi[\tau_j])u \otimes 0 \rangle =0.
\end{equation}
This establishes the first gauge condition on this subspace.

Next let us define
$$
\langle \Phi | u \otimes \psi_{-n_2-1/2} \psi_{-n_1-1/2}|0\rangle
$$
with $0\le n_2<n_1$.  Choose a meromorphic form
$\widetilde{\omega} \in H^0(C,\omega_C(*\sumQj+(n_2+1)Q_{N+1}))$ which has
the Laurent expansion
\begin{equation}
\label{Lofunc2}
\widetilde{\omega}_{N+1} =  \left(\xi^{-n_2-1}_{N+1} +
\sum_{n=0}^\infty b_n \xi^n_{N+1} \right) d\xi.
\end{equation}
Define
$$
\langle \Phi | u \otimes \psi_{-n_2-1/2} \psi_{-n_1-1/2}| 0 \rangle
= (-1)^{p_1+\cdots p_N+1} \sum_{j=1}^N
\langle \Phi | \rho_j(\psi[\widetilde{\omega}_j])u \otimes \psi_{-n_1-1/2} 0\rangle.
$$
We need to show that this is well-defined.
Choose another $\widetilde{\omega}'$ which has the  Laurent  expansion of the
same type \eqref{Lofunc2}.
Then $\tau=\widetilde{\omega} - \widetilde{\omega}'$ is holomorphic at $Q_{N+1}$, hence
$$
\psi[\tau_{N+1}]\psi_{-n_1-1/2}| 0 \rangle=0.
$$
Note that if $j>k$ we have that
\begin{eqnarray*}
&&\rho_j(\psi[\omega_j])\rho_k(\psi[\tau_k]) | u \rangle  \\
&& \phantom{XXX} = (-1)^{p_1+\cdots+p_{k-1}}\rho_j(\psi[\omega_j]) |u_1 \otimes \cdots
\otimes \psi[\tau_k]u_k \otimes \cdots u_N\rangle \\
&& \phantom{XXX} =
 (-1)^{(p_1+\cdots+p_{k-1} )+(p_1+\cdots + p_{k-1} +\cdots+p_{j-1}}
|u_1 \otimes \cdots \otimes \psi[\tau_k]u_k \otimes \cdots
\otimes \psi[\omega_j]u_j \otimes \cdots u_N\rangle \\
&& \phantom{XXX} = -\rho_k(\psi[\tau_k])\rho_j(\psi[\omega_j]) | u \rangle .
\end{eqnarray*}
We have the same result in the case $j \le k$.
Hence by \eqref{defpsi1} we have that
\begin{eqnarray*}
&&\sum_{j=1}^N
\langle \Phi | \rho_j(\psi[\tau_j])u \otimes \psi_{-n_1-1/2} 0\rangle \\
&& \phantom{X}=(-1)^{p_1+\cdots p_N+1} \sum_{k=1}^N  \sum_{j=1}^N
\langle \Phi | \rho_k(\psi[\omega_k])\rho_j(\psi[\tau_j]) u \otimes 0\rangle \\
&&\phantom{X}=(-1)^{p_1+\cdots p_N} \sum_{k=1}^N
\left\{  \sum_{j=1}^N \langle \Phi | \rho_j(\psi[\tau_j])
\rho_k(\psi[\omega_k]) u \otimes 0\rangle \right\}\\
&&\phantom{X} =(-1)^{p_1+\cdots p_N} \sum_{k=1}^N
\left\{  \sum_{j=1}^{N+1} \langle \Phi | \rho_j(\psi[\tau_j])
\rho_k(\psi[\omega_k]) u \otimes 0\rangle  \right\} \\
&&\phantom{X} = 0,
\end{eqnarray*}
where  the last equality is a consequence of \eqref{gauge0}.
Now by the same argument as above we can show that
for any $\tau \in H^0(C, \omega_C(*\sum_{j=1}^{N+1}))$ and any non-negative
integer $n$  the first gauge condition holds on the subspace spanned by
$|u \otimes \psi_{-n-1/2} 0\rangle$:
\begin{equation}
\label{gauge1}
\sum_{j=1}^{N+1} \langle \Phi | \rho_j(\psi[\tau_j])
 u \otimes  \psi_{-n-1/2} 0\rangle =0 .
\end{equation}
In this way, by induction on $r$ we can show that for any
non-positive integers $n_1>n_2>\cdots > n_r\ge 0$
$$
\langle \Phi |u \otimes \psi_{-n_1-1/2}\cdots \psi_{-n_r-1/2}0 \rangle
$$
is well-defined and we have the first gauge condition on the subspace spanned by
$|u \otimes \psi_{-n_1-1/2}\cdots \psi_{-n_r-1/2}0 \rangle$:
\begin{equation}
\label{gauger}
\sum_{j=1}^{N+1} \langle \Phi | \rho_j(\psi[\tau_j])
 u \otimes  \psi_{-n_1-1/2}\cdots \psi_{-n_{r-1}-1/2} 0\rangle =0.
\end{equation}

Next let us define
$$
\langle \Phi | u \otimes  \ovpsi_{-n-1/2}
\psi_{-n_1-1/2}\cdots \psi_{-n_r-1/2} 0\rangle.
$$
Choose  $f \in H^0(C,\cO_C(*\sumQj+(n+1)Q_{N+1}))$
in such a way that it has the Laurent expansion
\begin{equation}
\label{Laomega}
f_{N+1} =  \xi^{-n-1}_{N+1} +
\sum_{n=0}^\infty a_n \xi^n_{N+1}
\end{equation}
at the point $Q_{N+1}$. Then we have
\begin{eqnarray*}
\ovpsi[f_{N+1}]
\psi_{-n_1-1/2}\cdots \psi_{-n_{r-1}-1/2} | 0 \rangle
&= &\ovpsi_{-n-1/2} \psi_{-n_1-1/2}\cdots \psi_{-n_{r-1}-1/2} | 0\rangle \\
&= &
\oe_{n+1/2} \wedge |\psi_{-n_1-1/2}\cdots \psi_{-n_{r-1}-1/2} | 0\rangle.
\end{eqnarray*}
This is no-zero if  $n$ is a non-negative integer.
For non-negative integer $n$ define
\begin{eqnarray*}
&&\langle \Phi | u \otimes  \ovpsi_{-n-1/2}
\psi_{-n_1-1/2}\cdots \psi_{-n_{r-1}-1/2} 0\rangle \\
&& \phantom{XX}= (-1)^{p_1+\cdots p_N+1} \sum_{j=1}^N
\langle \Phi | \rho_j(\ovpsi[f_j])u \otimes
\psi_{-n_1-1/2}\cdots \psi_{-n_{r-1}-1/2} |0\rangle.
\end{eqnarray*}
This is independent of the choice of $f$ satisfying \eqref{Laomega},
since if $f'$ satisfies \eqref{Laomega} $f - f' $ is
holomorphic at $Q_{N+1}$.
Thus
$$
\langle \Phi | u \otimes  \ovpsi_{-n-1/2}
\psi_{-n_1-1/2}\cdots \psi_{-n_{r-1}-1/2} |0\rangle
$$
is well-defined for any integer $n$.  Then for any meromorphic function $h
\in H^0(C,\cO(*\sum_{n=1}^{N+1}Q_j))$ we have that
\begin{equation}
\label{gauge20}
\sum_{j=1}^{N+1} \langle \Phi | \rho_j(\ovpsi[h_j])u \otimes
\psi_{-n_1-1/2}\cdots \psi_{-n_{r-1}-1/2} |0\rangle =0,
\end{equation}
establishing the second gauge condition on this subspace.

Now by induction on $s$ we can show that
$$
\langle \Phi | u \otimes  \ovpsi_{-n_s-1/2} \cdots \ovpsi_{-n_1-1/2}
\psi_{-n_1-1/2}\cdots \psi_{-n_{r-1}-1/2} |0\rangle
$$
are well-defined and they satisfy the second gauge condition:
\begin{equation}
\label{gauge2s}
\sum_{j=1}^{N+1} \langle \Phi | \rho_j(\ovpsi[f_j])u \otimes
\ovpsi_{-n_{s-1}-1/2} \cdots \ovpsi_{-n_1-1/2}
\psi_{-n_1-1/2}\cdots \psi_{-n_{r-1}-1/2} |0\rangle =0.
\end{equation}
Thus $\langle\Phi| \in \cVd_{\ab}(\widetilde{\gX})$ has been uniquely constructed
from $\langle \phi| \in \cVd_{\ab}(\gX)$. By our construction we
have that
$$
\iota^*(\langle\Phi|) = \langle \phi|.
$$
{} \QED

Next let us consider a curve $C$ with a node $P$. Let $\widetilde{C}$ be
the curve obtained by resolving the singularity at $P$ and let
$\pi : \widetilde{C} \rightarrow C$ be the natural holomorphic mapping.
Then $\pi^{-1}(P)$ consists of two points $P_+$ and $P_-$.
Let
$$
\gX=(C;Q_1,\ldots,Q_N;\xi_1,\ldots,\xi_N)
$$
be an
$N$-pointed curve with formal coordinates and we let
$$
\widetilde{\gX} = (\widetilde{C};P_+,P_-,Q_1,\ldots,Q_N;z,w, \xi_1,\ldots,\xi_N)
$$
be the associated $N+2$-pointed curve with formal coordinates.
Define  an element $|0_{+,-}\rangle \in \cF\otimes \cF$ by
\begin{equation}
\label{0+-}
|0_{+,-}\rangle = |0\rangle \otimes |-1\rangle
 - |-1\rangle \otimes |0 \rangle .
\end{equation}
The natural  inclusion
\begin{eqnarray*}
 \cF_N & \hookrightarrow & \cF_{N+2} \\
 |u \rangle & \mapsto & |0_{+,-}\rangle  \otimes |u \rangle
\end{eqnarray*}
defines a natural linear mapping
$$
\iota_{+,-}^* : \cFd_{N+2}  \rightarrow   \cFd_N .
$$
\begin{Thm}
\label{thm3.5}
The natural mapping $\iota^*_{+,-}$  induces  a natural
isomorphism
$$
\cVd_{\ab}(\widetilde{\gX}) \cong \cVd_{\ab}(\gX).
$$
\end{Thm}

{\it Proof}. \quad First let us show that
$$
\iota_{+,-}^*(\cVd_{\ab}(\widetilde{\gX}) ) \subset \cVd_{\ab}(\gX).
$$
For $\langle \widetilde{\Phi}| \in \cVd_{\ab}(\widetilde{\gX}) $ put
$\langle \Phi|= \iota_{+,-}^*(\langle \widetilde{\Phi}| ) $.
For $|u \rangle \in \cF_N$ we have that
$$
\langle \Phi|u \rangle = \langle \widetilde{\Phi}| 0_{+,-}\otimes u \rangle.
$$
For a meromorphic one-form $\omega \in H^0(C,\omega_C(*\sumQj))$
put $\widetilde{\omega}= \pi^*(\omega)$. Then $\widetilde{\omega}$ has
a Laurent expansions at $P_+$ and $P_-$ of the following form:
\begin{eqnarray*}
\widetilde{\omega}_+ &=& \bigl( \frac{a_{-1}}{z} + a_0 + a_1z+a_2z^2 +\cdots
\bigr)dz \\
\widetilde{\omega}_- &=& \bigl( -\frac{a_{-1}}{w} + b_0 + b_1w+b_2w^2 +\cdots
\bigr)dw
\end{eqnarray*}
Then we have that
\begin{eqnarray*}
\rho_+(\psi[\widetilde{\omega}_+])|0_{+,-}\rangle &=&
\psi[\widetilde{\omega}_+]|0\rangle \otimes |-1\rangle -
\psi[\widetilde{\omega}_+]|-1\rangle \otimes |0 \rangle  \\
&=& a_{-1}|-1\rangle \otimes |-1\rangle
\end{eqnarray*}
Similarly we have that
\begin{eqnarray*}
\rho_-(\psi[\widetilde{\omega}_-])|0_{+,-}\rangle  &=&
|0\rangle \otimes \psi[\widetilde{\omega}_-]  |-1\rangle +
|-1\rangle \otimes  \psi[\widetilde{\omega}_-] |0 \rangle \\
&=& - a_{-1}|-1\rangle \otimes |-1\rangle.
\end{eqnarray*}
Since $\omega$ and $\widetilde{\omega}$ have the same Laurent
expansions at $Q_j$,  by the above results  we have that
\begin{eqnarray*}
&& \sum_{j=1}^N \langle \Phi |\rho_j(\psi[\widetilde{\omega}_j])u \rangle   \\
&&= \sum_{j=1}^N \langle \widetilde{\Phi} |
0_{+,-}\otimes  \rho_j(\psi[\widetilde{\omega}_j]) u \rangle \\
&& = \langle \widetilde{\Phi} |\rho_+(\psi[\widetilde{\omega}_+])0_{+,-}\otimes u \rangle
+ \langle \widetilde{\Phi} |\rho_-(\psi[\widetilde{\omega}_+])0_{+,-}\otimes u \rangle 
+ \sum_{j=1}^N \langle \widetilde{\Phi} |
0_{+,-}\otimes  \rho_j(\psi[\widetilde{\omega}_j]) u \rangle \\
&&= 0
\end{eqnarray*}
by the first gauge condition for $\langle \widetilde{\Phi} |$.
Thus $\langle\Phi|$ satisfies the first gauge condition. Similarly
we can show that $\langle\Phi|$ satisfies the second gauge condition.
Hence $\langle\Phi| \in \cVd_{\ab}(\gX)$.

Now let us show that $\iota_{+,-}^*$ is bijective. For that purpose it is
enough to show that $\langle \Phi | \in \cVd_{\ab}(\gX)$ determines
uniquely $\langle \widetilde{\Phi} | \in \cVd_{\ab}(\widetilde{\gX})$.
Choose $f \in H^0(\widetilde{C},\cO_{\widetilde{C}}(*\sumQj))$ such that
the Taylor expansions at $P_{\pm}$ have the forms
\begin{eqnarray}
\label{Tayl1}
f_+ &=& -1  + \sum_{n=1}^\infty a_nz^n , \\
\label{Tayl2}
f_- &=& \sum_{n=1}^\infty b_nw^n .
\end{eqnarray}
Then we have that
\begin{eqnarray*}
\ovpsi[f_+] &=& -\ovpsi_{1/2} +\sum_{n=1}^\infty a_n\ovpsi_{n+1/2},  \\
\ovpsi[f_-]&=& \sum_{n=1}^\infty b_n\ovpsi_{n+1/2} .
\end{eqnarray*}
Hence we have that
\begin{eqnarray*}
\rho_+(\ovpsi[f_+]) |0_{+,-}\rangle &=&
 \ovpsi[f_+]|0\rangle \otimes |-1 \rangle
- \ovpsi[f_+]|-1\rangle\otimes|0\rangle
=  |0\rangle\otimes |0\rangle, \\
\rho_-(\ovpsi[f_-]) |0_{+,-}\rangle &=&
 |0\rangle \otimes\ovpsi[f_-] |-1\rangle
+  |-1\rangle\otimes \ovpsi[f_-] |0\rangle
= 0.
\end{eqnarray*}
Now define
\begin{equation}
\label{def00}
\langle \widetilde{\Phi} | 0\otimes 0 \otimes u \rangle =
\sum_{j=1}^N \langle \widetilde{\Phi} | 0_{+,-} \otimes \rho_j(\ovpsi[f_j]) u \rangle
= \sum_{j=1}^N \langle \Phi |  \rho_j(\ovpsi[f_j] )u \rangle.
\end{equation}
This is well-defined,  since if $f' \in
H^0(\widetilde{C},\cO_{\widetilde{C}}(*\sumQj))$  has the same type of Taylor
expansions  \eqref{Tayl1} and \eqref{Tayl2}, then $h=f-f'$
is an element of $H^0(C,\cO_C(*\sumQj))$.
Hence
$$
\sum_{j=1}^N \langle \Phi |  \rho_j(\ovpsi[h_j] )u \rangle =0 .
$$
Now let us show that  $\langle \widetilde{\Phi} |$ just defined satisfies
the first gauge condition for any element
$h \in H^0(\widetilde{C}, \cO_{\widetilde{C}}(*\sumQj))$ and the second gauge
condition for any $\omega \in H^0(\widetilde{C}, \omega_{\widetilde{C}}(*\sumQj))$ applied
to the vectors
$|0 \otimes 0 \otimes u\rangle$.
Let
\begin{eqnarray*}
\omega_+ &=& \bigl(a_0+a_1z +\cdots \bigr)dz \\
\omega_- &=& \bigl(b_0+b_1w +\cdots \bigr)dw
\end{eqnarray*}
be the Taylor expansions of $\omega$ at $P_\pm$. Then
we have that
$$
\psi[\omega_+]|0 \rangle \otimes |0\rangle =0,
\quad \psi[\omega_-]|0 \rangle \otimes |0\rangle =0.
$$
Since  we may regard $\omega$ as an element of $H^0(C,\omega_C(*\sumQj))$
having a zero at $P$, we get by using \eqref{def00} that
\begin{eqnarray*}
&& \langle \widetilde{\Phi} |\rho_+(\psi[\omega_+] )0\otimes 0 \otimes u \rangle
 + \langle \widetilde{\Phi} |\rho_-(\psi[\omega_-] )0\otimes 0 \otimes u \rangle
+ \sum_{j=1}^N
 \langle \widetilde{\Phi} |\rho_j(\psi[\omega_j]) 0\otimes 0 \otimes u \rangle\\
&& \phantom{XX} =  - \sum_{k=1}^N \Big\{ \sum_{j=1}^N
 \langle \Phi|\rho_j(\psi[\omega_j] )\rho_k[f_k]u \rangle \Big\} =0,
\end{eqnarray*}
 by the first gauge condition for $\langle \Phi |$.

Next let
\begin{eqnarray*}
h_+ &=& a_0+a_1z +a_2z^2+\cdots \\
h_-&=&  b_0+b_1w+b_2w^2+\cdots
\end{eqnarray*}
be the Taylor expansion of $h \in H^0(\widetilde{C}, \cO_{\widetilde{C}}(*\sumQj))$
at $P_\pm$.
Then we have that
$$
\rho_+(\ovpsi[h_+])|0\rangle \otimes |0\rangle =0, \quad
\rho_-(\ovpsi[h_-])|0\rangle \otimes |0\rangle =0.
$$
If $a_0=b_0$ then $h \in H^0(C,\cO_C(*\sumQj))$.
Then by similar arguments as above, we can show that
$$
\langle \widetilde{\Phi} |\rho_+(\ovpsi[h_+] )0\otimes 0 \otimes u \rangle
 + \langle \widetilde{\Phi} |\rho_-(\ovpsi[h_-] )0\otimes 0 \otimes u \rangle
+ \sum_{j=1}^N
 \langle \widetilde{\Phi} |\rho_j(\ovpsi[h_j]) 0\otimes 0 \otimes u \rangle=0.
$$
If $a_0\ne b_0$, then by subtracting from $f$ the constant $b_0$ and multiplying by a
constant, we may assume that $h$ has Taylor expansions of type
\eqref{Tayl1} and \eqref{Tayl2}.  Then by the above argument we
have that
$$
\langle \widetilde{\Phi} |0\otimes 0 \otimes u \rangle =
\sum_{k=1}^N \langle \Phi | \rho_k(\ovpsi[h_k]) u \rangle.
$$
Then we have that
\begin{eqnarray*}
&& \sum_{j=1}^N \langle \widetilde{\Phi} |0\otimes 0 \otimes
\rho_j(\ovpsi[h_j])u \rangle =\sum_{j=1}^N
\sum_{k=1}^N \langle \Phi | \rho_j(\ovpsi[h_j])\rho_k(\ovpsi[h_k]) u \rangle\\
&&\phantom{X} = -\sum_{k=1}^N
\sum_{j=1}^N \langle \Phi | \rho_j(\ovpsi[h_j])\rho_k(\ovpsi[h_k]) u \rangle
= -\sum_{j=1}^N \langle \widetilde{\Phi} |0\otimes 0 \otimes
\rho_j(\ovpsi[h_j])u \rangle ,
\end{eqnarray*}
since we have that
$$
\rho_j(\ovpsi[h_j])\rho_k(\ovpsi[h_k]) | u \rangle =
- \rho_k(\ovpsi[h_k])\rho_j(\ovpsi[h_j]) \ |u \rangle.
$$
Hence
$$
\sum_{j=1}^N \langle \widetilde{\Phi} |0\otimes 0 \otimes
\rho_j(\ovpsi[h_j])u \rangle = 0.
$$
Thus we conclude that
$$
\langle \widetilde{\Phi} |\rho_+(\ovpsi[h_+] )0\otimes 0 \otimes u \rangle
 + \langle \widetilde{\Phi} |\rho_-(\ovpsi[h_-])0\otimes 0 \otimes u \rangle
+ \sum_{j=1}^N
 \langle \widetilde{\Phi} |\rho_j(\ovpsi[h_j]) 0\otimes 0 \otimes u \rangle=0.
$$

Now we can apply Theorem \ref{thm3.4} twice and obtained the desired
result. \QED

\begin{Rmk}
\label{rmk3.1}
The choice of $|0_{+,-}\rangle$ is not unique. Any non-zero
multiple of $|0_{+,-}\rangle$ gives a natural isomorphism of
the above theorem. But the above normalization of $|0_{+,-}\rangle$
will be seen to be
compatible with the preferred elements
 of the  ghost vacua on a nodal  curve and on its normalization. For details see Theorem
 \ref{thm6.5} below.
\end{Rmk}

\section{Sheaf of   Ghost Vacua}

  The arguments of this section are almost identical to those in [U2],  section 4,
but for the readers' convenience we shall repeat them here in this abelian case.

Let
$$
\gF = (\pi : \cC \rightarrow \cB; s_1, \ldots, s_N; \xi_1, \ldots, \xi_N)
$$
be a family of $N$-pointed semi-stable curves with formal coordinates. That is
$\cC$ and $\cB$ are  complex manifolds, $\pi$ is a proper holomorphic
mapping, and for each point
$b \in \cB$,  $\gF(b)= (C_b =\pi^{-1}(b); s_1(b),
\ldots, s_N(b); \xi_1, \ldots,\xi_N)$
is an $N$-pointed semi-stable curve with formal coordinates.
We let $\Sigma$
be the locus of double points of the fibers of $\gF$ and let $D$ be
$\pi(\Sigma)$. Note that $\Sigma$ is a non-singular submanifold of
codimension two in $\cC$,  and $D$ is a divisor
in $\cB$ whose irreducible components
$D_i$, $i = 1, 2, \dots ,m'$ are non-singular.

In this section  we use the following notation freely.
$$
S_j = s_j(\cB), \quad   S = \sum_{j=1}^N S_j.
$$
Put
$$
\cF_N(\cB)  = \cF_N \otimes_\bC \cO_\cB, \quad
\cFd_N(\cB)= \cO_\cB  \otimes_\bC    \cFd_N.
$$
\begin{Dfn}
\label{dfn4.1}
We define the subsheaf $\cVd_{\ab}(\gF)$ of $\cFd_N(\cB)$ by the gauge conditions:
\begin{eqnarray*}
&&\sum_{j=1}^N \langle \Phi | \psi[\omega_j] =  0 ,
\quad \hbox{\rm for all $\omega  \in \pi_*(\omega_{\cC/\cB}(*S))$}, \\
&&\sum_{j=1}^N \langle \Phi | \ovpsi[f_j] =  0,
\quad \hbox{\rm for all $f  \in \pi_*\cO_\cC(*S)$} .
\end{eqnarray*}
where $\omega_j$ and  $f_j$ are the Laurent expansion of
$\omega$ and $f$  along $S_j$ with respect to
the coordinate $\xi_j$.

The sheaf $\cVd_{\ab}(\gF)$ is called the sheaf of  ($j=0$) ghost vacua
or the sheaf of abelian vacua of the family $\gF$. Similarly the
sheaf $\cV_{\ab}(\gF)$ of  ($j=0$)  dual
ghost vacua of the family is defined by
$$
\cV_{\ab}(\gF)= \cF_N(\cB)/\cF_{\ab}(\gF).
$$
where  $\cF_{\ab}(\gF)$ is the $\cO_\cB$-submodule of $\cF_N(\cB)$ given by
$\cF_{\ab}(\gF) = \cF^0_{\ab}(\gF) + \cF^1_{\ab}(\gF)$,
where $\cF^0_{\ab}(\gF)$ is the span of $\ovpsi[f]\cF_N(\cB)$ for
all $f \in \pi_*\cO_\cC(*S)$ and $\cF^1_{\ab}(\gF)$ is the span of
$\psi[\omega]\cF_N(\cB)$ for all
$\omega \in \pi_*\omega_{\cC/\cB}(*S)$.
\end{Dfn}

Note that we have
\begin{equation}
\label{4.a}
\cVd_{\ab}(\gF) = \underline{\Hom}_{\cO_\cB}(\cV_{\ab}(\gF),
\cO_\cB).
\end{equation}
Moreover, by the right exactness of the tensor product we have
that
\begin{equation}
\label{4.1}
 \cV_{\ab}(\gF)\otimes_{\cO_\cB} \cO_{\cB, b}/\frak{m}_b
 \cong \cV_{\ab}(\gF(b)).
\end{equation}

\begin{Prop}
\label{prop4.1}
The sheaves $\cV_{\ab}(\gF)$ and $\cVd_{\ab}(\gF)$ are coherent
$\cO_\cB$-modules.
\end{Prop}

{\it Proof}\quad First we need to introduce a filtration on
$\cF$. Note that any element of $\cF$ is a finite linear combination of
elements of the form
$$
|u \rangle =\ovpsi_{\mu_1}\ovpsi_{\mu_2}\cdots\ovpsi_{\mu_r}
 \psi_{\nu_s}\psi_{\nu_{s-1}}\cdots\psi_{\nu_1}|0\rangle
$$
with $\mu_1<\mu_2<\cdots <\mu_r<0$,
$\nu_1<\nu_2<\cdots<\nu_s<0$. Let us define an alternative degree
$d(|u\rangle )$ by
$$
d(|u\rangle) = \sum_{i=1}^r (-\mu_i+\frac12) +
\sum_{j=1}^s(-\nu_j - \frac12).
$$
Note that this alternative degree $d(|u\rangle)$ is different form the degree
of the Maya diagram associated to $|u\rangle$. We shall only need
this alternative degree in this proof and here we
shall just refer to it as the degree.
Let $\cF_d$ be the subspace of $\cF$ spanned by  the elements
of degree $d$. The filtration
$F_\bullet$ on $\cF$ is defined by
$$
F_p \cF = \bigoplus_{d \le p} \cF_d .
$$
Note that $F_p\cF=0$ for $p<0$ and $F_0\cF= \bC|0\rangle +
\bC|-1\rangle$.  By definition we have that
\begin{equation}
\label{hoi4.1.1}
\psi_\nu F_p\cF \subset  F_{p-\nu-1/2} \cF, \quad
\ovpsi_\mu F_p\cF  \subset  F_{p-\mu+1/2} \cF .
\end{equation}
Moreover we have
$$
Gr_\bullet^F \cF = \bigoplus \cF_d.
$$

Also let us introduce filtrations on $\bC((\xi))$ and $\bC(\xi))d\xi$ by
$$
F_p \bC((\xi))= \bC[[\xi]]\xi^{-p}, \quad
F_p \bC(\xi))d\xi = \bC[[\xi]]\xi^{-p} d\xi.
$$
Then we have that
$$
Gr_\bullet^F\bC((\xi)) = \bC[\xi,\xi^{-1}],
\quad
Gr_\bullet^F\bC((\xi))d\xi  = \bC[\xi,\xi^{-1}] d\xi .
$$
Put
\begin{eqnarray*}
\psi[F_p\bC((\xi))d\xi ] &=& \{ \, \ovpsi[\omega(\xi)]\,|\,
\omega(\xi) \in F_p \bC((\xi))d\xi \,\} \\
\ovpsi[F_p\bC((\xi))] &=& \{ \, \ovpsi[f(\xi)]\,|\, f(\xi) \in F_p \bC((\xi))\,\}
\end{eqnarray*}
By \eqref{hoi4.1.1} it is easy to show that following facts
$$
\psi[F_p\bC((\xi))d\xi ] F_q\cF \subset F_{p+q}\cF, \quad
\ovpsi[F_p\bC((\xi))]  F_q\cF \subset F_{p+q}\cF .
$$

let $\gX=(C, Q, \xi)$ be a one-pointed curve of genus $g$ with coordinate.
The filtrations $F_\bullet$ on $\bC((\xi))$ and
$\bC((\xi))d\xi$ induce the ones on
$H^0(C, \cO_C(*Q))$ and $H^0(C, \omega_C(*Q))$ and we have that
$$
F_p H^0(C, \cO_C(*Q)) = H^0(C, \cO_C(pQ)), \quad
F_p H^0(C, \omega_C(*Q)) = H^0(C, \omega_C(pQ)) .
$$
Thus the filtrations are compatible with the actions of $\psi[\omega]$ and
$\ovpsi[f]$ to $\cF$ for  $\omega \in H^0(C, \omega_C(*Q))$ and $f \in
H^0(C, \cO_C(*Q))$. Hence the inclusion  $\cF(\gX) \subset \cF$ is
compatible with the filtrations and the quotient space $\cV(\gX)$ has
the induced filtration $F_\bullet$.  Thus we conclude that
$$
Gr_\bullet^F \cV(\gX) = Gr_\bullet^F \cF / Gr_\bullet^F \cF(\gX).
$$

Let $1=n_1<n_2<\cdots<n_g=2g-1$ be the Weierstrass gap values
 of the curve $C$ at $Q$.
Then $Gr_\bullet^F H^0(C, \cO_C(*Q)) \subset \bC[\xi^{-1}]$
does not contain $\xi^{-n_i}$, $i=1, \ldots , g$ but contain
$\xi^{-n}$, $n \ne n_i$,  $i=1, \ldots , g$.  Also
$Gr_\bullet^F H^0(C, \omega_C(*Q))$ contains all $\xi^{-n}d\xi$,
$n \le 0$. Thus
$$
|u \rangle =\ovpsi_{\mu_1}\ovpsi_{\mu_2}\cdots\ovpsi_{\mu_r}
 \psi_{\nu_s}\psi_{\nu_{s-1}}\cdots\psi_{\nu_1}|0\rangle
$$
may not be in $Gr_\bullet^F \cV(\gX)$ if $s=0$ and
$\{\mu_1, \ldots, \mu_r\} \subset \{ -n_2+1/2, \ldots, -n_g+1/2\}$.
Thus $Gr_\bullet^F\cV(\gX)$ is of finite dimensional. Moreover if
$V$ is the subspace spanned by
$$\ovpsi_{-m_1+1/2}\ovpsi_{-m_2+1/2}\cdots\ovpsi_{-m_r+1/2}|0\rangle
$$
with $2g-1\ge m_1> m_2>\cdots>m_r \ge0$ where $1 \le r \le 2g$, then
$V$ is a finite dimensional vector space and the
natural linear map $V \rightarrow Gr_\bullet^F\cV(\gX)$ is surjective.

Now for a family
$$
\gF = (\pi : \cC \rightarrow \cB; s ; \xi)
$$
of one-pointed curves of genus $g$ with coordinate,
we can introduce filtrations on $\cF(\cB)$, $\cF(\gF)$ , $\pi_*\cO_\cC(*S)$
and $\pi_*\omega_{\cC/\cB}(*S)$ of $\cO_\cB$-submodules
in  the same way as above so that
we can introduce a filtration on $\cV(\gF)$. Moreover, there is
an $\cO_\cB$-module homomorphism $V \otimes_\bC\cO_\cB
\rightarrow Gr_\bullet^F\cV(\gF)$ which is surjective by virtue of \eqref{4.1} and
the above argument. Thus $Gr_\bullet^F\cV(\gF)$ is a coherent $\cO_\cB$-module.
Therefore, $\cV(\gF)$ is also a coherent $\cO_\cB$-module. This proves the
proposition for the one-pointed case. The general case can be proved similarly or
we can use Theorem \ref{thm3.4} to reduce the general case to the one-pointed case.
\QED

Let us now show the local freeness of the sheaves $\cV_{\ab}(\gF)$ and
$\cVd_{\ab}(\gF)$.  For that purpose we first introduce a certain
$\cO_{\cB}$-submodule $\mathcal L(\gF)$ of
$$
\bigoplus_{j=1}^N  {\cO}_{\cB}((\xi_j^{-1}))
\frac{d}{d\xi_j}
$$
and an action of
$\cL(\gF)$
 on the sheaves $\cV_{\ab}(\gF)$ and $\cVd_{\ab}(\gF)$ as
first order twisted differential operators.
This action  will also be used to define a connection on the
sheaves $\cV_{\ab}(\gF)$ and $\cVd_{\ab}(\gF)$.

First recall that we have an exact sequence
$$
0 \rightarrow
 \Theta_{\cC/\cB}(-S)
\rightarrow \Theta_{\cC/\cB}( mS)
\rightarrow     \bigoplus_{j=1}^N \bigoplus_{k=0}^m
{\cO}_{\cB}\xi_j^{-k}
{\frac{d}{d\xi_j}}
 \rightarrow 0,
$$
for any positive integer $m$. From this exact sequence we obtain
the exact sequence
\begin{equation}
\label{4.2.1}
 0 \rightarrow \pi_*(\Theta_{\cC/\cB}( mS))
\overset{b_m}{\longrightarrow} \bigoplus_{j=1}^N \bigoplus_{k=0}^m
{\cO}_{\cB}\xi_j^{-k}
{\frac{d}{d\xi_j}}
\overset{\vartheta_m}{\longrightarrow} R^1\pi_* \Theta_{\cC/\cB}(-S)
 \rightarrow 0.
\end{equation}
Hence, we have
the following exact sequence of
${\cO}_{\cB}$-modules
\begin{equation}
\label{4.2.2}
 0 \rightarrow \pi_*(\Theta_{\cC/\cB}( * S))
\overset{b}{\longrightarrow} \bigoplus_{j=1}^N {\cO}_{\cB}[
\xi_j^{-1} ] {\frac{d}{d\xi_j}}
\overset\vartheta{\longrightarrow} R^1\pi_* \Theta_{\cC/\cB}(-S)
 \rightarrow 0.
\end{equation}
Note that the mappings  $b$ and $b_m$ correspond to the Laurent expansions
with respect to $\xi_j$ up to zero'th order. To define the first order
 differential operators acting
on  the sheaves $\cV_{\ab}(\gF)$ and $\cVd_{\ab}(\gF)$, we need to rewrite
the exact sequence \eqref{4.2.2} in the following way.

There is an exact sequence
$$
0 \rightarrow \Theta_{\cC/\cB} \rightarrow \Theta_{\cC}
\overset{d\pi}{\longrightarrow }
 \pi^*\Theta_{\cB}     \rightarrow  0,
$$
where $\Theta_{\cC/\cB}$ is the sheaf of
vector fields tangent to the fibers of $\pi$. Put
$$
\Theta'_{\cC,\pi} = d\pi^{-1}(\pi^{-1}\Theta_{\cB}(-\log D)).
$$
Hence, $\Theta'_{\cC,\pi}$ is the sheaf of
vector field on $\cC$  tangent to $\Sigma$
whose horizontal components are constant along
the fibers of $\pi$. That is, $\Theta'_{\cC,\pi}$
consists of germs of holomorphic vector fields of the form
$$
 a(z,u) \frac{\partial}{\partial z} +
 \sum_{i=1}^{m'} b_i(u) u_i \frac{\partial}{\partial u_i}
 + \sum_{i=m'+1}^n b_i(u)  \frac{\partial}{\partial u_i},
 $$
  where $(z,u_1 ,\dots, u_n)$ is a system of local
coordinates such that the mapping $\pi$ is expressed as the projection
  $$
 \pi (z,u_1 ,\dots, u_n) =(u_1, \dots, u_n)
 $$
and $\pi(\Sigma) = D$ is given by the equation
$$
u_1 \cdot u_2 \cdots u_{m'} = 0  .
$$
  More generally, we can define the sheaf $\Theta'_{\cC}(mS)_\pi$
as the one consisting of germs of meromorphic vector fields of
the form
$$
 A(z,u) \frac{\partial}{\partial z} +
 \sum_{i=1}^{m'} B_i(u) u_i  \frac{\partial}{\partial u_i}
 + \sum_{i=m'+1}^n B_i(u)   \frac{\partial}{\partial u_i},
 $$
  where $A(z,u)$ has a poles of order at most $m$
 along $S$. Now we have an exact sequence
\begin{equation}
\label{4.2.3}
 0 \rightarrow \Theta_{\cC/\cB}(mS) \rightarrow \Theta'_{\cC}(mS)_\pi
 \overset{d\pi}{\longrightarrow }
 \pi^{-1}\Theta_{\cB}( - \log D)    \rightarrow  0 .
\end{equation}
Note that $\Theta'_{\cC}(mS)_\pi $ has the structure of  a
sheaf of Lie algebras by the usual bracket operation of  vector
fields and the above exact sequence is one of sheaves of Lie algebras.

For $m > \frac1N(2g - 2)$, by \eqref{4.2.3} we have an exact sequence of
$\cO_\cB$-modules.
\begin{equation}
\label{4.2.4}
 0 \rightarrow \pi_* \Theta_{\cC/\cB}(mS)  \longrightarrow
\pi_*  \Theta'_{\cC}(mS)_\pi  \overset{d\pi}{\longrightarrow }
\Theta_{\cB}( - \log D)    \rightarrow  0 ,
\end{equation}
which is also an exact sequence of sheaves of Lie algebras.
Taking $m \rightarrow \infty$ we obtain the exact sequence
\begin{equation}
\label{4.2.5}
 0 \rightarrow \pi_* \Theta_{\cC/\cB}(*S)  \longrightarrow
\pi_*  \Theta'_{\cC}(*S)_\pi  \overset{d\pi}{\longrightarrow }
\Theta_{\cB}( - \log D)    \rightarrow  0 .
\end{equation}

 The exact sequences \eqref{4.2.2} and \eqref{4.2.5} are related by the
following commutative diagram.
$$\begin{array}{ccccccccc}
 0 & \rightarrow &  \pi_* \Theta_{\cC/\cB}(*S)  & \rightarrow  &
\pi_*  \Theta'_{\cC}(*S)_\pi  & \overset{d\pi}{\longrightarrow} &
\Theta_{\cB}( - \log D)  & \rightarrow & 0 \\
&&&&&&&&\\
&& || &&\phantom{\, p}\downarrow \,p &&
 \phantom{\, \rho}\downarrow \,\rho && \\
&&&&&&&&\\
0 &\rightarrow & \pi_* \Theta_{\cC/\cB}(*S)  &\rightarrow &
\bigoplus_{j=1}^N \cO_\cB[ \xi_j^{-1}] \frac d{d \xi_j}
&\overset{\vartheta}{\longrightarrow} &
R^1\pi_* \Theta_{\cC/\cB}( -S)  &\rightarrow &  0,
\end{array}
$$
where $\rho$ is the Kodaira-Spencer mapping
of the family $\gF$ and $p$ is
given by taking the non-positive degree  part of
the $\frac d {d\xi_j}$ part of the Laurent expansions
of the vector fields in $\pi_*\Theta'_\cC(*S)_\pi$
at $s_j(\cB)$.

In the following, for simplicity we assume that
the Kodaira-Spencer mapping $\rho$
is an injective  homomorphism of $\cO_\cB$-modules.
Therefore, $p$ is an injection.  Let
$$
\widetilde{p }  :   \pi_*  \Theta'_{\cC}(*S)_\pi  \longrightarrow
\bigoplus_{j=1}^N \cO_\cB(( \xi_j))\frac d{d \xi_j}
$$
be defined by taking the $\dfrac d{d \xi_j}$ part of the
Laurent expansions at $s_j(\cB)$. Since $p$ is
injective, $\widetilde{p}$ is injective. Put
\begin{equation}
\label{lF}
\cL(\gF) := \widetilde p ( \pi_*  \Theta'_{\cC}(*S)_\pi).
\end{equation}
Then, we have the exact sequence
\begin{equation}
\label{4.2.6}
0 \rightarrow \pi_* \Theta_{\cC/\cB}(*S)  \longrightarrow
\cL(\gF)  \overset{\theta}{\longrightarrow }
\Theta_{\cB}( - \log D)    \rightarrow  0
\end{equation}
of $\cO_\cB$-modules.  This is the exact sequence corresponding
to \eqref{4.2.2}.

The Lie bracket  $[\phantom{e}, \phantom{e}]_d$ on
$\cL(\gF)$ is obtained from the bracket on
$\pi_*  \Theta'_{\cC}(*S)_\pi$  by the mapping $\widetilde{p}$.
Thus,  for $\vec \ell$, $\vec m \in \cL(\gF)$ we have that
\begin{equation}
\label{4.2.7}
[\vec \ell, \vec m]_d = [\vec \ell, \vec m]_0 +
\theta(\vec \ell)(\vec m) - \theta(\vec m)(\vec \ell)
\end{equation}
where $[\cdot,\cdot]_0$ is the usual bracket of
formal vector fields and
the action of $\theta(\vec \ell)$ on
$$
\vec{m} = \left(m_1\frac{d}{d\xi_1}, \dots,
m_N \frac{d}{d\xi_N} \right)
$$
is defined by
$$
  \left( \theta(\vec \ell)(m_1) \frac{d}{d\xi_1}, \dots,
           \theta(\vec \ell)(m_N)\frac{d}{d\xi_N} \right).
$$
Then, the exact sequence \eqref{4.2.6} is also an exact sequence of
sheaves of Lie algebras.

 Let us define an action of $\cL(\gF)$ on
$\cF_N(\cB)$.
\begin{Dfn}
\label{dfn4.2}
{\rm
For ${\vec \ell} = ({\ell}_1, \dots,
{\ell}_N) \in \cL(\gF)$, the action of $D(\vec \ell)$ on
$\cF_N(\cB)$ is defined by
\begin{equation}
\label{4.2.8}
     D(\vec \ell) (F \otimes | u  \rangle)
    = \theta(\vec \ell)(F) \otimes | u \rangle -
F \cdot \Bigl(\sum_{j=1}^N \rho_j(T[{\ell}_j] )\Bigr) | u \rangle,
\end{equation}
where
$$
    F \in {\cO}_{\cB}, \quad | u \rangle \in \cF_N,
$$
and
$$
T[\ell] = \Res_{z=0}(T(z) \ell(z) dz).
$$}
\end{Dfn}

The following proposition can be proved in the same manner as
Proposition 4.2.2 of [U2].
\begin{Prop}
\label{prop4.2}
The action $D(\vec \ell)$ of $\vec \ell \in
\cL(\gF)$ on $\cF_N(\cB)$ defined above has the following
properties.
\begin{enumerate}
\item For any $F  \in {\cO}_{\cB}$ we have
$$
 D(F \vec \ell) = F  D(\vec \ell).
$$
\item For $\vec \ell$, $\vec  m \in  \cL(\gF)$ we have
$$
[ D(\vec \ell), D(\vec m) ] = D([ \vec \ell, \vec m ]_d)
 -  \frac{1}{ 6} \sum_{j=1}^N \Res_{\xi_j=0}
\left(  \frac{d^3 \ell_j}{d \xi_j^3} m_j d\xi_j \right)\cdot {\text id}.
$$
\item For $F  \in {\cO}_{\cB}$ and $|\phi\rangle \in
\cF_N(\cB)$ we have
$$
  D(\vec \ell)(F |\phi\rangle) = (\theta(\vec \ell)(F ))|\phi\rangle
 + F  D(\vec \ell)|\phi \rangle.
$$
Namely, $D(\vec \ell)$ is a first order differential operator, if
$\theta(\vec \ell ) \neq 0$.
\end{enumerate}
\end{Prop}
We define the dual action of $\cL(\gF)$ on $\cFd_N(\cB)$
by
\begin{equation}
\label{4.2.9}
  D(\vec \ell) (F \otimes \langle \Phi|) =
   \theta(\vec \ell)(F) \otimes \langle \Phi|  +
   \sum_{j=1}^N F  \cdot \langle \Phi| \rho_j(T[{\ell}_j]).
\end{equation}
where
$$
 F \in \cO_{\cB}, \quad \langle \Phi | \in \cFd_N(\cB).
$$
  Then, for any $| u \rangle \in \cF_N(\cB)$ and
$\langle\Phi | \in \cFd_N(\cB)$, we have
\begin{equation}
\label{4.2.10}
 \{D(\vec \ell) \langle \Phi| \}|\widetilde{\Phi}\rangle +
 \langle \Phi|\{D(\vec \ell)|\widetilde{\Phi}\rangle\}
=  \theta(\vec \ell)\langle \Phi|\widetilde{\Phi}\rangle .
\end{equation}
This agrees with the usual definition of the dual connection.

Now we shall show that the  operator $D(\vec \ell)$
acts on $\cV_{\ab}(\gF)$.
\begin{Prop}
\label{prop4.3}
For any $\vec \ell \in \cL(\gF)$
we have
$$
D(\vec \ell)(\cF_{\ab}(\gF))  \subset
   \cF_{\ab}(\gF).
$$
Hence, $D(\vec \ell)$ operates on $\cV_{\ab}(\gF)$. Moreover, it is a
first order differential operator, if $\theta(\vec \ell) \neq 0$.
\end{Prop}

{\it Proof}. \quad  An element of $\cF_{\ab}(\gF)$ is an
$\cO_\cB$-linear combination of elements of the form
$$
F \otimes \Bigl(\sum_{j=1}^N \rho_j(\psi[\omega_j])|u \rangle\Bigr), \quad
F' \otimes \Bigl(\sum_{j=1}^N \rho_j(\ovpsi[f_j])|u \rangle\Bigr)
$$
where
$$
  F , F' \in \cO_{\cB}, \quad \omega \in \pi_*(\omega_{\cC/\cB}(*S)),
  \quad f \in \pi_*\cO_\cC(*S),
\quad |u \rangle \in \cF_N
$$
and  $f_j$ is the Laurent expansion of $f$ along $S_j = s_j(\cB)$
with respect to the  coordinate $\xi_j$ and likewise for $\omega$.
First we shall show the following equality as  operators on $\cF_N(\cB)$
\begin{equation}
\label{4.2.11}
  \Bigl[ D(\vec \ell), \sum_{j=1}^N \rho_j(\ovpsi[f_j]) \Bigr]
 = \sum_{j=1}^N \rho_j\left(\ovpsi\bigl[\theta(\vec \ell)(f_j)
   +  \underline{\ell}_j(f_j)\bigr]\right),
\end{equation}
where $\theta(\vec{\ell})$
operates on the coefficients of $f_j$.
By Proposition \ref{prop4.2},~(3) it is enough to show the
equality \eqref{4.2.11} as operators on $\cF_N$.
For $| u \rangle \in \cF_N$, by \eqref{4.2.8}
we have
\begin{eqnarray*}
&& D(\vec \ell)  \Bigl(\sum_{j=1}^N  \rho_j(\ovpsi[f_j])| u \rangle \Bigr) -
    \sum_{j=1}^N \rho_j(\ovpsi[f_j]) (D(\vec \ell)| u \rangle) \\
&&\qquad= \sum_{j=1}^N \{\rho_j(\ovpsi[ \theta(\vec{\ell})(f_j)])
-  \rho_j(T[\underline \ell_j])\rho_j(\ovpsi[f_j])\}| u \rangle
 + \sum_{j=1}^N
\rho_j(\ovpsi[f_j])\rho_j(T[\underline{\ell}_j])| u \rangle \\
&&\qquad = \sum_{j=1}^N \{\rho_j\Bigl(\ovpsi\bigl[\theta(\vec \ell)(f_j) +
\underline{\ell}_j(f_j))\}\bigr] \Bigr) | u \rangle .
\end{eqnarray*}
This implies equation \eqref{4.2.11}. Now $\theta(\vec \ell)(f_j )+ \underline{\ell}_j(f_j)$ is
nothing but the Laurent expansion at $s_j(\cB)$ of
the meromorphic function $\tau(h)$ where $\tau = \widetilde{p}^{-1}({\vec \ell})
\in  \pi_*(\Theta'_\cC(*S)_\pi)$.
 Hence, we have the result that $D(\vec \ell) \cF^0_{(0,1}(\gF) \subset \cF^0_{(0,1}(\gF)$.

 Next let us consider $\omega \in \pi_*(\omega_{\cC/\cB}(*S))$.
For $\tau = \widetilde{p}^{-1}({\vec \ell})
\in  \pi_*(\Theta'_\cC(*S)_\pi)$ put
$$
 \widetilde{\omega} =\frac{d}{dt} \Bigl(\exp\{t \tau \} ^*(\omega)\Bigr).
$$
Then the Laurent expansion of $\widetilde{\omega}$ along $S_j$ is
written as
$$
\widetilde{\omega}_j = {\ell}_j(A_j(\xi_j)) +  \omega^{(\ell_j)},
$$
where $\omega_j = A_j(\xi_j)d\xi_j$ and
$$
\omega^{(\ell_j)} = A_j(\xi_j)\frac{d \ell_j(\xi_j)}{d\xi_j} .
$$
Now we are ready to prove the following equality as  operators on $\cF_N(\cB)$
\begin{equation}
\label{4.2.11a}
  \Bigl[ D(\vec \ell), \sum_{j=1}^N \rho_j(\psi[\omega_j]) \Bigr]
 = - \sum_{j=1}^N \rho_j\left(\psi\bigl[\widetilde{\omega}_j \bigr]\right).
\end{equation}
Note that $\widetilde{\omega} \in \pi_*(\omega_{\cC/\cB}(*S))$. Hence, the equality
\eqref{4.2.11a} implies the desired result.  To prove \eqref{4.2.11a} it is
enough to show that for
\begin{eqnarray*}
\underline{\ell} &=& \sum_{n=-n_0}^\infty l_nz^n \frac{d}{dz} \\
\omega &=& A(z)dz = (\sum_{m=-m_0}^\infty a_mz^m)dz
\end{eqnarray*}
we have that
$$
\Bigl[ T[\underline{\ell}], \psi[\omega] \Bigr] = -\psi[\{\ell(z)A'(z) +A(z)\ell'(z)\}dz].
$$
This can be proved as follows. First note that by \eqref{3.2}
 we have that
$$
\Bigl[ L_n, \psi_{k+1/2}\Bigr] = -(n+k+1)\psi_{n+k+1/2}.
$$
Hence
\begin{eqnarray*}
\Bigl[ T[\underline{\ell}], \psi[\omega] \Bigr]
&=& \sum_{n,k} a_k  l_{n+1}  \Bigl[ L_n, \psi_{k+1/2}\Bigr] \\
&=& - \sum_{n, k} (n+k+1)a_k l _{n+1} \psi_{n+k+1/2} \\
&=& -\psi[\{\ell(z)A'(z) +A(z)\ell'(z)\}dz].
\end{eqnarray*}
{}From this we conclude that $D(\vec \ell)\cF^1_{\ab}(\gF) \subset
\cF^1_{\ab}(\gF)$.
\QED

\begin{Prop}
\label{prop4.4}  For each element
$\vec \ell \in
\cL(\gF)$, $D(\vec \ell)$ acts on $\cVd_{\ab}(\gF)$.  Moreover, if
$\theta(\vec \ell) \neq 0$, then $D(\vec \ell)$ acts on  $\cVd_{\ab}(\gF)$
as a twisted first order
differential operator.
\end{Prop}

{\it Proof}. \quad  Choose
$$
\langle \widetilde{\Psi} | \in \cVd_{\ab}(\gF).
$$
For any element  $f \in \pi_*(\cO_{\cC}(*S))$ and
$|  u \rangle \in \cF_N(\cB)$, by Proposition \ref{prop4.3} and
\eqref{4.2.10} we have that
$$
\{ D(\vec \ell) \langle \widetilde{\Psi } | \}(\ovpsi[ f]
|  u \rangle  =
  \theta(\vec \ell)(\langle \widetilde{\Psi } |
\ovpsi[f] | u \rangle) -
  \langle  \widetilde{\Psi }  |\{D(\vec \ell)(\ovpsi[ f]
 | u \rangle)\}  = 0.
$$
Thus we conclude that
$$
    D(\vec \ell) \langle \widetilde{\Psi}  | \in \cVd_{\ab}(\gF).
$$
The remaining statement is an easy consequence of definition~\eqref{4.2.9}.
{}\QED

\begin{Cor}
\label{cor4.1}
The $\cO_{\cB}$-module $\cVd_{\ab}(\gF)$ is locally free on $\cB \setminus D$.
\end{Cor}

For a proof see [U2], Proposition~4.2.4.

For a coherent $\cO_{\cB}$-module $\mathcal G$, the locus $M$ consisting of points at
which
$\mathcal G$ is not locally free, is a closed analytic subset of $\cB$ of
codimension
at least~$2$.  Therefore, we have the following corollary.

\begin{Cor}
\label{cor4.2}
Let $W$ be the maximal subset of $\cB$
over which $\cV_{\ab} (\gF)$ is not locally free. Then, $W$ is an
analytic subset of $\cB$ and
$$
    W \subsetneqq D.
$$
\end{Cor}

  Since  we defined
$$
\cVd_{\ab}(\gF) = \underline{\Hom}_{\cO_{\cB}}(\cV_{\ab}(\gF), \cO_{\cB}),
$$
we have the following corollary.
\begin{Cor}
\label{cor4.3}
$\cVd_{\ab}(\gF)|_{\cB \setminus D}$ is a
locally free $\cO_{\cB}$-module
and for any subvariety $Y$ of $\cB \setminus D$  we have
an $\cO_Y$-module isomorphism
$$
   \cO_Y   \otimes_{\cO_{\cB}}   \cVd_{\ab}(\gF)   \simeq \cVd_{\ab}(\gF|_Y).
$$
\end{Cor}

These two corollaries play a crucial role in proving locally freeness in
general. The above corollaries imply the following theorem.

\begin{Thm}
\label{thm4.1}
If $\gF$ is a family of $N$-pointed smooth curves with formal coordinates,
then  $\cV_{\ab}(\gF)$ and $\cVd_{\ab}(\gF)$ are
locally free $\cO_{\cB}$-modules and they are dual to each other.
\end{Thm}

 Note that for the bilinear pairing
$$
\langle \cdot|\cdot\rangle
 :  \cV_{\ab}(\gF) \times \cVd_{\ab} (\gF) \rightarrow \cO_{\cB},
$$
we have the equality
\begin{equation}
\label{4.2.14}
 \{ D(\vec{\ell}) \langle\Psi | \}| u \rangle +
  \langle \Psi |\{D(\vec \ell)| u \rangle  \}= \theta (\vec{\ell})
 (\langle\Psi | u \rangle).
\end{equation}

Now let us define a connection of $\cVd_{\ab}(\gF)$.
First note that by direct calculations we have that
\begin{eqnarray*}
&& \langle \Phi |\rho_j\Bigl(\frac{d\psi(\zeta_j)}{d\zeta_j}
\ovpsi(\xi_j)\Bigr)|u\rangle \\
&&\phantom{XX} =- \sum_{m=-\infty}^\infty
\sum_{n=1}^\infty(m+1)\langle \Phi |\rho_j(\psi_{m+1/2}
\ovpsi_{n-1/2}) | u \rangle \zeta_j^{-m-2}\xi_j^{-n} \\
&& \phantom{XXXXXX}+ \sum_{m=-\infty}^\infty
 \sum_{n=0}^\infty(m+1)\langle \Phi |\rho_j(\ovpsi_{-n-1/2}
\psi_{m+1/2}) | u \rangle \zeta_j^{-m-2}\xi_j^{n}
  - \frac{\langle \Phi|u\rangle}{(\zeta_j -\xi_j)^2} .
\end{eqnarray*}

This implies
\begin{eqnarray}
&& \langle \Phi |\rho_j(T(\xi_j))| u \rangle (d\xi_j)^2 \nonumber \\
&& \label{emj}
\phantom{XX} =
\lim_{\zeta_j \rightarrow \xi_j}\left\{
 \langle \Phi |\rho_j \Bigl(\frac{d \psi(\zeta_j)}{d\zeta_j}
\ovpsi(\xi_j)\Bigr)|u\rangle  d\zeta_jd\xi_j
 +\frac{\langle \Phi|u\rangle}{(\zeta_j -\xi_j)^2}
  d\zeta_jd\xi_j \right\} .
\end{eqnarray}
Note that we can as in [U2, 3.4] define the
correlation functions
$$
 \langle \Phi |\frac{d \psi(w)}{dw}
\ovpsi(z)|u\rangle  dwdz, \quad  \langle \Phi |T(z) |u\rangle  dz^2
$$
such that
$$
 \langle \Phi |T(z)|u\rangle  dz^2  =
\lim_{w \rightarrow z}\left\{
 \langle \Phi |\frac{d \psi(w)}{dw}
\ovpsi(z)|u\rangle  dwdz
 +\frac{\langle \Phi|u\rangle}{(w -z)^2}
  dwdz \right\} .
$$

Note that $\langle \Phi |\frac{d \psi(w)}{dw}
\ovpsi(z)|u\rangle dwdz$ is a meromorphic section of
$p_1^*\omega_{\cC/\cB} \otimes_{\cO_{\cC\times_{\cB}}}
p_2^*\omega_{\cC/\cB}$
 over $\cB$ whose Laurent expansion in a neighbourhood of
$S_j \times S_j$ is equal to the above expansion of $\langle \Phi
|\rho_j\Bigl(\frac{d\psi(\zeta_j)}{d\zeta_j}
\ovpsi(\xi_j)\Bigr)|u\rangle$, where $p_i: \cC \times_\cB \cC
\rightarrow \cC$ is the projection to the $i$-th factor.

Now choose a meromorphic bidifferential
$$
\omega \in H^0(\cC\times_{\cB}\cC,
\omega_{\cC\times_{\cB}\cC}(2\Delta))
$$
which has the form
$$
\omega(w,z) = \frac{dw dz}{(w-z)^2} + \hbox{\rm holomorphic }
$$
in a neighbourhood of the diagonal $\Delta$ of $\cC\times_{\cB}\cC$.
Define $\widetilde{T}(z)$ by
\begin{equation}
\label{tildeem}
\langle \Phi |\rho_j(\widetilde{T}(z))| u \rangle (dz)^2
= \lim_{w \rightarrow z}\left\{
 \langle \Phi |\frac{d \psi(w)}{dw}
\ovpsi(z)|u\rangle  dwdz
 +\langle \Phi|u\rangle\omega(w,z) \right\} .
\end{equation}
Then $\langle \Phi |\rho_j(\widetilde{T}(z))| u \rangle (dz)^2$ is a global meromorphic
form, that is a meromorphic section of $\omega_{\cC/\cB}$ over $\cB$.
The {\it projective connection\/}  $S_\omega(z)(dz)^2$ associated to
the bidifferential $\omega$ is defined by
\begin{equation}
\label{projective}
S_\omega(z)(dz)^2 = 6 \lim_{w \rightarrow z}
\Big\{ \omega(w,z) - \frac{dwdz}{(w - z)^2} \Bigr\}.
\end{equation}
Then,  by \eqref{tildeem}  we have that
\begin{equation}
\label{diffem}
\langle \Phi |\rho_j(T(z))| u \rangle (dz)^2 =
\langle \Phi |\widetilde{T}(z))| u \rangle (dz)^2
-\frac16\langle \Phi| u \rangle S_\omega(z)^2(dz)^2 .
\end{equation}
Then for any element $| u \rangle \in \cF_N$ by \eqref{4.2.8}
\begin{eqnarray*}
\langle \Phi |\{D(\vec{\ell})| u \rangle \}
 &=&
  - \sum_{j=1}^N \Res_{\xi_j=0}
\Big(\ell_j(\xi_j)\langle \Phi |\rho_j(\widetilde{T}(\xi_j))| u \rangle d\xi_j\Bigr)
\nonumber \\
&&
\phantom{XX}+ \frac16\langle \Phi | u \rangle  \sum_{j=1}^N \Res_{\xi_j=0}
\Big(\ell_j(\xi_j) S_\omega(\xi_j)d\xi_j\Bigr).
\end{eqnarray*}
If $\theta(\vec{\ell})=0$, that is $\vec{\ell}$ is the image of
a global vector field $\tau \in \pi_*(\Theta_{\cC/\cB}(\*S))$, then
$$
\ell_j(\xi_j)\langle \Phi |\rho_j(\widetilde{T}(\xi_j))| u \rangle d\xi_j
$$
is the Laurent expansion of a global meromorphic relative one-form which has poles
only at $S_j$. Hence the first term of the right hand side of the above
equality is zero. Therefore, in this case we have that
$$
\langle \Phi |\{D(\vec{\ell})| u \rangle \} =
 \frac16\langle \Phi | u \rangle  \sum_{j=1}^N \Res_{\xi_j=0}
\Big(\ell_j(\xi_j) S_\omega(\xi_j)d\xi_j\Bigr).
$$
Now put
\begin{equation}
\label{aomega}
b_\omega(\vec{\ell}) =\sum_{j=1}^N \Res_{\xi_j=0}
\Big(\ell_j(\xi_j) S_\omega(\xi_j)d\xi_j\Bigr).
\end{equation}
Then this defines an $\cO_\cB$-module homomorphism
$$
b_\omega : \cL(\gF) \rightarrow \cO_\cB,
$$
and if $\theta(\vec{\ell})=0$ then we have that
$$
\langle \Phi |\{D(\vec{\ell})| u \rangle \} =\frac16 b_\omega(\vec{\ell})
\langle \Phi | u \rangle .
$$
For a vector field $X \in \Theta(-\log D)$ on $\cB$ tangent to $D$
choose $\vec{\ell} \in \cL(\gF)$ such that $\theta(\vec{\ell})=X$.
Then the connection on $\cV_{\ab}(\gF)$ is defined by
\begin{equation}
\label{4.18}
\nabla_X^{(\omega)}([| \Psi\rangle]) = D(\vec{\ell})([| \Psi\rangle])
-\frac16  b_\omega(\vec{\ell}) ([| \Psi\rangle]),
\end{equation}
for $| \Psi\rangle \in \cV_{\ab}(\gF)$. Dually the connection
on $\cVd_{\ab}(\gF)$ is defined by
\begin{equation}
\label{4.19}
\nabla_X^{(\omega)}(\langle\Phi|) =  D(\vec{\ell})(\langle\Phi|)
+ \frac16 b_\omega(\vec{\ell})    (\langle\Phi|),
\end{equation}
for $\langle\Phi| \in \cVd_{\ab}(\gF)$.
These are well-defined, that is the right hand sides of \eqref{4.18}
and \eqref{4.19} are independent of the choice of $\vec{\ell} \in\cL(\gF)$
with $\theta(\vec{\ell})=X$.
Just like for the non-abelian conformal field theory (see for example [U2], section 5)
we can prove the following theorem.
\begin{Thm}
\label{thm4.2}
The operator $\nabla^{(\omega)}$ defines a
projectively flat holomorphic connection of
the sheaves $\cV_{\ab}(\gF)$ and $\cVd_{\ab}(\gF)$.
Moreover, the connection has a regular singularity along the locus
$D \subset \cB$ which is the locus of the singular curves.
The connection $\nabla^{(\omega)}$ depends on the choice of
bidifferential $\omega$ and if we choose another bidifferential
$\omega'$ then there exists a holomorphic one-form $\phi_{\omega,\omega'}$
on $\cB$ such that
\begin{equation}
\label{4.20}
\nabla_X^{(\omega)} - \nabla_X^{(\omega')}= \frac16\langle
\phi_{\omega,\omega'}, \, X\rangle.
\end{equation}
Moreover, the curvature form $R$ is given by
\begin{equation}
\label{4.21}
R(X,Y) = \frac16 \Big\{b_\omega(\vec{n}) - X(b_\omega(\vec{m}) )
+Y(b_\omega(\vec{\ell}) )- \sum_{j=1}^N
\Res_{\xi_j=0}\big(\frac{d^3\ell_j}{d\xi_j}m_jd \xi_j\bigr)\Bigr\},
\end{equation}
where $X$,$Y \in \Theta_{\cC/\cB}(*S))$,
$\vec{\ell}$, $\vec{m} \in \cL(\gF)$ with $X=\theta(\vec{\ell})$,
$Y= \theta(\vec{m})$, and $\vec{n} \in \cL(\gF)$ is defined by
$\vec{n} = [\vec{\ell}, \vec{m}]_d$  $($see \eqref{4.2.7}$)$.
\end{Thm}
\begin{Rmk}
\label{rmk4.1}
For the connection on the
sheaf of non-abelian  vacua $\cVd_{\vec{\lambda}}(\gF)$
of level $l$ over $\cB$ with gauge symmetry $\frak{g}$,
a complex simple Lie algebra,
the curvature form $R^\frak{g}$ is given by
$$ R^\frak{g}(X,Y) = \frac{c_v}{2}R(X,Y)\cdot \hbox{\rm id} $$
where
$$ c_v = \frac{l\cdot \dim \frak{g}}{l + g^*}, \quad
\hbox{\rm $g^*$ is the dual Coxeter number of $\frak{g}$.} $$
For
the details see [U2], section 5.
\end{Rmk}

\section{Degeneration of Curves and  Sewing}

{}From a curve with one ordinary double point,
we construct via the usual sewing a one-parameter family of curves over
a disk, whose fibers over the punctured disk are smooth curves.
We shall show that the sheaf of $j=0$  ghost vacua gives an invertible sheaf
over the whole disk.
This implies that the dimension of the space of ghost vacua is one for any pointed curve.
The arguments of this section are almost the same, as those from
[U2] section 5, we just need to modify a few arguments in our setting.

First we recall certain basic facts about the energy-momentum tensor
and the fermion operators.
In the $j=0$ ghost system,  the energy-momentum tensor $T(z)$ is
given  by
$$
T^{(0)}(z)= \normalord \frac{d\psi(z)}{dz} \ovpsi(z) \normalord.
$$
Note that we have the formula
\begin{equation}
L_n = \sum_{k=0}^\infty k \ovpsi_ {n-k+1/2}\psi_{k-1/2} +
\sum_{k=1}^\infty k \psi_{-k-1/2}\ovpsi_{n+k+1/2}.
\end{equation}
If $n \ne 0$ then by the anti-commutation relation \eqref{anticomm3} we
have that
$\ovpsi_ {n-k+1/2}\psi_{k-1/2}= -\psi_{k-1/2} \ovpsi_ {n-k+1/2}$
and $\psi_{-k-1/2}\ovpsi_{n+k+1/2}= - \ovpsi_{n+k+1/2}\psi_{-k-1/2}$.
Thus the normal ordering is non-trivial only for $L^{(0)}_0$ and we
have that
\begin{equation}
\label{lzero1}
L^{(0)}_0= \sum_{\nu>0}(\nu+1/2)\ovpsi_ {-\nu}\psi_\nu +
\sum_{\nu>0}(\nu -1/2)\psi_{-\nu}\ovpsi_{\nu}.
\end{equation}
In the following, we often use the notation $T(z)$, $L_n$ instead of
$T^{(0)}(z)$ and $L_n^{(0)}$. By \eqref{3.2} and \eqref{3.3} we
have that
\begin{eqnarray}
\label{ln}
[ L_n, \psi_\nu]  &=&  -(n+ \nu+1/2)\psi_{n+\nu}, \\
\label{lnbar}
[ L_n, \ovpsi_\nu ] & =& - (n+\nu+1/2)\ovpsi_{n+\nu} + (n+1)\ovpsi_{n+\nu}.
\end{eqnarray}
In particular we have that
\begin{equation}
\label{lzero2}
[L_0, \psi_\nu]= -(\nu+1/2)\psi_\nu, \quad
[L_0, \ovpsi_\nu]= -(\nu-1/2)\ovpsi_\nu .
\end{equation}

Note that the degree of  the element
$$
|u \rangle =\ovpsi_{\mu_1}\ovpsi_{\mu_2}\cdots\ovpsi_{\mu_r}
 \psi_{\nu_s}\psi_{\nu_{s-1}}\cdots\psi_{\nu_1}|0\rangle
$$
with $\mu_1<\mu_2<\cdots <\mu_r<0$,
$\nu_1<\nu_2<\cdots<\nu_s<0$
is given by
\begin{equation}
\label{degform}
d(|u \rangle) = -\sum_{i=1}^r \mu_i - \sum_{j=1}^s \nu_j -p^2/2
\end{equation}
where $p=r-s$ is the charge of the element.  By \eqref{lzero1}  and
\eqref{lzero2} we thus see that
$$
L_0| u\rangle = (d+p(p+1)/2) | u\rangle.
$$
Hence we have the following lemma.
\begin{Lem}\label{lem5.1}
On $\cF_d(p)$ the operator $L_0$ acts by multiplication by the scalar
$d+p(p+1)/2$.
\end{Lem}
\begin{Lem}
\label{lem5.2}
For any element $| u \rangle \in \cF_d(p)$ and any half integer
$\nu \in \hZ$ we have that
$$
\deg(\psi_\nu |u\rangle ) = d +p -(\nu +1/2),\quad
\deg (\ovpsi_\nu | u \rangle ) = d -p -(\nu +1/2).
$$
\end{Lem}

{\it Proof}. \quad Note that the charge of $\psi_\nu |u\rangle$ is $p-1$ and
the charge of $\ovpsi_\nu | u \rangle$ is $p+1$. By \eqref{lzero2} we have
$$
L_0 (\psi_\nu | u \rangle) =\psi_\nu ( L_0 | u \rangle) - (\nu +1/2)\psi_\nu | u \rangle
=\{ d+p(p+1)/2 - (\nu +1/2)\}\psi_\nu | u \rangle .
$$
Hence, we get that
$$
\deg (\psi_\nu | u \rangle)= d+p(p+1)/2 - (\nu +1/2)- p(p-1)/2 =
d +p -(\nu +1/2).
$$
By the same way we can show that
$$
L_0 (\ovpsi_\nu | u \rangle)
=\{ d+p(p+1)/2 - (\nu -1/2)\}\psi_\nu | u \rangle.
$$
Hence
$$
\deg (\ovpsi_\nu | u \rangle ) = d -p -(\nu +1/2).
$$
{}\QED

Let
\begin{eqnarray*}
&&|u \rangle =\ovpsi_{\mu_1}\ovpsi_{\mu_2}\cdots\ovpsi_{\mu_r}
 \psi_{\nu_s}\psi_{\nu_{s-1}}\cdots\psi_{\nu_1}|0\rangle\\
&&\mu_1<\mu_2<\cdots<\mu_r<0, \quad
\nu_1<\nu_2<\cdots<\nu_s<0,
\end{eqnarray*}
and $M \in {\mathcal M}(p)$ be the Maya diagram corresponding  to $|u \rangle$. Then
$$
| M \rangle = (-1)^{\sum_{i=1}^s \nu_i +s/2}
\ovpsi_{\mu_1}\ovpsi_{\mu_2}\cdots\ovpsi_{\mu_r}
 \psi_{\nu_s}\psi_{\nu_{s-1}}\cdots\psi_{\nu_1}|0\rangle.
$$

Let us introduce a bilinear pairing
$$
(\phantom{X}|\phantom{X}) : \cF(p) \times \cF(p) \rightarrow \bC
$$
by putting
$$
(|M \rangle | | N \rangle ) = \langle M |N \rangle = \delta_{M,N}
$$
for Maya diagrams $M,N \in {\mathcal M}(p)$. The pairing is perfect
on $\cF_d(p) \times \cF_d(p)$ and zero on $\cF_d(p) \times \cF_{d'}(p)$
for $d \ne d'$.  Note that the bilinear pairing $(\phantom{X}|\phantom{X})$
is symmetric.
\begin{Lem}
\label{lem5.3}
\begin{eqnarray*}
(\psi_\nu |u\rangle \, | \, | v \rangle ) &= &
(| u \rangle \, | \, \ovpsi_{-\nu} | v \rangle ), \quad
\hbox{for $|u\rangle \in \cF(p+1)$, $|v \rangle \in \cF(p)$,} \\
(\ovpsi_\nu |u\rangle \, | \, | v \rangle ) &= &
(| u \rangle \, | \, \psi_{-\nu} | v \rangle ), \quad
\hbox{for $|u\rangle \in \cF(p-1)$, $|v \rangle \in \cF(p)$.}
\end{eqnarray*}
\end{Lem}

{\it Proof}. \quad It is enough to consider the case in which
$|u \rangle$ and $|v \rangle$ are Maya diagrams.
Assume
$$
|M\rangle = (-1)^{\sum_{i=1}^s \nu_i +s/2}
\ovpsi_{\mu_1}\ovpsi_{\mu_2}\cdots\ovpsi_{\mu_r}
 \psi_{\nu_s}\psi_{\nu_{s-1}}\cdots\psi_{\nu_1}|0\rangle,
$$
with
$\mu_1<\mu_2<\cdots<\mu_r<0$, $\nu_1<\nu_2<\cdots<\nu_s<0$.
Then we have that
$$
\psi_\nu |M \rangle =
(-1)^{\sum_{i=1}^s \nu_i +s/2+r}
\ovpsi_{\mu_1}\ovpsi_{\mu_2}\cdots\ovpsi_{\mu_r}
\psi_\nu \psi_{\nu_s} \psi_{\nu_{s-1}}\cdots\psi_{\nu_1}|0\rangle.
$$
For simplicity assume that $\nu > \nu_s$. Then
$(\psi_\nu |M\rangle \, | \, |N\rangle) \ne 0$ if and only if the Maya diagram
$N$ has the form
$$
|N\rangle = (-1)^{\sum_{i=1}^s \nu_i +\nu +s/2+1/2}
\ovpsi_{\mu_1}\ovpsi_{\mu_2}\cdots\ovpsi_{\mu_r}
\psi_\nu \psi_{\nu_s} \psi_{\nu_{s-1}}\cdots\psi_{\nu_1}|0\rangle.
$$
In this case
$$
(\psi_\nu|M \rangle \, |\, |N\rangle ) = (-1)^{\nu + 1/2+r}.
$$
On the other hand
$$
\ovpsi_{-\nu }| N \rangle = (-1)^{\nu +1/2+r}|M \rangle.
$$
Hence
$$
(\psi_\nu | M \rangle \, |\, |N\rangle ) = (|M \rangle \, |\, \ovpsi_{-\nu} | N\rangle ).
$$
Other cases can be treated similarly. \QED

Let $r : \cF(p) \rightarrow \cF(-p)$ be a mapping defined
by taking the mirror image
of the Maya diagram with respect to $0$ and interchange white and black.
Namely, if
$$
| M \rangle = (-1)^{\sum_{i=1}^s \nu_i +s/2}
\ovpsi_{\mu_1}\ovpsi_{\mu_2}\cdots\ovpsi_{\mu_r}
 \psi_{\nu_s}\psi_{\nu_{s-1}}\cdots\psi_{\nu_1}|0\rangle,
$$
with
$\mu_1<\mu_2<\cdots<\mu_r<0$, $\nu_1<\nu_2<\cdots<\nu_s<0$,
then  $r(M)$ is defined by
$$
|r(M) \rangle = (-1)^{\sum_{j=1}^r \mu_j + r/2}
\ovpsi_{\nu_1}\ovpsi_{\nu_2}\cdots\ovpsi_{\nu_s}
\psi_{\mu_r}\psi_{\mu_{r-1}}\cdots\psi_{\mu_1}|0\rangle .
$$
Note that we have
$$
\charge (r(M))= - \charge (M)\mbox{ and } \deg(r(M))= \deg(M).
$$
Hence, we have a linear map
$$
r : \cF_d(p) \rightarrow \cF_d(-p).
$$
Now for any integer $p$ let us  introduce a bilinear pairing
$$
\{\phantom{X}|\phantom{X}\} : \cF_d(p) \times \cF_d(-p) \rightarrow  \bC.
$$
For that purpose put
$$
\alpha(n) = \left\{
\begin{array}{cl}
1, & n \equiv 1, \, 2 \pmod{4} \\
-1, & n \equiv 0, \, 3 \pmod{4}
\end{array}
\right.
$$
Then it is easy to show that we for any $n \in \bZ$ get that
\begin{equation}
\label{eqalpha}
\alpha(n+1) = (-1)^{n+1}\alpha(n).
\end{equation}
The bilinear pairing $\{\phantom{X}|\phantom{X}\}$ is defined by
\begin{equation}
\label{bilinear2}
\{| u \rangle \, |\, | v \rangle  \}  = \alpha(\charge(| u \rangle
)) (|u\rangle \, | \, r(|v\rangle )).
\end{equation}
In particular for Maya diagrams $M \in {\mathcal M}(p)$ and $N \in {\mathcal M}(-p)$ we
have that
$$
\{|M \rangle | | N \rangle \} = \alpha(p) \langle M |
 r(N) \rangle =\alpha(p)  \delta_{M,r(N)}.
$$

Note that the pairing $\{\phantom{X}|\phantom{X}\}
 : \cF_d(p) \times \cF_d(-p) \rightarrow  \bC$ is perfect and that
it  can be extended to $\cF(p) \times \cF(-p)  \rightarrow  \bC$ by requiring that
the pairing is zero on $\cF_d(p) \times \cF_{d'}(-p) $ for $d\ne d'$.
Now from the above Lemma \ref{lem5.3} we obtain the following results.
\begin{Lem}
\label{lem5.4}
\begin{eqnarray*}
\{\psi_\nu | u \rangle \,|\, | v \rangle \} &=& (-1)^{-\nu-1/2}
\{ | u \rangle \,|\, \psi_{-\nu} | v \rangle\},  \\
\{ \ovpsi_\nu | u \rangle \,|\, | v \rangle \} &=&
(-1)^{-\nu-1/2}
\{ | u \rangle \,|\, \ovpsi_{-\nu} | v \rangle\}.
\end{eqnarray*}
\end{Lem}

{\it Proof}. \quad
By the definition of the mapping $r$ we have that
\begin{eqnarray*}
\psi_{\nu}r(|v \rangle ) &=&  (-1)^{-\nu+1/2+\charge(|v\rangle)}
r(\ovpsi_{-\nu}|v \rangle), \\
\ovpsi_{-\nu}r(|v \rangle ) &=&  (-1)^{-\nu+1/2+\charge(|v\rangle) }
r(\psi_{-\nu}|v \rangle).
\end{eqnarray*}
Hence by  Lemma \ref{lem5.3} and \eqref{eqalpha} for $|u\rangle \in \cF(p+1)$
and $|v \rangle \in \cF(-p)$ we see that
\begin{eqnarray*}
\{\psi_\nu | u \rangle \,|\, | v \rangle \} &=& \alpha(p) (\psi_\nu
| u \rangle \, |\, r( |v \rangle)) =
\alpha(p) (
| u \rangle \, |\, \ovpsi_{-\nu} r( |v \rangle)) \\
&=& (-1)^{-\nu +1/2+p}\alpha(p)(
| u \rangle \, |\, r(\psi_{-\nu}  |v \rangle)) \\
&=&  (-1)^{-\nu -1/2}\alpha(p+1)
( | u \rangle \, |\, r(\psi_{-\nu}  |v \rangle))  \\
&=& (-1)^{-\nu -1/2}\{| u \rangle \, |\, \psi_{-\nu}  |v \rangle\}.
\end{eqnarray*}
By a similar argument we can show the second equality.
 \QED

Let us introduce the shift operator $s : \cF(p) \rightarrow \cF(p+1)$
by defining
$$
s(\mu)(\nu) = \mu(\nu - 1) + 1,
$$
for a Maya diagram  $M$ with characteristic function $\mu$.
E.g.
\begin{eqnarray*}
&& s(\ovpsi_{\mu_1}\ovpsi_{\mu_2}\cdots\ovpsi_{\mu_r}
 \psi_{\nu_s}\psi_{\nu_{s-1}}\cdots\psi_{\nu_1}|0\rangle)   \\
&& \quad =\ovpsi_{\mu_1-1}\ovpsi_{\mu_2-1}\cdots\ovpsi_{\mu_r-1}
 \psi_{\nu_s+1}\psi_{\nu_{s-1}+1}\cdots\psi_{\nu_1+1}| 1\rangle \\
&& \quad =(-1)^s \ovpsi_{\mu_1-1}\ovpsi_{\mu_2-1}\cdots\ovpsi_{\mu_r-1}
 \ovpsi_{-1/2} \psi_{\nu_s+1}\psi_{\nu_{s-1}+1}\cdots\psi_{\nu_1+1}|
 0\rangle\\
 && \quad \quad + \delta_{\nu_s,-\frac12} \ovpsi_{\mu_1-1}\ovpsi_{\mu_2-1}\cdots\ovpsi_{\mu_r-1}
  \psi_{\nu_{s-1}+1}\cdots\psi_{\nu_1+1}|
 0\rangle
\end{eqnarray*}
with
$\mu_1<\mu_2<\cdots<\mu_r<0$, $\nu_1<\nu_2<\cdots<\nu_s<0$.

Now the bilinear pairings
$$
\{\phantom{X}|\phantom{X}\}_+  : \cF(p) \times \cF(-p -1)
\rightarrow \bC
$$
is defined by
\begin{equation}
\label{inn+}
 \{| u \rangle \, | \,  | v \rangle\}_+  =  \{|u\rangle \,|\, s(| v \rangle) \}
=\alpha(\charge(|u \rangle))(| u \rangle,\, | \, r( s( | v \rangle))  .
\end{equation}
Note that the pairing $\{\phantom{X}|\phantom{X}\}_+$ is perfect
on $\cF_d(p) \times \cF_d(-p-1)$ and zero on
$\cF_d(p) \times \cF_{d'}(-p-1)$ for $d \ne d'$.

Finally we obtain the following result which play an important role
in the sewing procedure.
\begin{Lem}
\label{lem5.5}
We have the following equalities
\begin{eqnarray*}
\{\psi_\nu | u \rangle \, | \, | v \rangle \}_+ &=&  (-1)^{-\nu -1/2}
\{ | u \rangle \,|\, \psi_{-\nu-1} | v \rangle\}_+, \\
\{\ovpsi_\nu | u \rangle \,|\, | v \rangle \}_+  &=& (-1)^{-\nu - 1/2}
\{ | u \rangle \,|\, \ovpsi_{-\nu+1} | v \rangle\}_+.
\end{eqnarray*}
\end{Lem}

{\it Proof}. \quad These results follow from
 Lemma \ref{lem5.4} and \eqref{inn+}. \QED

Let us now consider a one-parameter family of
curves which is a family of smooth curves degenerating to a single nodal curve.

Let $C_0$ be a complete curve with only one ordinary double point $P$
such that $C_0 \setminus \{P\}$ is non-singular.
Let $Q_1$, $Q_2, \ldots, Q_N$ be distinct non-singular points on $C_0$.
Let $\nu : \widetilde{C}_0
\rightarrow C_0$ be the normalization of the singular curve. Put $\{P_+, P_-\}
= \nu^{-1}(P)$.

\begin{Lem}
\label{lem5.6}
There exist a meromorphic vector field $\tilde{l} \in H^0(\widetilde{C}_0,
\Theta_{C_0}(*\sum_{j=1}^N Q_j))$ and a coordinate neighbourhood $U_+$ (resp. $U_-$) of
the point $P_+$ (resp. $P_-$) with local coordinate $z$ with center $P_+$ (resp.
$w$ with center $P_-$) such that we have
\begin{eqnarray*}
U_+ =  \{ \; z \; | \; |z| < 1 \;\}, &\quad &
U_-= \{ \; w \; | \; |w| < 1 \;\} \\
\tilde{l}|_{U_+} = \frac12 z \frac{d}{dz}, & \quad &
\tilde{l}|_{U_-} = \frac12 w \frac{d}{dw} .
\end{eqnarray*}
\end{Lem}

This lemma is a special case of lemma 5.3.1 in \cite{U2}. For the
convenience of the read we give the proof in this special case
here.

{\it Proof}. \quad Choose coordinate neighbourhood $U_+$ (resp. $U_-$) and
a local coordinate $x$ with center $P_+$ ($y$ with center $P_-$). Then we can find
a meromorphic vector field $\tilde{l} \in H^0(\widetilde{C}_0,
\Theta_{C_0}(*\sum_{j=1}^N Q_j))$
such that
$$
\tilde{l}|_{U_+} = (\frac12 x + a_2x^2 + \cdots )\frac{d}{dx}, \quad
\tilde{l}|_{U_-} = (\frac12 y + b_2 y^2 + \cdots ) \frac{d}{dy} .
$$
Then, by solving the following differential equations
\begin{eqnarray*}
\tilde{l} (z) &=& \frac12  z \frac{dz}{dx}, \quad z(0) = 0 , \\
\tilde{l} (w) &=& \frac12  w \frac{dw}{dy}, \quad w(0) = 0,
\end{eqnarray*}
we obtain  holomorphic functions $z$ and $w$ in neighbourhoods of
 $P_+$ and $P_-$, respectively.
It is easy to show that $z$ and $w$ give local coordinates with center $P_+$ and
$P_-$, respectively. Moreover, choosing $U_{\pm}$
smaller if necessary,  by the differential equations
we have
$$
\tilde{l}|_{U_+} = \frac12 z \frac{d}{dz}, \quad
\tilde{l}|_{U_-} = \frac12 w \frac{d}{dw} .
$$
Replacing $z$ by $cz$ with $c\ne 0$ and $w$ by $c'w$ with $c' \ne 0$
we may assume that
$$
U_+= \{ \; z \; | \; |z| < 1 \;\}, \quad
U_-= \{ \; w \; | \; |w| < 1 \;\}
$$
{}\QED

Let us now construct a flat one-parameter deformation of $C_0$. We
do this just like in section 5.3 in \cite{U2}.
Let
\begin{eqnarray*}
D &=&  \{ \; q \in {\bf C}\; | \; |q| < 1\; \} \\
S &=& \{ \; (x,y,q) \in {\bf C}^3\; | \; xy=q, |x|<1, |y|<1, |q|<1\; \} \\
Z &=& \{ \; (P,q) \in \widetilde{C}_0 \times D\; | \; P \in \widetilde{C}_0
\setminus (U_+ \cup U_-), \hbox{ or } P \in U_+, |z(P)|>|q|\;\} \\
W &=& \{ \; (Q,q) \in \widetilde{C}_0 \times D\; | \; Q \in \widetilde{C}_0
\setminus (U_+ \cup U_-), \hbox{ or } Q \in U_-, |w(Q)|>|q|\;\}
\end{eqnarray*}
Now introduce an equivalence relation $\sim$ on $Z\sqcup S \sqcup W$ as follows.

\begin{eqnarray*}
\hbox{\rm 1)} &\quad& Z \cap (U_+ \times D) \ni (P,q)  \sim (x,y,q') \in S
\Longleftrightarrow
(x,y,q') = (z(P), q/z(P), q) \\
\hbox{\rm 2)} &\quad&  W \cap (U_- \times D) \ni (Q,q)  \sim (x,y,q') \in S
\Longleftrightarrow
(x,y,q') = (q/w(P), w(P), q) \\
\hbox{\rm 2)} &\quad&  Z \ni (P,q)  \sim (Q,q') \in W
\Longleftrightarrow
(P,q) = (Q, q')
\end{eqnarray*}
Let ${\mathcal C}$ be the two-dimensional complex manifold obtained as the quotient of
$Z\sqcup S \sqcup W$ by this equivalence relation. There is a proper holomorphic
map $\pi : {\mathcal C} \rightarrow D$ such that the fiber over the origin is
$C_0$ and for $q \ne 0$ the fiber $C_q=\pi^{-1}(q)$ is a non-singular curve of genus
$g+1$ or $g$ according as $\widetilde{C}_0$ is connected or not, where
the genus of $\widetilde{C}_0$ is $g$.
The family $\pi : {\mathcal C} \rightarrow
D$ is a flat deformation (or smoothing) of $C_0$.
Let us assume that $Q_j \notin U_+ \cup U_-$. Then,
the point $Q_j\in C_0 \setminus \{P\}$ defines a holomorphic section
$\sigma_i : D \rightarrow {\mathcal C}$.
Let us choose a local coordinate $\xi_j$ with center $Q_i$.
Put $\gF = (\pi : \cC \rightarrow D;
\sigma_1, \ldots, \sigma_N; \xi_1, \ldots, \xi_N)$.
The meromorphic vector field $\tilde{l}$ given in Lemma \ref{lem5.6}
defines  local vector fields
$$
\vec{l} =(l_1(\xi_1)\frac{d}{d\xi_1}, l_2(\xi_2)\frac{d}{d\xi_2},
\ldots, l_N(\xi_N)\frac{d}{d\xi_N})
$$
where $l_j(\xi_j)\frac{d}{d\xi_j}$ is the Laurent expansion of $\tilde{l}$
with respect to the coordinate $\xi_j$ with center $Q_j$.
\begin{Lem}
\label{lem5.7}
The local vector fields $\vec{l}$ is an element of $\cL(\gF)$ and
we have
$$
\theta(\vec{l}) = q \frac{d}{dq}.
$$
\end{Lem}

{\it Proof}. \quad The meromorphic vector field
$$
\tilde{l} -q\frac{q}{dq}
$$
is a holomorphic section of $\pi_*\Theta_{\cC}'(*\sum_{j=1}^N \sigma_j(D))_\pi$ over $D$.
This show the first part of the proposition.
A proof of the  second part can be found in
[U2,  Corollary 5.3.3 ]. \QED

\begin{Thm}
\label{thm5.1}
Put
$$
\widetilde{\gX}= (\widetilde{C}_0; P_+,P_-, Q_1,
\ldots, Q_N, z,w,\xi_1,\ldots,\xi_N)
$$
and
$$
\gX = (C_0; Q_1, \ldots, Q_N, \xi_1,\ldots,\xi_N).
$$
Then for any element $\langle \Phi | \in \cVd_{\ab}(\widetilde{\mathfrak{X}})$
there exists a unique global flat section $\langle \widetilde{\Phi}(q) |$ of the
sheaf of ghost vacua $\cVd_{\ab}(\gF)$ over the disk such that
$\langle \widetilde{\Phi}(0) |$ is the element of $\cVd_{\ab}(\gX)$ corresponding
to the element $\langle \Phi | $, i.e.

$$
\langle \widetilde{\Phi}(0)|u\rangle = \langle \Phi|0_{+,-}\otimes u \rangle
$$
for any $|u \rangle \in \cF_N$.
\end{Thm}
{\it Proof}. \quad
Let $\{v_i(d,p)\}_{i=1,\ldots,m_d}$ be a basis of $\cF_d(p)$
for any $p \in \bZ$ and
$\{v^i(d,p)\}_{i=1,\ldots,m_d}$ be the dual basis of $\cF_d(-p-1)$
with respect to the pairing $\{\phantom{X}|\phantom{X}\}_+$.
Note that $\cF_0(p)$ is spanned by $|p\rangle$ over $\bC$ and
we choose $v_1(0,p) = | p\rangle$ for all $p \in \bZ$. Then
$v^1(0,-p-1) = \alpha(p)|-p-1\rangle$.

Define $\langle \widetilde{\Phi}(q) |$ by
\begin{equation}
\label{solution}
\langle  \widetilde{\Phi}(q) | u \rangle =
\sum_{p \in \bZ}\Bigl\{ \sum_{d=0}^\infty \sum_{i=1}^{m_d}(-1)^{p+d}
\langle \Phi |v_i(d,p)\otimes v^i(d,-p-1)\otimes u\rangle \Bigr\}q^{d+p(p+1)/2}.
\end{equation}
Then by Theorem 3.5 $\langle \widetilde{\Phi}(0)|  \in \cVd_{\ab}(\gX)$
is the element corresponding to $\langle \Phi | $.

Let us show that $\langle \widetilde{\Phi}(q) |$  satisfies the gauge conditions.
Since $\langle \widetilde{\Phi}(q) |$ is defined as a formal power series,  we
will first prove the formal gauge conditions. Later we shall show that
$\langle \widetilde{\Phi}(q) |$ satisfies a differential equation of Fuchsian
type so that it converges. Then, the formal gauge conditions
implies the usual gauge conditions.

Choose $\tau \in \pi_*(\omega_{\cC/\cB}(*S))$. In a neighbourhood of the double
point $P$ in $\cC$, $\tau$ can be expressed in the form
$$
\tau = \bigl(\sum_{n,m \geq 0} a_{n,m}z^nw^m\bigr)\bigl[dz \wedge dw/dq\bigr]
$$
where $\bigl[dz \wedge dw/dq\bigr]$ is a local basis of $\omega_{\cC/\cB}$
such that
$$
\nu^*(\bigl[dz \wedge dw/dq\bigr])|_{U_+} = -\frac{dz}{z}, \quad
\nu^*(\bigl[dz \wedge dw/dq\bigr])|_{U_-} = \frac{dw}{w}.
$$
Now put
\begin{eqnarray*}
\tau_+& = &\bigl(\sum_{n,m \ge 0} a_{n,m}z^n\bigl(\frac{q}{z}\bigr)^m\bigr)
(-\frac{dz}{z}) =
\sum_{n=0}^\infty \bigl\{-\sum_{m=0}^\infty a_{m,n}z^{m-n-1} dz\bigr\}q^n
\\
\tau_-&=&\bigl(\sum_{n,m \ge 0} a_{n,m}\bigl(\frac{q}{w}\bigr)^nw^m\bigr)
\frac{dw}{w}
=\sum_{n=0}^\infty \bigl\{\sum_{m=0}^\infty a_{n,m}w^{m-n-1}
dw)\bigr\}q^n,
\end{eqnarray*}
and
\begin{eqnarray*}
\tau^{(n)}_+ &= & -\bigl(\sum_{m=0}^\infty a_{m,n}z^{m-n-1}\bigr) dz, \\
\tau^{(n)}_-&=& \bigl(\sum_{m=0}^\infty a_{n,m}w^{m-n-1} \bigr)dw.
\end{eqnarray*}
Then
$$
\tau_+ = \sum_{n=0}^\infty \tau^{(n)}_+ q^n, \quad
\tau_- = \sum_{n=0}^\infty \tau^{(n)}_- q^n.
$$
It is easy to show that there is a unique $\tau^{(n)} \in
H^0(\widetilde{C}_0, \omega_{\widetilde{C}_0}(*(P_++P_-+\sumQj))$ such that $\{\tau^{(n)}_\pm$
is the Taylor expansion of $\tau^{(n)}$ at $P_\pm$.
the data $\{\tau^{(n)}_+, \tau^{(n)}_-, \tau^{(n)}(\xi)\}$
define a meromorphic one-form .
Hence by the first gauge condition for $\langle \Phi |$ we see
that
\begin{eqnarray*}
&&
\langle \Phi | \psi[\tau^{(n)}_+]v_i(d,p)\otimes v^i(d,-p-1)
\otimes u \rangle \\
&&\phantom{XXX} +
(-1)^p \langle \Phi | v_i(d,p)\otimes \psi[\tau^{(n)}_-]v^i(d,-p-1)
\otimes u \rangle \\
&&\phantom{XXX} - \sum_{j=1}^N
\langle \Phi | v_i(d,p)\otimes v^i(d,-p-1)
\otimes  \rho_j(\psi[\tau^{(n)}_j])  u \rangle  =0.
\end{eqnarray*}
Hence
\begin{eqnarray}
&&\sum_{j=1}^N \langle  \widetilde{\Phi}(q)  | \rho_j(\psi[
\sum_{n=0}^\infty \tau^{(n)}_j q^n] )  u \rangle =
\sum_{j=1}^N \sum_{n=0}^\infty  \langle  \widetilde{\Phi} (q) | \rho_j(
\psi[\tau^{(n)}_j ]   )  u \rangle q^n  \nonumber \\
&&=
\sum_{p \in \bZ} \sum_{n=0}^\infty
\sum_{d=0}^\infty \sum_{i=1}^{m_d} (-1)^{p+d}  \Big\{  \langle\Phi |
\psi[\tau^{(n)}_+] v_i(d, p) \otimes v^i(d,-p-1) \otimes u\rangle
\nonumber \\
&&+(-1)^p \langle\Phi |
v_i(d, p) \otimes \psi[\tau^{(n)}_-]  v^i(d,-p-1) \otimes u\rangle
\Bigr\}q^{d+n+p(p+1)/2}.  \label{fineq}
\end{eqnarray}
Note that
\begin{eqnarray*}
\psi[\tau^{(n)}_+] &= & - \sum_{m=0}^\infty a_{m,n}\psi_{m-n-1/2},  \\
\psi[\tau^{(n)}_-] &= & \sum_{m=0}^\infty a_{n,m}\psi_{m-n-1/2}.
\end{eqnarray*}
Then by Lemma \ref{lem5.5}
\begin{eqnarray*}
&& \sum_{i=1}^{m_d} \langle \Phi |
\psi[\tau^{(n)}_+] v_i(d, p) \otimes v^i(d,-p-1) \otimes u\rangle
q^{d+n+p(p+1)/2}  \\
&&= - \sum_{i=1}^{m_d}  \sum_{m=0}^\infty a_{m,n}
\langle\Phi |
\psi_{m-n-1/2} v_i(d, p) \otimes v^i(d,-p-1) \otimes u\rangle
q^{d+n+p(p+1)/2}   \\
&& =- \sum_{m=0}^\infty  \sum_{i=1}^{m_d}  \sum_{j=1}^{m_{d+p-m+n}} a_{m,n}
\{\psi_{m-n-1/2} v_i(d, p) \,|\, v^j(d+p-m+n, -p)\}_+ \\
&& \phantom{XX}\cdot \langle\Phi | v_j(d+p-m+n, p-1)
\otimes v^i(d,-p-1) \otimes u\rangle
q^{d+n+p(p+1)/2}\\
&&= - \sum_{m=0}^\infty \sum_{i=1}^{m_d} \sum_{j=1}^{m_{d+p-m+n}}
(-1)^{n-m} a_{m,n}
\{v_i(d , p) \,|\, \psi_{n-m-1/2} v^j(d+p-m+n, -p)\}_+ \\
&& \phantom{XX}\cdot \langle\Phi | v_j(d+p-m+n, p-1) \otimes
v^i(d,-p-1) \otimes u \rangle
q^{d+n+p(p+1)/2}\\
&& =  \sum_{m=0}^\infty \sum_{j=1}^{m_{d+p-m+n} }
(-1)^{n-m+1}a_{m,n}
\langle\Phi | v_j(d+p-m+n, p-1) \\
&& \phantom{XXXXXXXX} \otimes \psi_{n-m-1/2} v^j(d+p-m+n, -p)
\otimes  u \rangle q^{d+n+p(p+1)/2}.
\end{eqnarray*}
Hence
\begin{eqnarray*}
&&\sum_{n=0}^\infty  \sum_{i=1}^{m_d} (-1)^{d+p}\langle \Phi |
\psi[\tau^{(n)}_+] v_i(d, p) \otimes v^i(d,-p-1) \otimes u\rangle
q^{d+n+p(p+1)/2}  \\
&&  = \sum_{n=0}^\infty \sum_{m=0}^\infty \sum_{j=1}^{m_{d+p-n+m} }
(-1)^{d-n+m+p+1} a_{n,m} \langle\Phi | v_j(d+p-n+m, p-1) \\
&& \phantom{XXXXXXXX} \otimes \psi_{m-n-1/2} v^j(d+p-n+m, -p)
\otimes  u \rangle q^{d+m+p(p+1)/2} .
\end{eqnarray*}
Put
$$
\widetilde{d} = d+p-n+m.
$$
Then
$$
\widetilde{d} + n + p(p-1)/2 = d+m+p(p+1)/2.
$$
Hence the right hand side of the last equality can be rewritten in the form
\begin{eqnarray*}
&& \sum_{n=0}^\infty \sum_{m=0}^\infty \sum_{j=1}^{m_{\widetilde{d}} }
(-1)^{\widetilde{d}+1} a_{n,m}
\langle\Phi | v_j(\widetilde{d}, p-1) \otimes \psi_{m-n-1/2} v^j(\widetilde{d}, -p)
\otimes  u \rangle q^{\widetilde{d}+n+p(p-1)/2}\\
&& = \sum_{n=0}^\infty  \sum_{j=1}^{m_{\widetilde{d}} }
(-1)^{\widetilde{d}+1}  \langle\Phi | v_j(\widetilde{d}, p-1)
\otimes \psi[\tau^{(n)}_-] v^j(\widetilde{d}, -p)
\otimes  u \rangle q^{\widetilde{d}+n+p(p-1)/2}.
\end{eqnarray*}
Hence
\begin{eqnarray*}
&&
\sum_{p \in \bZ} \sum_{n=0}^\infty
\sum_{d=0}^\infty \sum_{i=1}^{m_d} (-1)^{p+d}   \langle\Phi |
\psi[\tau^{(n)}_+] v_i(d, p) \otimes v^i(d,-p-1) \otimes u\rangle
q^{d+n+p(p+1)/2}
\\
&&=-\sum_{p \in \bZ}
\sum_{n=0}^\infty \sum_{d=0}^\infty \sum_{i=1}^{m_d} (-1)^d
\langle\Phi | v_j(d, p-1) \otimes \psi[\tau^{(n)}_-] v^j(d, -p)
\otimes  u \rangle q^{d+n+p(p-1)/2}\\
&& = - \sum_{p \in \bZ}
\sum_{n=0}^\infty \sum_{d=0}^\infty \sum_{i=1}^{m_d} (-1)^d
\langle\Phi | v_j(d, p) \otimes \psi[\tau^{(n)}_-] v^j(d, -p-1)
\otimes  u \rangle q^{d+n+p(p-1)/2}
\end{eqnarray*}
Therefore, by \eqref{fineq} we conclude that
$$
\sum_{j=1}^N \langle \widetilde{\Phi}(q)|
\rho_j(\psi[\sum_{n=0}^\infty \tau^{(n)}_j q^n]) u \rangle =0.
$$
This proves the first formal gauge condition. Similarly, we can
show that the second formal gauge condition also holds.

Next we shall show that the formal power series
\eqref{solution} is a formal solution of a Fuchsian differential equation.
Hence the power series \eqref{solution} indeed converges. Therefore the
formal gauge conditions are nothing but the usual gauge conditions.

Note that on the
punctured disk $D^*$ we have the connection introduced in section 4.
For the local vector field $\vec{l}$ defined by the meromorphic vector field
$\widetilde{l}$ in Lemma \ref{lem5.6} the corresponding connection has the form
\begin{equation}
\label{eqn5.11}
q\frac{d}{dq}(\langle\widetilde{\Phi}(q)|u \rangle) -
\sum_{j=1}^N \langle\widetilde{\Phi}(q)|\rho_j(T[\underline{l}_j])u \rangle
- \frac16 b_\omega(\vec{l})\langle\widetilde{\Phi}(q)|u \rangle =0
\end{equation}
where $\omega$ is an element of
$$
H^0({\mathcal C} \times_{D}{\mathcal C}, \omega_{{\mathcal C}/D}^{\boxtimes 2}
 (2\Delta))
$$
with
$$
\omega= \frac{dv du}{(v-u)^2} + \hbox{\rm holomorphic at $\Delta$},
$$
where $\Delta$ is the diagonal of ${\mathcal C} \times_{D}{\mathcal C}$
and if  $p_j : {\mathcal C} \times_{D}{\mathcal C}
\rightarrow {\mathcal C}$ is the projection onto the $j$'th factor,
then
$$
 \omega_{{\mathcal C}/D}^{\boxtimes 2} = p_1^*\omega_{{\mathcal C}/D}
 \otimes p_2^*\omega_{{\mathcal C}/D}.
$$

The formal correlation function $\langle \widetilde{\Phi}(q)|\widetilde{T}(z)|
u  \rangle dz^2$ which will be later proved to be an element of
$ H^0(({\mathcal C}, \omega_{{\mathcal C}/D})^{\otimes 2}   (*\sum_{j=1}^N s_j(D))$
is defined in \eqref{tildeem}:
$$
\langle \widetilde{\Phi}(q)|\widetilde{T}(z)|
u  \rangle dz^2 =
\lim_{w \rightarrow z} \frac12 \big\{
\langle \widetilde{\Phi}(q)|
\frac{d\psi(w)}{dw}\ovpsi(z)|u\rangle dw dz -
\omega(w,z)\langle \widetilde{\Phi}(q)|u\rangle dw dz\big\}.
$$
Then, by \eqref{diffem} we have that
\begin{equation}
\label{EMrel}
\langle \widetilde{\Phi}(q)|\widetilde{T}(z)|
u \rangle dz^2 =
\langle \widetilde{\Phi}(q)|T(z)|u \rangle dz^2
+ \frac{1}{6} \langle \widetilde{\Phi}(q)|u\rangle S_\omega(z)dz^2.
\end{equation}

Let $\tilde{l}= l(z)\displaystyle{\frac{d}{dz}}$ be the meromorphic
vector field defined in Lemma \ref{lem5.6}. For $q \ne 0$
$$
\tilde{l}\cdot\langle \widetilde{\Phi}(q)|\widetilde{T}(z)|
u \rangle dz^2 = l(z)\langle \widetilde{\Phi}(q)|\widetilde{T}(z)|
u \rangle dz
$$
is a meromorphic one-form on
$$
C_q' = C_q \setminus \{ \; (x,y,q) \in S\;| \; |x| \le \epsilon,
\; |y| \le \epsilon\;\}
$$
where $\epsilon<1$ is a sufficiently small positive number.

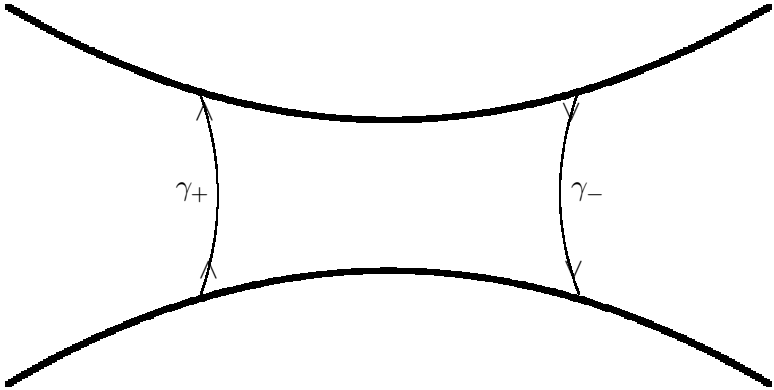
\begin{picture}(100,70)(-10,-10)
\linethickness{2pt}
\qbezier(0,50)(50, 20)(100,50)
\qbezier(0,0)(50, 30)(100, 0)
\thinlines
\qbezier(25, 11)(30, 25)(25,39)
\qbezier(75, 12)(70, 25)(75,39)
\put(72.7, 35){$\vee$}
\put(73, 14){$\vee$}
\put(25, 14){$\wedge$}
\put(24.5, 35){$\wedge$}
\put(22,25){$\gamma_+$}
\put(74,25){$\gamma_-$}
\end{picture}

The boundary of $C_q'$ consists of two disjoint curves
$\gamma_\pm$. We choose the orientation of $\gamma_\pm$ in such a way that
$C_q'$ lies to the left of $\gamma_\pm$. Then, we have that
\begin{eqnarray*}
&&\frac{1}{2\pi \sqrt{-1}}\int_{\gamma_+}
l(z)\langle \widetilde{\Phi}(q)|\widetilde{T}(z)|
u \rangle dz +
\frac{1}{2\pi \sqrt{-1}}\int_{\gamma_-}
l(w)\langle \widetilde{\Phi}(q)|\widetilde{T}(w)|
u \rangle dw \\
&&\phantom{XX}= \sum_{j=1}^N \res_{Q_j} \{
l(z)\langle \widetilde{\Phi}(q)|\widetilde{T}(z)|
u \rangle dz \} \\
&&
\phantom{XX}= \sum_{j=1}^N
\langle \widetilde{\Phi}(q)|\Res_{\xi_j =0}
(l (\xi_j) T(\xi_j)d\xi_j)u \rangle
+ \frac{1}{6} \sum_{j=1}^N \Res_{\xi_j = 0} \{ l (\xi_j) S_\omega(\xi_j)d\xi_j
\langle\widetilde{\Phi}(q)|u\rangle \} \\
&&
\phantom{XX}= \sum_{j=1}^N
\langle \widetilde{\Phi}(q)| \rho_j(T[l(\xi_j)])u\rangle +
b_\omega (\tilde{l}) \langle \widetilde{\Phi}(q)|u \rangle.
\end{eqnarray*}
On the curve $\gamma_+$
$$
\tilde{l} = \frac{1}{2} z \frac{d}{dz} .
$$
Hence
\begin{eqnarray*}
&&\frac{1}{2\pi \sqrt{-1}}\int_{\gamma_+}
l(z)\langle \widetilde{\Phi}(q)|\widetilde{T}(z)|
u \rangle dz  \\
&&\phantom{XX}=
\frac{1}{4\pi \sqrt{-1}}\int_{\gamma_+} z \big\{
\langle \widetilde{\Phi}(q)|T(z)|u \rangle
+ \frac{1}{6}S_\omega(z)\langle \widetilde{\Phi}(q) | u \rangle \big\}dz .
\end{eqnarray*}
Since $S_\omega(z)dz^2$ is holomorphic at $z=0$
$$
\frac{1}{2\pi \sqrt{-1}}\int_{\gamma_+}
l(z)\langle \widetilde{\Phi}(q)|\widetilde{T}(z)|
u \rangle dz
=
\frac{1}{4\pi \sqrt{-1}}\int_{\gamma_+}
z \langle \widetilde{\Phi}(q)|T(z)|u \rangle dz
$$
On the other hand
\begin{eqnarray*}
&&\frac{1}{4\pi \sqrt{-1}}\int_{\gamma_+}
z \langle \widetilde{\Phi}(q)|T(z)|u \rangle dz \\
&&\phantom{XX}=
\frac12 \sum_{p \in \bZ} \sum_{d=0}^\infty \sum_{i=1}^{m_d}
\frac{(-1)^{p+d}}{2\pi \sqrt{-1}} \int_{\gamma_+} z \langle \Phi|
T(z)|v_i(d,p)\otimes v^i(d,-p-1)\otimes u \rangle q^{d + p(p+1)/2} dz \\
&&\phantom{XX}=
\frac12 \sum_{p \in \bZ} \sum_{d=0}^\infty \sum_{i=1}^{m_d}
(-1)^{p+d}  \langle \Phi|L_0(v_i(d,p))\otimes v^i(d,-p-1)\otimes
 u \rangle q^{d + p(p+1)/2}\\
&&\phantom{XX}=
\frac12 \sum_{p \in \bZ} \sum_{d=0}^\infty \sum_{i=1}^{m_d}
(-1)^{p+d} (d + \frac{p(p+1)}{2})  \langle \Phi|v_i(d,p)\otimes v^i(d,-p-1)\otimes
 u \rangle q^{d + p(p+1)/2} .
\end{eqnarray*}
Similarly
\begin{eqnarray*}
&&
\frac{1}{2\pi \sqrt{-1}}\int_{\gamma_-}
l(w)\langle \widetilde{\Phi}(q)|\widetilde{T}(w)|
u  \rangle dw \\
&&\phantom{XX}=
\frac12 \sum_{p \in \bZ} \sum_{d=0}^\infty \sum_{i=1}^{m_d}
(-1)^{p+d }(d + \frac{p(p+1)}{2})  \langle \Phi|v_i(d,p)\otimes v^i(d,-p - 1)\otimes
 u \rangle q^{d + p(p+1)/2} .
\end{eqnarray*}
Thus we obtain that
\begin{eqnarray*}
&&\sum_{j=1}^N \langle \widetilde{\Phi}(q)|\rho_j(T[l(\xi_j)])u \rangle
+ b_\omega (\tilde{l})\langle \Phi|u \rangle \\
&&\phantom{XX}=
\sum_{p \in \bZ} \sum_{d=0}^\infty \sum_{i=1}^{m_d}
(-1)^{p+d} (d + \frac{p(p+1)}{2})  \langle \Phi|v_i(d,p)\otimes v^i(d,-p-1)\otimes
 u \rangle q^{d + p(p+1)/2} .
\end{eqnarray*}
On the other hand
$$
q\frac{d}{dq} \langle \widetilde{\Phi}(q)|u \rangle =
\sum_{p \in \bZ} \sum_{d=0}^\infty \sum_{i=1}^{m_d}
(-1)^{p+d }(d + \frac{p(p+1)}{2})  \langle \Phi|v_i(d,p)\otimes v^i(d,-p-1)\otimes
 u \rangle q^{d + p(p+1)/2} .
$$
Hence $\langle \widetilde{\Phi}(q)|u \rangle$ is a formal solution
of the differential equation
$$
q\frac{d}{dq} \langle \widetilde{\Phi}(q)|u \rangle -
\sum_{j=1}^N \langle \widetilde{\Phi}(q)|\rho_j(T[l(\xi_j)])u \rangle
-  \frac16 b_\omega(\tilde{l})\langle \Phi|u \rangle =0
$$
which is of Fuchsian type. Hence $\langle \widetilde{\Phi}(q)|u \rangle$
converges for all $u  \in \cF_N$.
\QED

Now let us give a proof of Theorem 3.2.

By \eqref{4.1} and Theorem \ref{thm4.1} $\dim \cVd_{\ab}(\gX) $ depends only on the genus
of the curve $C$ and the number $N$ of the points, if the curve $C$ is
non-singular.  Moreover, by
Theorem \ref{thm3.4} $\dim \cVd_{\ab}(\gX) $ depends only on
the genus of the curve $C$ provided $C$ is non-singular.
Put
$$
d(g) = \dim \cVd_{\ab}(\gX)
$$
for a non-singular curve $C$ of genus $g$.

Now assume that the curve $C$ has one double point $P$.  Let
$\gF = (\pi : \cC \rightarrow D; s_+,s_-, \xi_1, \ldots,\xi_N)$ be
the family constructed above so that $\gF(q)$, $q \ne 0$ is a
non-singular curve of genus $g$. Then, by \eqref{4.1} and
 Proposition \ref{prop4.1}
$$
d(g) = \dim \cVd_{\ab}(\gX) = \dim \cV_{\ab}(\gF(0)) \ge \dim \cV_{\ab}(\gF(q))
= \dim \cVd_{\ab}(\gF(q)).
$$
Moreover, by Theorem \ref{thm3.5} we have that
$$
\dim \cVd_{\ab}(\gX) = \dim \cVd_{\ab}(\widetilde{\gX}) = d(g-1).
$$
On the other hand by Theorem \ref{thm5.1}
$$
\dim \cVd_{\ab}(\gF(0))  \le \dim \cVd_{\ab}(\gF(q)).
$$
Hence
$$
d(g-1) = \dim \cVd_{\ab}(\gF(0))  \le \dim \cVd_{\ab}(\gF(q)) = d(g).
$$
If the curve $C$ has $m$ double points, applying Theorem \ref{thm3.5}
$m$ times  we conclude that $\dim \cVd_{\ab}(\gX)$ is independent of $\gX$.
In Example \ref{exmp3.1}, we showed that $\dim \cVd_{\ab}((\bP^1; 0; z)) =1$.
Hence, by Theorem  \ref{thm3.4} for any $\gX$ we have that
$\dim \cVd_{\ab}(\gX)=1$.  \QED

\begin{Cor}
\label{cor5.1}  The connection \eqref{eqn5.11} extends
holomorphically over $q=0$.
\end{Cor}

{\it Proof}. \quad  Since $\langle \widetilde{\Phi}(q)|$ is a single valued holomorphic
function and the differential equation \eqref{eqn5.11} is of Fuchsian type,
$$
\sum_{j=1}^N \langle\widetilde{\Phi}(q)|\rho_j(T[\underline{l}_j])u \rangle
+ \frac16 b_\omega(\vec{l})\langle\widetilde{\Phi}(q)|u \rangle
$$
must be divisible by $q$. Hence the differential equation has no singularity at $q=0$.
\QED

Theorem\ref{thm3.2} implies the following important fact.
\begin{Thm}
\label{thm5.1a}
For a family
$$
\gF = (\pi : \cC \rightarrow \cB; s_1, \ldots, s_N; \xi_1, \ldots, \xi_N)
$$
of $N$-pointed semi-stable curves with formal neighbourhoods the sheaf
$\cVd_{\ab}(\gF)$ of ghost vacua and the sheaf $\cV_{\ab}(\gF)$ of
dual ghost vacua are invertible $\cO_\cB$-modules. Moreover, they are dual
to each other.
\end{Thm}

{\it Proof}. \quad By Theorem\ref{thm3.2} and \eqref{4.1}, for any point
$b \in \cB$ the vector space $ \cV_{\ab}(\gF)\otimes_{\cO_\cB}
\cO_{\cB, b}/\frak{m}_b$ is one-dimensional. As $\cV_{ab}(\gF)$ is a coherent
$\cO_\cB$-module, this implies that $\cV_{ab}(\gF)$ is an invertible
$\cO_\cB$-module. Then, by \eqref{4.a} $\cVd_{ab}(\gF)$ is also
an invertible $\cO_\cB$-module.  \QED

Next we shall consider smoothing of  a family of nodal curves with sections and
formal coordinates. We need to generalize the above construction.

Let
$$
\widetilde{\gF} = ( \widetilde{\pi} : \widetilde{\cC} \rightarrow E ; s_1, \ldots,s_N,
t_+, t_-; \xi_1, \ldots, \xi_N, z, w)
$$
be a family of $(N+2)$-pointed smooth curves with formal coordinates.  Assume that $N \ge 1$
and $E$ is a small polydisk
$\{(u_1, \ldots, u_m)\in \bC^m \, | \, |u_i| < \varepsilon, \, 1\le i \le m\,\}$,
and $z$, $w$ are holomorphic coordinates such that there exists a meromorphic
vector field $\widetilde{l} \in H^0(\widetilde{C}, \Theta_{\widetilde{\cC}/E}(*S))$
with
$$
\widetilde{l}|_X = \frac12z \frac{d}{dz}, \quad
\widetilde{l}|_Y = \frac12w \frac{d}{dw}.
$$
There always exist such coordinates by Lemma \ref{lem5.6}.  Put
$$
S = \sum_{j=1}^N s_j(E), \quad S_\pm = t_+(E)  + t_-(E) .
$$
Then for any positive integer $M$ we have and exact sequence
$$
0 \rightarrow \Theta_{\widetilde{\cC}/E}(-S_\pm + MS)  \rightarrow
\Theta_{\widetilde{\cC}}(- S_\pm + MS)
\stackrel{d\widetilde{\pi}}{\rightarrow}  \widetilde{\pi}^* \Theta_E \rightarrow 0
$$
where $\Theta_{\widetilde{\cC}/E}$ is the sheaf of relative holomorphic vector fields and
$\Theta_Z$ denotes the sheaf of holomorphic vector fields on some complex manifold $Z$.
Put
$$
\Theta_{\widetilde{\cC}}(-S_\pm + MS)_{\widetilde{\pi}}
= d\widetilde{\pi}^{-1}(\widetilde{\pi}^{-1}\Theta_E).
$$
Taking the inductive limit on $M$ we can also define
$\Theta_{\widetilde{\cC}}(-S_\pm + *S)_{\widetilde{\pi}}$.

Similarly for any positive integer $M$ we have an exact sequence
$$
0 \rightarrow \Theta_{\widetilde{\cC}/E}(-S_\pm - S)  \rightarrow
\Theta_{\widetilde{\cC}/E}(- S_\pm + MS) \rightarrow
\bigoplus_{j=1}^N  \bigoplus_{k=0}^M \cO_E \xi_j^{-k}\frac{d}{d\xi_j} \rightarrow 0.
$$
{}From these exact sequences we have a commutative diagram
$$
\begin{array}{ccccccccc}
0 & \rightarrow & \widetilde{\pi}_*\Theta_{\widetilde{\cC}/E}(-S_\pm + *S) &
\rightarrow & \widetilde{\pi}_*\Theta_{\widetilde{\cC}}(- S_\pm +*S)_{\widetilde{\pi}}
 & \stackrel{d\widetilde{\pi}}{\rightarrow}
& \Theta_E  & \rightarrow &  0\\
&&&&&&&& \\
&&\parallel && \phantom{\scriptstyle{p}}\downarrow \scriptstyle{p} &&
\phantom{\scriptstyle{\rho}}\downarrow \scriptstyle{\rho} && \\
&&&&&&&& \\
0 & \rightarrow & \widetilde{\pi}_*\Theta_{\widetilde{\cC}/E}(-S_\pm + *S) &
\rightarrow & \bigoplus_{j=1}^N\cO_E[\xi_j^{-1}] \frac{d}{d\xi_j}  & \rightarrow &
R^1\widetilde{\pi}_*\Theta_{\widetilde{C}/E}(-S_\pm -S) &\rightarrow & 0
\end{array}
$$
For simplicity assume that
the $\cO_E$-homomorphism $\rho$ is injective.
Therefore, $p$ is also injective.
Let
$$
\widetilde{p} : \widetilde{\pi}_*\Theta_{\widetilde{\cC}}(- S_\pm + *S)_{\widetilde{\pi}}
\rightarrow \bigoplus_{j=1}^N\cO_E((\xi_j)) \frac{d}{d\xi_j}
$$
be the Laurent expansions of the $\frac{d}{d\xi_j}$ part along $s_j(E)$ and put
$$
\cL(\widehat{\gF}) = \widetilde{p}(\widetilde{\pi}_*
\Theta_{\widetilde{\cC}}(- S_\pm + *S)_{\widetilde{\pi}} )
\subset \bigoplus_{j=1}^N\cO_E((\xi_j)) \frac{d}{d\xi_j}.
$$
Then we have an exact sequence
$$
0   \rightarrow  \widetilde{\pi}_*\Theta_{\widetilde{\cC}/E}(-S_\pm + *S) \rightarrow
\cL(\widehat{\gF}) \stackrel{\widehat{\theta}}{\rightarrow}  \Theta_E \rightarrow 0.
$$
Hence we can find local vector fields
\begin{equation}
\label{mvector}
\vec{m}^{(i)} = \big(m_1^{(i)}\frac{d}{d\xi_1}, \ldots, m_N^{(i)}\frac{d}{d\xi_N}\big)
\end{equation}
which are sections of $\cL(\widehat{\gF})$ such that
$$
\widehat{\theta}(\vec{m}^{(i)}) = \frac{d}{du_i}.
$$

Note that by \eqref{lF} we can define
$$
\cL(\widetilde{\gF}) \subset \left(\bigoplus_{j=1}^N\cO_E((\xi_j)) \frac{d}{d\xi_j} \right)
\oplus \cO_E((z))\frac{d}{dz} \oplus  \cO_E((w))\frac{d}{dw}
$$
and we have the exact sequence
$$
0   \rightarrow  \widetilde{\pi}_*\Theta_{\widetilde{\cC}/E}( *S) \rightarrow
\cL(\widehat{\gF}) \stackrel{\widetilde{\theta}}{\rightarrow}  \Theta_E \rightarrow 0.
$$
Then the above argument shows that if we put
$$
\vec{\widetilde{m}}^{(i)} =
\big(m_1^{(i)}\frac{d}{d\xi_1}, \ldots, m_N^{(i)}\frac{d}{d\xi_N}, 0, 0 \big),
$$
then
\begin{equation}
\label{mtildevector}
\widetilde{\theta}(\vec{\widetilde{m}}^{(i)}) = \frac{d}{du_i}.
\end{equation}
Let us identify $t_+(u)$ with $t_-(u)$ for each $u\in E$ and obtain a family
$\widehat{\gF}=(\widehat{\pi}:\widehat{\cC} \rightarrow E; s_1, \ldots, s_N;
\xi_1, \ldots, \xi_N)$ of $N$-pointed nodal curves. Let
$\nu :\widetilde{\cC} \widehat{\cC}$ be this identification map. A flat
deformation of this family $\widehat{\gF}$ can be constructed as follows.  For further
details see section 5.3 in \cite{U2}.

For $ 0< \varepsilon \le 1$ put
\begin{eqnarray*}
X_\varepsilon =\{\, P \in \widetilde{\cC}\,|\, |z(P)|<\varepsilon\,\}\\
Y_\varepsilon=\{\, P \in \widetilde{\cC}\,|\, |w(P)|<\varepsilon\,\}
\end{eqnarray*}
and also put $X=X_1$, $Y=Y_1$.  Choose $0<\varepsilon_1<\varepsilon_2<1$ and  choose an
open covering $\{U_\alpha\}_{4 \le \alpha \le s}$ of $\widetilde{\cC} \setminus (X_{\varepsilon_2}
\cup Y_{\varepsilon_2})$ such that
$$
U_\alpha \cap X_{\varepsilon_1} = \emptyset, \quad
U_\alpha \cap Y_{\varepsilon_1} = \emptyset .
$$
  Put
\begin{eqnarray*}
D&=&\{\, q \in \bC\,|\, |q|<1\,\} \\
S_0&=&\{ \, (x,y,q) \in \bC^3\,|\, xy=q, \, |x|<1, \, |y|<1,\, |q|<1\, \} \\
\cS &=& S_0 \times E\\
Z&=& \{\,(P,q) \in \widetilde{\cC}\times D\,|\,
P \in \widetilde{\cC} \setminus(X\cup Y)\, \hbox{\rm , or } P \in X \, \&\,  |z(P)|>|q|\,\} \\
W&=& \{\,(P,q) \in \widetilde{\cC}\times D\,|\,
P \in \widetilde{\cC} \setminus(X\cup Y)\, \hbox{\rm , or } P \in Y \, \&\,  |w(P)|>|q|\,\} .
\end{eqnarray*}
On $Z\cup \cS \cup W$ let us introduce an equivalence relation $\sim$ as follows:.
\begin{enumerate}
\item A point $(P,q) \in Z \cap (X \times D)$ and a point $(x,y,q',u) \in \cS$ are equivalent if
and only if
$$
(x,y,q',u)=\big( z(P),\frac{q}{z(P)}, q, \widetilde{\pi}(P)\big) .
$$
\item A point $(P,q) \in W \cap (Y \times D)$ and a point $(x,y,q',u) \in \cS$ are equivalent if
and only if
$$
(x,y,q',u)=\big( \frac{q}{w(P)}, w(P),q, \widetilde{\pi}(P)\big) .
$$
\item A point $(P,q) \in Z$ and a point $(Q, q') \in W$ are equivalent if and only if
$$
(P, q) = (Q,q').
$$
\end{enumerate}
Now put $\cC = Z \cup \cS\cup W/\sim$. Then, $\cC$ is a complex manifold and there is a
natural holomorphic mapping $\pi : \cC \rightarrow E \times D$. Moreover $\big(\widetilde{\cC}
\setminus (X\cup Y)\big) \times D$ is contained in $\cC$ as an open subset.
Hence we can define holomorphic sections  $s_j$ by
\begin{eqnarray*}
s_j :  E \times D & \rightarrow &\cC \\
\phantom{s_j : }(u,q) &\rightarrow & (s_j(u),q)\in Z.
\end{eqnarray*}
Then it is easy to show that $\gF=(\pi : \cC \rightarrow \cB= E \times D;
s_1, \ldots, s_N; \xi_1, \ldots, \xi_N)$ is a family of $N$-pointed curves with formal coordinates.
It is also easy to show that the family $\pi : \cC \rightarrow \cB=E \times D$ is a versal family
of $N$-pointed curves. Note that each fiber over $q=0$ is a nodal
curve. We identify $E$ with $E\times \{0\}$.

Then we have an exact sequence.
$$
0   \rightarrow  \pi_*\Theta_{\cC/\cB}(*S) \rightarrow
\cL(\gF) \stackrel{\theta}{\rightarrow}  \Theta_\cB(-\log E ) \rightarrow 0.
$$

\begin{Lem}
\label{lem5.8}
We can choose a
local vector field $\vec{m}^{(i)}$  in \eqref{mvector} in such a way that it
can be regarded as a section of  $\cL(\gF)$ over $\cB$ with
$$
\theta(\vec{m}^{(i)}) = \frac{d}{du_i}
$$
\end{Lem}

{\it Proof}. \quad
Since $\frac{d}{du_i}\in H^0(\cB,\Theta_\cB(-\log E))$, we can
choose a $\tau^{(i)}\in \pi_*\Theta_{\cC,\pi}'(*S) $ such that
$\theta(\widetilde{p}(\tau^{(i)})) = \frac{d}{du_i}$, provided we choose $\cB$ small enough.

Let $\{U_\alpha\}_{4 \le \alpha \le s}$ be the open covering of $\widetilde{C} \setminus
\big\{ X_\varepsilon \cup Y_{\varepsilon_2}\big\}$ chosen above.  Put
\begin{eqnarray*}
&&\cU_0 = \{(x,y, u_1, \ldots, u_n, q) \in \cS \,|\, x \ne 0\,\}, \quad
\cU_1 = \{(x,y, u_1, \ldots, u_n, q) \in \cS \,|\, y \ne 0\,\} \\
&&\cU_2 = X \times D, \quad \cU_3= Y \times D, \quad
\cU_\alpha = U_\alpha \times D, \quad 4 \le \alpha \le s.
\end{eqnarray*}
Then, $\{\cU_\alpha\}_{1\le \alpha \le s}$ is an open covering on $\cC$.
Let $(w_\alpha, u_1, \ldots, u_m, q)$ be local coordinates of
$\cU_\alpha$, where
$$
w_0 =y, \quad w_1=x, \quad w_2=x, \quad w_3=y.
$$
In these coordinates our vector field $\tau^{(i)}$ is expressed as
$$
\tau^{(i)} |_{\cU_\alpha}= A_\alpha(w_\alpha ,q,u)\frac{d}{dw_\alpha} + \frac{d}{du_i}
+ B(u, q)q\frac{d}{dq} .
$$
On the other hand the Kodaira-Spencer class $\rho(\frac{d}{du_i})$ is expressed
by a $\check{C}$ech cocycle
$$
\theta_{\alpha \beta} = \frac{\partial f_{\alpha \beta}}{\partial u_i} \frac{d}{dw_\alpha}
$$
where on $\cU_\alpha \cap \cU_\beta \ne \emptyset$ we write
$$
w_\alpha = f_{\alpha\beta}(w_\beta, u_1, \ldots, u_m,  q).
$$
Hence, if we put
$$
\theta_\alpha = - A_\alpha(w_\alpha ,q,u)\frac{d}{dw_\alpha} - B(u, q)q\frac{d}{dq}
$$
then we have
$$
\theta_{\alpha \beta}  = \theta_\beta  - \theta_\alpha.
$$

Note that for $\alpha, \beta \le 4$, $f_{\alpha\beta}$ is independent of $q$.
Therefore on
$\big\{\widetilde{\cC} \setminus \big\{ X_{\varepsilon_2} \cup Y_{\varepsilon_2}\big\}
\big\} \times D$
$$
 A_\alpha(w_\alpha ,0,u)\frac{d}{dw_\alpha} + \frac{d}{du_i}
+ B(u, 0)q\frac{d}{dq} , \quad \alpha \ge 4
$$
define a global meromorphic vector filed $\widetilde{\tau}^{(i)'}$.

On the other hand for $\alpha, \beta \le 3$, $f_{\alpha\beta}$ is independent of
$u_1, \ldots, u_m$.  On $\cU_0 \cap \cU_2$ we have
$$
w_2 = q/w_0.
$$
Hence we have
$$
A_0(w_0,u,q)\frac{d}{dw_0} + \frac{d}{du_i}+ B(u, q)q\frac{d}{dq}
= (-\frac{ w_2}{w_0} A_0(q/w_2, u, q ) + w_2 B(u, q))\frac{d}{dw_2}
+ \frac{d}{du_i}+ B(u, q)q\frac{d}{dq}.
$$
Therefore
$$
A_2(w_2,u,q) = -\frac{ w_2}{w_0} A_0(q/w_2, u, q ) + w_2 B(u, q).
$$
Since $A_0(w_0, u,q)$ and $A_2(w_2, u,q)$ are holomorphic, we conclude
that $A_0(w_1,u,q)$ has the form
$$
A_0(w_0,u,q) = a_0(u,q)w_0,
$$
and on $\cU_0 \cap \cU_2$
$$
A_2(w_2, u,q) =(-a_0(u, q) + B(u,q))w_2,
$$
where $a_0(u,q)$ is holomorphic in $u$ and $q$.
In this way we can show that on $\cU_\alpha$, $\alpha\le 3$, $A_\alpha$ has the form
$$
A_\alpha(w_\alpha, u,q) = a_\alpha(u,q)w_\alpha,
$$
where $a_\alpha(u,q)$ is holomorphic in $u$ and $q$.
Also  the above calculations show that
$$
a_\alpha(u, 0)w_\alpha \frac{d}{dw_\alpha}+ \frac{d}{du_i}+ B(u, 0)q\frac{d}{dq}, \quad
\alpha \le 3
$$
define  a holomorphic vector field on $\cup_{\alpha=0}^3\cU_\alpha$.  This vector field
coincides with $\widetilde{\tau}^{(i)'}$ on the intersection of
$\big\{\widetilde{\cC} \setminus \big\{ X_\varepsilon \cup Y_{\varepsilon_2}\big\}
\big\} \times D$ and $\cup_{\alpha=0}^3\cU_\alpha$. Thus we have a global
meromorphic vector field $\widetilde{\tau}^{(i)}$ on $\cC$. Then, by our construction
if we put $\vec{m}^{(i)}= \widetilde{p}(\widetilde{\tau}^{(i)})$, it is independent
of $q$ and we have $\theta(\vec{m}^{(i)}) = \frac{d}{du_i}$.

Now restrict $\widetilde{\tau}^{(i)}$ to $q=0$
 and denote this restriction $\widetilde{\tau}^{(i)}(0)$.  Then the above
calculation shows that $\widetilde{\tau}^{(i)}(0)I \subset I$ where $I$  is
the ideal defining the double point.
Hence,  $\nu^*(\widetilde{\tau}^{(i)}(0))$ is
a meromorphic vector field on $\widetilde{C}$ and
 $\tilde{\theta}(\tilde{p}(\nu^*\widetilde{\tau}^{(i)}(0))) =
\frac{d}{du_i}$.  Note that
$\tilde{p}(\nu^*\widetilde{\tau}^{(i)}(0))= \vec{m}^{(i)}$.
 $\phantom{x}$ \QED

Now from $\langle \Phi| \in \cVd_{\ab}(\widetilde{\gF})$ we can construct $\langle \Psi(q)|
\in \cVd_{\ab}(\gF)$ using formula \eqref{solution}:
\begin{equation}
\label{solfam}
\langle  \widetilde{\Phi}(q)| u \rangle =
\sum_{p \in \bZ}\Bigl\{ \sum_{d=0}^\infty \sum_{i=1}^{m_d}(-1)^{p+d}
\langle \Phi |v_i(d,p)\otimes v^i(d,-p-1)\otimes u\rangle \Bigr\}q^{d+p(p+1)/2}.
\end{equation}
Then by \eqref{4.2.8}, \eqref{mvector} and \eqref{mtildevector}
\begin{eqnarray*}
&&\langle\widetilde{\Phi}(q)|D(\vec{m}^{(i)})|u\rangle =
-\sum_{j=1}^N \langle \widetilde{\Phi}(q)|
\rho_j(T[m_j^{(i)}])|u\rangle \\ && \phantom{X} = -\sum_{j=1}^N
\sum_{p \in \bZ}\Bigl\{ \sum_{d=0}^\infty
\sum_{i'=1}^{m_d}(-1)^{p+d} \langle \Phi | \rho_j(T[m_j^{(i)}])
v_{i'}(d,p)\otimes v^{i'}(d,-p-1) \otimes u\rangle
\Bigr\}q^{d+p(p+1)/2} \\ && \phantom{XXXX}= -\sum_{p \in
\bZ}\Bigl\{ \sum_{d=0}^\infty \sum_{i'=1}^{m_d}(-1)^{p+d} \langle
\Phi |D(\vec{\widetilde{m}}^{(i)} ) v_{i'}(d,p)\otimes
v^{i'}(d,-p-1) \otimes u\rangle \Bigr\}q^{d+p(p+1)/2} .
\end{eqnarray*}
Let $\omega(q) = \omega(x,y, u.q)dxdy \in H^0(\cC\otimes_{E\times D} \cC,
\omega^{\boxtimes 2}(2 \Delta))$ be a bidifferential with $\res^2(\omega) =1$.
Note that $\omega(0)$ is a bidifferential on $\widetilde{\cC} \times _E
\widetilde{\cC}$.

Then, by \eqref{4.19}  the connection $\nabla^{(\omega(q))}_{\frac{d}{du_i}}$ is defined by
$$
\nabla^{(\omega(q))}_{\frac{d}{du_i}}\langle \widetilde{\Phi}(q)|
= D(\vec{m}^{(i)})\langle \widetilde{\Phi}(q)|  + \frac16 b_{\omega(q)}(\vec{m}^{(i)})
\langle \widetilde{\Phi}(q)|.
$$
Therefore, if we start from $\langle \Psi |=\nabla^{(\omega(0))}_{\frac{d}{du_i}}
\langle \Phi|  \in \cVd_{\ab}(\widetilde{\gF})$ and construct
$\langle \widetilde{\Psi}(q)| \in \cVd_{\ab}(\gF)$ by the sewing procedure, we
have that
$$
\langle \widetilde{\Psi}(q)|=  D(\vec{m}^{(i)})\langle \widetilde{\Phi}(q)|+
\frac16 b_{\omega(0)}(\vec{m}^{(i)})\langle \widetilde{\Phi}(q)|.
$$
Therefore, we have
\begin{equation}
\label{compare}
\nabla^{(\omega(q))}_{\frac{d}{du_i}}\langle \widetilde{\Phi}(q)|
- \langle \widetilde{\Psi}(q)|  =
\frac16\big(b_{\omega(q)}(\vec{m}^{(i)}) - b_{\omega(0)}(\vec{m}^{(i)})\big)
\langle \widetilde{\Phi}(q)| .
\end{equation}

Thus we obtain the following theorem.
\begin{Thm}
\label{thm5.2}
Let $\gF=(\pi : \cC \rightarrow \cB= E \times D;
s_1, \ldots, s_N; \xi_1, \ldots, \xi_N)$ be a family of $N$-pointed
curves with formal coordinates
such that the restriction of the family to $E$ is the  family of
nodal curves
$\gF=(\pi : \widehat{\cC} \rightarrow E ;
s_1, \ldots, s_N; \xi_1, \ldots, \xi_N)$. Let
$$\widetilde{\gF} = ( \widetilde{\pi} : \widetilde{\cC} \rightarrow E ; s_1, \ldots,s_N,
t_+, t_-; \xi_1, \ldots, \xi_N, z, w)
$$
be the family of $(N+2)$-pointed curves obtained by normalization of each fiber of $\widehat{\gF}$.
For $\langle \Phi | \in \cVd_{\ab}(\widetilde{\gF})$ let $\langle \widetilde{\Phi}(q)| \in
\cVd_{\ab}(\gF)$ be obtained by sewing.
Let $\omega(q)\in H^0(\cC\otimes_{E\times D} \cC,
\omega^{\boxtimes 2}(2 \Delta))$ be a bidifferential with $\res^2(\omega) =1$.
Let $\langle \widetilde{\Psi}(q)| \in \cVd_{\ab}(\gF)$ be constructed  from
$\langle \Psi | =\nabla^{(\omega(0))}_{\frac{d}{du_i}}  \langle \Phi|  \in \cVd_{\ab}(\widetilde{\gF})$
by the sewing. Then we have that
$$
\nabla^{(\omega(q))}_{\frac{d}{du_i}}\langle \widetilde{\Phi}(q)|
- \langle \widetilde{\Psi}(q)|  =
\frac16\big(b_{\omega(q)}(\vec{m}^{(i)}) - b_{\omega(0)}(\vec{m}^{(i)})\big)
\langle \widetilde{\Phi}(q)| .
$$
\end{Thm}

\begin{Rmk}
For the non-abelian case the same result holds, provided one replaces the
 coefficient $1/6$ by $c_v/12$.
\end{Rmk}

\section{Formal Coordinates and Preferred Sections}

In this section we shall study the behavior of the ghost vacua
under coordinate changes. We shall also define a preferred section which
is a local holomorphic section of the sheaf of ghost vacua $\cVd_{\ab}(\gF)$
of a family of one-pointed smooth  curves with formal coordinates.

First we recall some basic facts about the relationship between
coordinates change and local vector fields.

We let $\cD$ be the automorphism
group $\Aut \bC((\xi))$ of the field $\bC((\xi))$ of formal Laurent
series as a $\bC$-algebra. The group $\cD$ may be regarded
as the automorphism group
of the ring $\bC[[\xi]]$ of formal power series.
There is a natural isomorphism
\begin{eqnarray*}
\cD & \simeq & \{ \; \sum_{n=0}^\infty a_n \xi^{n+1} \; | \; a_0 \ne 0\; \}\\
 h & \mapsto & h(\xi)
\end{eqnarray*}
where the composition $h \circ g$ of $h$, $g \in \cD$ corresponds to
a formal power series $h(g(\xi))$. In the following we identify $\cD$
with $\{ \; \sum_{n=0}^\infty a_n \xi^{n+1} \; | \; a_0 \ne 0\; \}$.

Put
$$
\cD^p = \{ \;h \in  \cD\; | \; h(\xi) = \xi + a_p\xi^{p+1} + \cdots\;\}
$$
for a positive integer $p$. Then we have a filtration
$$
\cD= \cD^0 \supset \cD^1 \supset \cD^2 \supset \ldots
$$
Put also
\begin{eqnarray*}
\underline{d} &=& \bC[[\xi]]\xi\frac{d}{d\xi} \\
\underline{d}^p &=& \bC[[\xi]]\xi^{p+1}\frac{d}{d\xi} \quad p=0,1,2,\ldots
\end{eqnarray*}
Then, we have a filtration
$$
\underline{d} = \underline{d}^0 \supset \underline{d}^1 \supset \underline{d}^2
\supset \cdots
$$

For any $\underline{l} \in \underline{d}$ and $f(\xi) \in\bC[[\xi]]$ define
$\exp (\underline{l})(f(\xi))$ by
$$
\exp (\underline{l})(f(\xi)) = \sum_{k=0}^\infty \frac{1}{k!}
(\underline{l}^k f(\xi)).
$$
Also put
\begin{eqnarray*}
\cD_+^0 &=& \{ \;h \in  \cD\; | \; h(\xi) = a\xi +
a_1\xi^{2} + \cdots,  \quad a>0 \;\} \\
\underline{d}_+^0 &=& \{ \; l(\xi)\frac{d}{d\xi}\;|\;
l(\xi) = \alpha \xi + \alpha_1 \xi^2 +\cdots, \quad \alpha \in \bR\;\}
\end{eqnarray*}
 Then, we have the following result.
 \begin{Lem}
 \label{lem6.1}
 The exponential map
 \begin{eqnarray*}
 \exp \; : \; \underline{d} & \rightarrow & \cD \\
 \phantom{\exp \; : \;{}} \underline{l} & \mapsto & \exp(\underline{l})
 \end{eqnarray*}
 \label{exponential}
 is surjective. Moreover, the exponential map induces an
 isomorphism
 $$
 \exp \; : \; \underline{d}_+^0 \simeq \cD_+^0 .
 $$
 \end{Lem}

 Since,  for any integer $n$,  we have
 $$
 \exp(2\pi n \sqrt{-1} \xi \frac{d}{d\xi}) = id,
 $$
 the exponential mapping is not injective on $\underline{d}$.

For the energy-momentum tensor $T(z)$ and any element
$\underline{l} \in \cD_+^0$ we define $\exp(T[\underline{l}])$ by
$$ \exp(T[\underline{l}])= \sum_{k=0}^\infty
\frac{1}{k!}T[\underline{l}]^k. $$ Then, $\exp(T[\underline{l}])$
operates on $\cF$ from the left and on $\cFd$ from the right.

By Lemma \ref{exponential},  for any automorphism $h \in \cD_+^0$,
there exist a unique $\underline{l} \in \underline{d}_+^0$ with
$\exp(\underline{l}) = h$. Now for $h \in \cD_+^0$ define
the operator $G[h]$ by
$$
  G[h] = \exp (- T[\underline{l}])
$$
where $\exp (\underline{l}) = h$. Then, $G[h]$ operates on $\cF$
from the left and on $\cFd$ from the right.
Then we have the following important results.
\begin{Thm}
\label{thm6.1}
For any $h \in \cD_+^0$, $f(\xi)d\xi \in \bC((\xi))d\xi$,
$g(\xi) \in \bC((\xi))$   and
$\underline{l} = l(\xi)\frac{d}{d\xi} \in \bC((\xi))\frac{d}{d\xi}$,
we have the following equalities as operators on $\cF$ and $\cFd$.
\begin{eqnarray*}
G[h](\psi[f(\xi)d\xi])G[h]^{-1}&=& \psi[h^*(f(\xi)d\xi)]  =
\psi[f(h(\xi))h'(\xi)d\xi] \\
G[h](\ovpsi[g(\xi)])G[h]^{-1}&=& \ovpsi[h^*(g(\xi))]  =
\ovpsi[g(h(\xi))] \\
G[h_1\circ h_2] &=& G[h_1]G[h_2] \\
G[h]T[\underline{l}]G[h]^{-1} &=& T[\ad(h)(\underline{l})]+
\frac16 \Res_{\xi=0}\big(\{h(\xi);\xi\}l(\xi)d\xi\big) .
\end{eqnarray*}
where $\{f(\xi); \xi\}$ is the Schwarzian derivative.
\end{Thm}

A proof is easily given by applying \eqref{ln} and \eqref{lnbar}.  From this theorem
we infer easily the following proposition.

\begin{Prop}
\label{prop6.1}
For any $h_j \in {\mathcal D}_+^0$, $j=1,2, \ldots, N$ and $N$-pointed curve
$$\frak{X}= (C; Q_1,Q_2, \ldots, Q_n; \xi_1, \xi_2, \ldots, \xi_N)$$
with formal coordinates, put
$$
\frak{X}_{(h)} =(C; Q_1,Q_2, \ldots, Q_N;
h_1(\xi_1), h_2(\xi_2), \ldots, h_N(\xi_N)).
$$
Then, the isomorphism $G[h_1] \widehat{\otimes}\cdots
\widehat{\otimes}G[h_N]$
\begin{eqnarray*}
\cFd_N  & \rightarrow & \cFd_N \\
\langle \phi_1\widehat{\otimes} \cdots \widehat{\otimes}\phi_N|
& \mapsto & \langle \phi_1G[h_1] \widehat{\otimes}\cdots
\widehat{\otimes}\phi_NG[h_N]|
\end{eqnarray*}
induces the canonical isomorphism
$$
G[h_1] \widehat{\otimes}\cdots \widehat{\otimes}G[h_N]:
\cVd_{\ab}(\gX) \rightarrow
\cVd_{\ab}(\frak{X}_{(h)})
$$
\end{Prop}

Let $\gX=(C;Q;\xi)$ be a one-pointed smooth curve of genus $g$
with a formal coordinate. We shall show that if we fix a
symplectic basis $\{\alpha_1, \ldots, \alpha_g, \beta_1, \ldots,
\beta_g\}$ of $H_1(C, \bZ)$, then there is a canonical preferred non-zero
vector $\langle \omega(\gX,\{\alpha, \beta\})| \in
\cVd_{\ab}(\gX)$ which is a refinement of the construction given
in Lemma \ref{lem3.1}. Let us choose a normalized basis
$\{\omega_1, \ldots, \omega_g\}$ of holomorphic one-forms on $C$
which is characterized by
\begin{equation}
\label{betaone}
\int_{\beta_i}\omega_j = \delta_{i j}, \quad 1 \le i,j \le g.
\end{equation}
The matrix $$ \tau = (\tau_{ij}), \quad \tau_{ij} =
\int_{\alpha_i}\omega_j $$ is then called the period matrix of the
curve $C$. Now the numbers $I_n^i$, $n = 1,2, \ldots$, $i=1,\ldots
g$ are defined by
$$ \omega_i = (\sum_{n=1}^\infty I_n^i
\xi^{n-1})d\xi. $$
Note that the numbers $I_n^i$ depend on the
symplectic basis $\{\alpha_1, \ldots, \alpha_g, \beta_1, \ldots,
\beta_g\}$ and the coordinate $\xi$.

For a positive integer $n \ge 1$ let $\omega_Q^{(n)}$ be a meromorphic
one-form on $C$ which has a pole of order $n+1$ at $Q$ and holomorphic elsewhere
such that
\begin{eqnarray}
\label{omegaQ1}
\int_{\alpha_i} \omega_Q^{(n)} &= &
-\frac{2 \pi \sqrt{-1}I_n^i}{n},
\quad \int_{\beta_i} \omega_Q^{(n)} = 0, \quad 1\le i \le g  \\
\label{omegaQ}
\omega_Q^{(n)} &=& \bigl( \frac{1}{\xi^{n+1}} +
\sum_{m=1}^\infty q_{n,m}\xi^{m-1}\bigr)d\xi.
\end{eqnarray}
These conditions uniquely determine $\omega_Q^{(n)}$. Note that the second equality
of \eqref{omegaQ1} and \eqref{omegaQ} imply the first equality of \eqref{omegaQ1}.
The preferred  element $\langle \omega(\gX,\{\alpha, \beta\})| \in \cVd_{\ab}(\gX)$
is defined by
$$
\langle \omega(\gX,\{\alpha, \beta\}) |  = \cdots e(\omega_{g+2} )\wedge e(\omega_{g+1})
\wedge e(\omega_g) \wedge  \cdots \wedge e(\omega_1),
$$
where
$$
\omega_{g+n} = \omega_Q^{(n)}.
$$
For details see Lemma \ref{lem3.1}.
We call $\{\omega_n\}$, $n=1,2,\ldots$ a normalized basis for $\gX$.
Note that the normalized basis depends on the choice of a symplectic
basis of $H_1(C,\bZ)$ and the coordinate $\xi$.

\begin{Thm}
\label{thm6.2}
For $h(\xi) \in \cD_+^0$ put $\gX_h = \{ C;Q;\eta= h(\xi)\}$.
Then
$$
\langle \omega(\gX,\{\alpha, \beta\})| G[h] = \langle \omega(\gX_h,\{\alpha, \beta\})|,
$$
where $G[h] : \cVd(\gX_h) \rightarrow \cVd(\gX)$ is  the canonical isomorphism
given in Proposition \ref{prop6.1}
\end{Thm}

{\it Proof}. \quad  First consider the case in which $h(\xi) = a\xi$
for a positive number $a$.  Put $\eta= a \xi$.
Let $\omega_i$, $i=1,\ldots,g$ and $\omega_{g+n} = \omega_Q^{(n)}$
be chosen for $\gX$ as above. Then
\begin{eqnarray*}
\omega_i &=& \bigl(\sum_{n=1}^\infty a^{-n-1}I_n^i \eta^{n-1}\bigr)d\eta, \\
\omega_Q^{(n)} &=&
\bigl( \frac{a^n}{\eta^{n+1}} +
\sum_{m=1}^\infty  a^{-m-1}q_{n,m}\eta^{m-1}\bigr)d\eta.
\end{eqnarray*}
Put
\begin{eqnarray*}
\widetilde{\omega}_i &=&  \omega_i, \quad 1 \le i \le g \\
\widetilde{\omega}_Q^{(n)} &=& a^{-n}\omega_Q^{(n)}
 = \bigl( \frac{1}{\eta^{n+1}} +
\sum_{m=1}^\infty \frac{q_{n,m}}{ a^{-n -m-1}} \eta^{m-1}\bigr)d\eta.
\end{eqnarray*}
Then, the  element $\langle \gX_h, \{\alpha,\beta\}| \in \cVd_{\ab}(\gX_h)$
is given by
$$
\langle \widetilde{\omega}(\gX,\{\alpha, \beta\}) |  =
\cdots e(\widetilde{\omega}_{g+2} )\wedge e(\widetilde{\omega}_{g+1})
\wedge e(\widetilde{\omega}_g) \wedge  \cdots e(\widetilde{\omega}_1),
$$
where
$$
\widetilde{\omega}_{g+n} = \widetilde{\omega}_Q^{(n)}.
$$
Note that $\langle \omega(\gX,\{\alpha, \beta\})|$,
$\langle \widetilde{\omega}(\gX,\{\alpha, \beta\}) |
\in \cF(g-1)$.

Let $\alpha$ be the
positive number such that $\exp(\alpha) = a$.  Put
$$
\underline{l} = \alpha \xi \frac{d}{d\xi}.
$$
Then we have
$$
T[\underline{l} ] = \alpha L_0, \quad G[h] = \exp(-\alpha L_0).
$$
On $\cFd_d(g-1)$, $G[h]$ operates by the multiplication by
$a^{-(d + g(g-1)/2)}$.
The coefficient of
\begin{equation}
\cdots e^{-l_1-5/2}\wedge e^{-l_1-3/2} \wedge e^{m_1+1/2}
\wedge e^{-l_1+1/2} \wedge \cdots \wedge e^{m_k +1/2} \wedge
\cdots \wedge e^{-3/2}\wedge e^{n_g+1/2}\wedge \cdots \wedge
e^{n_1+1/2}
\end{equation}
in $\langle \omega(\gX,\{\alpha, \beta\})|$ is
$$
A= q_{l_1,m_1+1} \cdots q_{l_k, m_k+1}I_{n_g+1}^g\cdots I_{n_1+1}^1.
$$
The degree $d$ of this term is given by formula (\ref{degform})
$$
d=  \sum_{i=1}^g n_i + \sum_{j=1}^l (m_j + l_j +1)  - \frac{g(g-1)}{2} .
$$
Thus we have that
$$
-d-g(g-1)/2= -\sum_{i=1}^g n_i  - \sum_{j=1}^l (m_j + l_j +1) .
$$
On the other hand the coefficient of $\langle \omega(\gX_h,\{\alpha, \beta\})|$ is
$$
a^{-\sum_{i=1}^g n_i - \sum_{j=1}^l (m_j + l_j +1)}A= a^{-(d+g(g-1)/2)}A.
$$
This implies
$$
\langle \omega(\gX,\{\alpha, \beta\})| G[h] =
\langle \omega(\gX_h,\{\alpha, \beta\})|.
$$
Now let us consider the  case in which $h(\xi)$ is an element of
$\cD^1 \subset \bC((\xi))$.  Let $\{\widetilde{\omega}_n\}$,
$n=1,2,\ldots$ be a normalized basis of $\gX_h$.
Put
$$
\omega_n = h^*(\widetilde{\omega}_n).
$$
Since $h(\xi) = \xi+ a_1\xi^2+\cdots $, the normalized basis of
$\langle \gX|$ is given by $\{\omega_n\}$,  $n=1,2,\ldots$.
For any positive integer $m$ all terms in
$$
\langle \phi_m| = \cdots e^{-m-5/2} \wedge e^{-m -3/2} \wedge e(\omega_{m+g})
\wedge \cdots \wedge e(\omega_1)
$$
are also present in $\langle \omega(\gX,\{\alpha, \beta\})|$ by its definition.
Moreover,  we can express
\begin{eqnarray*}
\langle \phi_m| &=& \langle  -m-1|  \psi[\omega_{m+g}]
\cdots \psi[\omega_1] \\
&=& \langle  -m-1|  \psi[h^*(\widetilde{\omega}_{m+g})]
\cdots \psi[h^*(\widetilde{\omega}_1)] .
\end{eqnarray*}
Then by Theorem \ref{thm6.1}
$$
\langle \phi_m| = \langle  -m-1|
G[h]\psi[\widetilde{\omega}_{m+g}]
\cdots \psi[\widetilde{\omega}_1] G[h]^{-1}.
$$

Note that
$$
\langle -m-1| 
\equiv  \langle -m-1|  \mod \cF_1(-m-1).
$$
It is easy to show that each term appearing in
$$
\langle  \phi_m |  G[h]  - \langle  -m-1|
\psi[\widetilde{\omega}_{m+g}]
\cdots \psi[\widetilde{\omega}_1] G[h]^{-1}.
$$
does not appear in
$$
\langle \phi_{m+1} | = \langle  -m-2| G[h]
\psi[\widetilde{\omega}_{m+1+g}]
\cdots \psi[\widetilde{\omega}_1] G[h]^{-1}.
$$
Taking the limit $m \rightarrow \infty$  we conclude
$$
\langle \omega(\gX,\{\alpha, \beta\})| =
\langle \omega(\gX_h,\{\alpha, \beta\}) |G[h]^{-1}.
$$
{} \QED

Next let us study  the dependence of the section
$\langle \omega(\gX, \{\alpha,\beta\})|$ on the choice of
symplectic basis.
\begin{Thm}
\label{thm6.3}
Let $\{\alpha_1, \ldots \alpha_g,\beta_1, \ldots, \beta_g \}$
and $\{\widetilde{\alpha}_1, \ldots \widetilde{\alpha}_g,
\widetilde{\beta}_1, \ldots, \widetilde{\beta}_g \}$ be symplectic
bases of $H^1(C, \bZ)$ of the non-singular curve $C$. Assume that
$\{\beta_1, \ldots, \beta_g\}$ and
$\{\widetilde{\beta}_1, \ldots, \widetilde{\beta}_g \}$ span the same
Lagrangian sublattice in $H^1(C, \bZ)$. Then
$$
\langle \omega(\gX, \{\alpha,\beta\}) |  = \det U
\langle \omega(\gX, \{\widetilde{\alpha}, \widetilde{\beta}\}) | ,
$$
where $U \in GL(g,\bZ)$ is defined by
$$
\left( \begin{array}{c}
\widetilde{\beta}_1 \\ \vdots \\ \widetilde{\beta}_g
\end{array} \right)
= U
\left( \begin{array}{c}
\beta_1 \\ \vdots \\ \beta_g \end{array} \right) .
$$
\end{Thm}

{\it Proof}. \quad  By the assumption on the $\beta$-cycles we
have that
$$
\left(\begin{array}{c}
\widetilde{\alpha} \\ \widetilde{\beta}
\end{array}\right)
= \left( \begin{array}{cc}
{}^t U^{-1} &B \\ 0 & U
\end{array} \right)
\left(\begin{array}{c}
\alpha \\  \beta
\end{array}\right) ,
$$
where $B$ is a $g \times g$ integral matrix.
For a normalized basis $\{\omega_1, \ldots, \omega_g \}$
of holomorphic one-forms with respect to the
symplectic basis $\{ \alpha,\beta\}$ we have that
$$
\Bigl(\int_{\widetilde{\beta}_i}\omega_j\Bigr) = U.
$$
Hence the normalized basis
of holomorphic one-forms
$\{\widetilde{\omega}_1, \ldots, \widetilde{\omega_g}  \}$ with
respect to $\{\widetilde{\alpha},\widetilde{\beta}\}$ is given by
$$
(\widetilde{\omega}_1, \ldots, \widetilde{\omega_g} ) =
(\omega_1, \ldots, \omega_g )U^{-1}.
$$
Thus
$$
(\widetilde{I}_n^1, \ldots,\widetilde{I}_n^g) =
(I_n^1, \ldots, I_n^g)U^{-1},
$$
for all $n$'s. Therefore,  for the normalized  meromorphic
one-form $\omega_Q^{(n)}$
with respect to the symplectic basis  $\{ \alpha,\beta\}$ we have
that
\begin{eqnarray*}
\left(
\begin{array}{c}
\int_{\widetilde{\alpha}_1}\omega_Q^{(n)} \\
\vdots \\
\int_{\widetilde{\alpha}_g}\omega_Q^{(n)}
\end{array} \right)
&=& {}^t U^{-1} \left(
\begin{array}{c}
\int_{\alpha_1}\omega_Q^{(n)}  \\
\vdots \\
\int_{\widetilde{\alpha}_g}\omega_Q^{(n)}
\end{array} \right) +
B \left(
\begin{array}{c}
\int_{\beta_1}\omega_Q^{(n)}  \\
\vdots \\
\int_{\beta_g}\omega_Q^{(n)}
\end{array} \right) \\
&= &  {}^t U^{-1} \left(
\begin{array}{c}
I_n^1 \\ \vdots  \\ I_n^g \end{array} \right)
= \left(
\begin{array}{c}
\widetilde{I}_n^1 \\ \vdots  \\ \widetilde{I}_n^g \end{array} \right)
\\
\left( \begin{array}{c}
\int_{\widetilde{\beta}_1}\omega_Q^{(n)} \\
\vdots \\
\int_{\widetilde{\beta}_g}\omega_Q^{(n)}
\end{array}  \right)
&=& U \left( \begin{array}{c}
 \int_{\beta_1}\omega_Q^{(n)} \\ \vdots \\
\int_{\beta_g}\omega_Q^{(n)}  \end{array}  \right) =0
\end{eqnarray*}
Thus $\omega_Q^{(n)}$ is also the normalized meromorphic one-form
with respect to the symplectic basis
$\{\widetilde{\alpha}, \widetilde{\beta}\}$. Hence we have that
\begin{equation}
\langle \omega(\gX, \{\alpha,\beta\}) |  = \det U
\langle \omega(\gX, \{\widetilde{\alpha}, \widetilde{\beta}\}) | .
\end{equation}
{}\QED

Note that $\det U= \pm 1$.
Next we shall show that $\langle\omega(\gX)|$ is independent of
the point on the curve $C$. For that
purpose we need  the following lemma.  We use the basis $\{\omega_1, \ldots,
\omega_g\}$ of holomorphic one-forms on the curve $C$ normalized by \eqref{betaone}.
\begin{Lem}
\label{lem6.2}
The numbers $1 = n_1<n_2<\cdots<n_g \le 2g-1$ are the Weierstrass gap values
if and only if
$$
\det \left|
\begin{array}{ccc}
I_{n_1}^1&\cdots & I_{n_g}^g \\
\vdots&\ddots&\vdots \\
I_{n_g}^1 & \cdots & I_{n_g}^g
\end{array}
\right| \ne 0
$$
\end{Lem}

{\it Proof}. \quad If
$$
\det \left|
\begin{array}{ccc}
I_{n_1}^1&\cdots & I_{n_g}^g \\
\vdots&\ddots&\vdots \\
I_{n_g}^1 & \cdots & I_{n_g}^g
\end{array}
\right| =  0,
$$
then there exists a vector $(a_1, \ldots,a_g) \ne (0,\ldots,0)$ such that
$$
(a_1, \ldots,a_g) \left(
\begin{array}{ccc}
I_{n_1}^1&\cdots & I_{n_g}^g \\
\vdots&\ddots&\vdots \\
I_{n_g}^1 & \cdots & I_{n_g}^g
\end{array}
\right) =(0,\ldots, 0) .
$$
Form this we infer
$$
\Res_{\xi=0} \bigl( (a_g \xi^{-n_g} + \cdots + a_1\xi^{-n_1})\omega_i
\bigr) = 0,
\quad i=1, \ldots,g.
$$

This means that there is a meromorphic function
$f \in H^0(C, \cO_C(*Q))$ whose principal part at the point $Q$ is
$a_g \xi^{-n_g} + \cdots + a_1\xi^{-n_1}$.  Thus $n_g$ is not a Weierstrass
gap value.

On the other hand if
$$
\det \left|
\begin{array}{ccc}
I_{n_1}^1&\cdots & I_{n_1}^g \\
\vdots&\ddots&\vdots \\
I_{n_g}^1 & \cdots & I_{n_g}^g
\end{array}
\right| \ne  0,
$$
for any vector $(a_1, \ldots, a_g) \ne 0$,  we have
$$
(a_1, \ldots,a_g) \left(
\begin{array}{ccc}
I_{n_1}^1&\cdots & I_{n_1}^g \\
\vdots&\ddots&\vdots \\
I_{n_g}^1 & \cdots & I_{n_g}^g
\end{array}
\right) \ne (0, \ldots,0) .
$$
This means that
$$
\Res_{\xi=0} \bigl( (a_g \xi^{-n_g} + \cdots + a_1\xi^{-n_1})\omega_i)
\bigr) \ne 0,
$$
for a suitable $1 \le i \le g$. Hence there does not exists a meromorphic
function $f \in H^0(C,\cO_C(*Q))$ whose principal part at $Q$ is
$a_g \xi^{-n_g} + \cdots + a_1\xi^{-n_1}$. \QED

Let $\{P, Q\}$ be two smooth points on the curve $C$ with formal coordinates
$\xi$, $\eta$, respectively. Put
$\gX_0 = (C; P,Q; \xi , \eta)$, $\gX_1= (C; P; \xi)$, $\gX_2=(C;Q;\eta)$.
Then the natural imbeddings
\begin{eqnarray*}
\iota_1&:&   \cF \hookrightarrow \cF_2 \\
  &&  |u\rangle \mapsto |u \rangle \otimes |0\rangle \\
\iota_2&:&   \cF \hookrightarrow \cF_2 \\
  &&  |u\rangle \mapsto |0\rangle  \otimes  |u \rangle\\
\end{eqnarray*}
induce canonical isomorphisms
$$
\begin{array}{ccccc}
&&\cVd_{\ab}(\gX_0)& &\\
&&  &&   \\
& {}^{\iota_1^*} \swarrow&& \searrow^{\iota_2^*} & \\
&&  &&   \\
&\cVd_{\ab}(\gX_1) && \cVd_{\ab}(\gX_2) &
\end{array}
$$ by Theorem \ref{thm3.4}.

\begin{Thm}
\label{thm6.4}
Under the above notation we have
$$
\iota_2^*\circ (\iota^*_1)^{-1}(\langle \omega(\gX_1,\{\alpha, \beta\})| )
= \langle\omega(\gX_2,\{\alpha, \beta\})|.
$$
\end{Thm}

{\it Proof}. \quad Put $\langle \widetilde{\Phi}|
= (\iota^*_1)^{-1}(\langle \omega(\gX_1,\{\alpha, \beta\})| )$.
Then, by Theorem \ref{thm3.4} we have
$$
\langle \widetilde{\Phi}| u \otimes 0\rangle =
\langle \omega(\gX_1,\{\alpha, \beta\}) | u \rangle
$$
for any $|u\rangle \in \cF$. Since the space of  ghost vacua is one-dimensional,
we have that
\begin{equation}
\label{cc}
\langle \widetilde{\Phi}|0 \otimes  u \rangle =
c \langle \omega(\gX_2,\{\alpha, \beta\})|  u \rangle
\end{equation}
for a constant independent  of $|u\rangle  \in \cF$. To determine $c$ it is enough to
choose a  $|u\rangle $ such that we can calculate both sides of the equality (\ref{cc}).
Let $1=n_1<n_2< \cdots < n_g$ be the gap values at the point $Q$ and put
$$
|u \rangle = \overline{e}_{n_g -1/2}\wedge \cdots
\overline{e}_{n_1 -1/2}\wedge |-1\rangle .
$$
Let $f_j \in H^0(C, \cO_C(*(P+Q))$, $j=1,2,\ldots, g$  be chosen in such a way that
$f_j$ has the Laurent expansion
$$
f_j = \eta^{-n_j}+a_0^{(j)} + a_1^{(j)}\eta +\cdots
$$
at $Q$. Also let $\omega_{Q,P}$ be a meromorphic one-form having poles of
order one at $Q$ with residue 1 and at $P$ with residue $-1$ and else holomorphic.
Then we have that
\begin{eqnarray*}
\langle \widetilde{\Phi}|0 \otimes  u \rangle &= &
\langle \widetilde{\Phi}|0 \otimes  \rho_Q (\ovpsi[f_g])
\rho_Q (\ovpsi[f_{g-1}])\cdots \rho_Q (\ovpsi[f_1]) |-1\rangle \\
&= &\langle \widetilde{\Phi}|0 \otimes  \rho_Q (\ovpsi[f_g])
\rho_Q (\ovpsi[f_{g-1}])\cdots \rho_Q (\ovpsi[f_1])
\rho_Q(\psi[\omega_{Q,P}])0\rangle \\
&= &(-1)^g
\langle \widetilde{\Phi}| 0 \otimes \rho_Q(\psi[\omega_{Q,P}]) \rho_Q (\ovpsi[f_g])
\rho_Q (\ovpsi[f_{g-1}])\cdots \rho_Q (\ovpsi[f_1])
0\rangle \\
&=& (-1)^{g+1} \langle \widetilde{\Phi}| \rho_P (\ovpsi[\omega_{Q,P}])
0 \otimes  \rho_Q (\ovpsi[f_g])
\rho_Q (\ovpsi[f_{g-1}])\cdots \rho_Q (\ovpsi[f_1]) 0\rangle \\
&=& (-1)^g \langle \widetilde{\Phi}| -1 \otimes  \rho_Q (\ovpsi[f_g])
\rho_Q (\ovpsi[f_{g-1}])\cdots \rho_Q (\ovpsi[f_1]) 0\rangle \\
&=& (-1)^{g} \langle \widetilde{\Phi}| \rho_P (\ovpsi[f_g]) |-1 \otimes
\rho_Q (\ovpsi[f_{g-1}])\cdots \rho_Q (\ovpsi[f_1]) 0\rangle \\
&=& (-1)^{g+g(g-1)/2} \langle \widetilde{\Phi}| \rho_P (\ovpsi[f_1])  \cdots
\rho_P (\ovpsi[f_{g-1}])\rho_P (\ovpsi[f_g]) |-1 \otimes  0\rangle \\
&=& (-1)^{g} \langle \widetilde{\Phi}| \rho_P (\ovpsi[f_g])
\rho_P (\ovpsi[f_{g-1}]) \cdots \rho_P (\ovpsi[f_1]) |-1 \otimes  0\rangle \\
&=& (-1)^g \langle \omega(\gX_1,\{\alpha, \beta\}) | \rho_P (\ovpsi[f_g])
\rho_P (\ovpsi[f_{g-1}]) \cdots \rho_P (\ovpsi[f_1])| -1 \rangle
\end{eqnarray*}
The last term can be expressed in a matrix form
$$
(-1)^g \left|
\begin{array}{cccc}
\res_P(f_g\omega_1)&\res_P(f_g\omega_{2})&\cdots & \res_P(f_g\omega_g) \\
\res_P(f_{g-1}\omega_1)&\res_P(f_{g-1}\omega_{2})&\cdots &
 \res_P(f_{g-1}\omega_g) \\
\vdots&\vdots& \ddots & \vdots \\
\res_P(f_{2}\omega_1)&\res_P(f_{2}\omega_{2})&\cdots &
 \res_P(f_{2}\omega_g) \\
\res_P(f_{1}\omega_1)&\res_P(f_{1}\omega_{2})&\cdots &
 \res_P(f_{1}\omega_g) \\
\end{array}
\right|
$$
But $f_j\omega_k$ has poles only at $P$ and $Q$, thus
$$
\res_P(f_j\omega_k) = -\res_Q(f_j\omega_k)= - I_{n_j}^k(Q)
$$
where
$$
\omega_k = \sum_{n=1}^\infty I_n^k(Q) \eta^{n-1} d\eta.
$$
Hence the above determinant  is
$$\left|
\begin{array}{ccc}
I_{n_g}^1(Q)&\cdots & I_{n_g}^g(Q) \\ \vdots&\ddots&\vdots \\
I_{n_1}^1(Q) & \cdots & I_{n_1}^g(Q)
\end{array} \right|= \langle \omega(\gX_2,\{\alpha, \beta\})|
\overline{e}_{n_g -1/2}\wedge \cdots
\wedge\overline{e}_{n_1 -1/2}\wedge |-1\rangle .
$$
Hence we have that
$$
\langle \widetilde{\Phi}| 0 \otimes \overline{e}_{n_g -1/2}\wedge \cdots
\wedge\overline{e}_{n_1 -1/2}\wedge| -1 \rangle =
\langle \omega(\gX_2,\{\alpha, \beta\}) |  \overline{e}_{n_g -1/2}\wedge \cdots
\wedge\overline{e}_{n_1 -1/2}\wedge | -1 \rangle.
$$
This shows  $c=1$. \QED

Let us consider a curve $C$ with a node $P$. Let $\widetilde{C}$ be
the curve obtained by resolving the singularity at $P$ and let
$\pi : \widetilde{C} \rightarrow C$ be the natural holomorphic
mapping. Then $\pi^{-1}(P)$ consists of two points $P_+$ and $P_-$.
Assume that $\widetilde{C}$ is connected.

Let
$$
\gX=(C;Q;\xi)
$$
be a one-pointed curve with formal coordinates and we let
$$
\widetilde{\gX} = (\widetilde{C};P_+,P_-,Q;z,w, \xi)
$$
be the associated 3-pointed curve with formal coordinates.
Then  Theorem \ref{thm3.5} says that the natural isomorphism
$$
|0_{+,-}\rangle  \otimes \cF \cong \cF
$$
induces the natural isomorphism
$$
\iota_{+,-}^* : \cVd_{\ab}(\widetilde{\gX})   \cong   \cVd_{\ab}(\gX)  .
$$
where
$$
|0_{+,-}\rangle = |0\rangle \otimes |-1\rangle
 - |-1\rangle \otimes |0 \rangle .
$$
We can define $\langle \omega(\gX,\{\alpha, \beta\})|$ similar
to the non-singular case, by choosing a basis $\{\alpha,
\beta\}=\{\alpha_1, \ldots, \alpha_{g-1},\alpha_g, \beta_1,
\ldots, \beta_{g-1}\}$ of $H_1(C,\bZ)$, in such a way that
$\alpha_1$, $\alpha_2, \ldots,\alpha_{g-1}$ and $\beta_1$,
$\beta_2, \ldots, \beta_{g-1}$ is the image of a symplectic basis
of $H_1(\widetilde{C},\bZ)$ under natural map to $H_1(C,\bZ)$ and
$\alpha_g$ corresponds to the invariant cycle of a flat
deformation of the curve $C$. Then we can choose a basis
$\{\omega_1, \ldots, \omega_{g-1}, \omega_g, \omega_{g+1},
\omega_{g+2}, \ldots \}$ of $H^0(C, \omega(*Q))$ such that
$\{\pi^*\omega_1, \ldots, \pi^*\omega_{g-1}, \pi^*\omega_{g+1},
\pi^*\omega_{g+2}, \ldots \}$ is a normalized basis of
$H^0(\widetilde{C}, \omega_{\widetilde{C}}(*Q))$ as in
(\ref{betaone}), (\ref{omegaQ1}) and (\ref{omegaQ}) where we put $$
\pi^*\omega_{g+n} = \omega_Q^{(n)}, \quad n=1,2, \ldots, $$ and
$\pi^*\omega_g$ is a meromorphic one-form on $\widetilde{C}$ which
has poles of order one at $P_+$ and $P_-$ with residue $-1$ and 1,
respectively, is holomorphic outside $P_\pm$ and $$
\int_{P_+}^{P_-}\pi^*\omega_g = 1. $$
Then put $$ \langle
\omega(\gX,\{\alpha, \beta\})| = \langle \cdots \wedge
e(\omega_{m} )\wedge \cdots \wedge e(\omega_2)\wedge e(\omega_1) |. $$
The proof of Lemma \ref{lem3.1} applies also in this case and shows that
$\langle \omega(\gX,\{\alpha, \beta\})|$ is an element of
$\cVd_{\ab}(\gX)$. Let
$$ \widehat{\gX} = (\widetilde{C}, Q;
\xi) . $$
Then, by applying Theorem \ref{thm3.4} at the points
$P_\pm$ we have a canonical isomorphism
$$ \iota^* :
\cVd_{\ab}(\widetilde{\gX}) \cong \cVd_{\ab}(\widehat{\gX}) .
$$
\begin{Thm}
\label{thm6.5} Under the above assumptions and notation we have
that
$$
 \iota^* \circ  (\iota_{+,-}^*)^{-1}  (\langle \omega(\gX,\{\alpha, \beta\}) |) =
 (-1)^g  \langle \omega(\widehat{\gX},\{\widehat{\alpha}, \widehat{\beta}\})|
$$
where $\{\widehat{\alpha}, \widehat{\beta}\}= \{\alpha_1, \ldots, \alpha_{g-1},
\beta_1, \ldots, \beta_{g-1}\}$.
\end{Thm}

{\it Proof}. \quad Put
$$
\langle \widetilde{\Phi}|= (\iota_{+,-}^* )^{-1}
(\langle \omega(\gX,\{\alpha, \beta\})|).
$$
By the definition of $\iota^*$ we have that
$$
\langle \widetilde{\Phi}|0\otimes 0\otimes u\rangle =
\langle \iota^*\circ (\iota_{+,-}^* )^{-1}(\langle \omega(\gX,\{\alpha, \beta\})|)|u\rangle.
$$
Choose a meromorphic function
$f \in H^0(\widetilde{C}, \cO_{\widetilde{C}}(*Q))$ such that $f(P_+)=-1$
and $f(P_-)=0$. Then by (\ref{def00}) (using the notation $\langle \omega(\gX)| =
\langle \omega(\gX,\{\alpha, \beta\})| $), we have that
$$
\langle \widetilde{\Phi}|0\otimes 0\otimes u\rangle =
\langle \omega(\gX)|    \rho_Q(\ovpsi[f]) u\rangle
$$
for any $ |u \rangle \in \cF$.  The proof of lemma \ref{lem3.1} show that
for a meromorphic form $\omega_j$ on $C$ we have that
$$
\wedge e(\omega_j)\ovpsi[f] = \res_Q(f\omega_j )=  \res_Q(f\pi^*\omega_j ) =
\wedge e(\pi^*\omega_j)\ovpsi[f] .
$$
For $j \ne g$, the meromorphic forms $\pi^*\omega_j$ have only poles at
$Q$. Hence
$$
\wedge e(\omega_j)\ovpsi[f] =0, \quad j \ne g.
$$
On the other hand, since $\pi^*\omega_g$ has a pole of order one at $P_+$ and
$P_-$ with residues $-1$ and $1$, respectively, and is holomorphic outside
the points $P_\pm$ and $Q$  and $f(P_+)=-1$, $f(P_-)=0$ we have
that
$$
\res_Q(f\pi^*\omega_g)= -\res_{P_+}(f \pi^*\omega_g) = -1.
$$
Thus we conclude that
\begin{eqnarray*}
\langle \widetilde{\Phi}|0\otimes 0\otimes u\rangle &=&
\langle \omega(\gX)|   \rho_Q(\ovpsi[f]) u\rangle =
\langle \omega(\gX)\rho_Q(\ovpsi[f]) |   u\rangle \\
&=&(-1)^{g-1} \res_Q(f\omega_g)\langle \cdots\wedge e(\omega_{g+1} )\wedge
e(\omega_{g-1}) \wedge \cdots \wedge e(\omega_1) |u\rangle \\
&=& (-1)^{g}\langle \omega(\widehat{\gX},\{\widehat{\alpha}, \widehat{\beta}\})| u\rangle.
\end{eqnarray*}
{}\QED

For an $N$-pointed curve $\gX=(C;Q_1, \ldots, Q_N;\xi_1,\ldots,\xi_N)$ we
define the preferred element $\langle \omega(\gX, \{\alpha,\beta\})|$ by
\begin{equation}
\label{Nsection}
\langle \omega(\gX, \{\alpha,\beta\})|= {\iota^*}^{-1}
\langle \omega(\gX_1, \{\alpha,\beta\})|
\end{equation}
where $\gX_1=(C; Q_1, \xi_1)$ and $\iota^* : \cVd_{\ab}(\gX)
\rightarrow \cVd_{\ab}(\gX_1)$ is defined by applying  Theorem \ref{thm3.4}
several times.
Note by Theorem \ref{thm6.4} that if we choose to use the pair $(Q_j, \xi_j)$
instead of $(Q_1, \xi_1)$ in defining $\gX_1$, we get the same preferred
element.

For an $N$-pointed nodal curve $\gX=(C;Q_1, \ldots, Q_N;\xi_1,\ldots,\xi_N)$
with formal coordinates such that the normalization $\widetilde{C}$ of $C$
is an irreducible curve, we can define the preferred element by
generalizing the discussion just before Theorem \ref{thm6.5} in
the obvious way.

We can however not apply the same
method to define the preferred element for one-pointed nodal
curves $(C;Q, \xi)$ whose normalization $\widetilde{C}$ is disconnected.
This is simply because a holomorphic one-form with support
on one component will have zero Taylor expansion at all points of the curve not
contained in that component.

Therefore, we need to find another definition. For simplicity,  in the
following,
we shall  consider a two-pointed nodal curve $\gX=(C; Q_1,Q_2; \xi_1,\xi_2)$
with formal coordinates such that  the normalized curve $\widetilde{C}$ has
two connected components $C_1$ and $C_2$. Assume that $Q_1$, $P_+ \in
C_1$ and
$Q_2$, $P_- \in C_2$. Put $\gX_1=(C_1; Q_1, \xi_1)$, $\gX_2=(C_2; Q_2, \xi_2)$
$\widetilde{\gX}_1=(C_1; Q_1, P_+,\xi_1, z)$,
$\widetilde{\gX}_2=(C_2; Q_2, P_-,\xi_2, w)$, $\widetilde{\gX}
= (C_1\cup C_2; Q_1,Q_2,P_+,P_-, \xi_1,\xi_2, z,w)$. Then by Theorem \ref{thm3.4}
and Theorem \ref{thm3.5} we have isomorphisms
\begin{eqnarray*}
\iota_j^* :\cVd_{\ab}(\widetilde{\gX}_j) &\cong &\cVd_{\ab} (\gX_j), \quad j=1,2 \\
\iota_{+,-}^* : \cVd_{\ab}(\widetilde{\gX}) &\cong & \cVd_{\ab} (\gX).
\end{eqnarray*}
Moreover, by Proposition \ref{prop3.1} we have
$$
\cVd_{\ab}(\widetilde{\gX}) = \cVd_{\ab}(\widetilde{\gX}_1)\otimes
\cVd_{\ab}(\widetilde{\gX}_2).
$$
Choose symplectic bases $\{\alpha^{(i)}, \beta^{(i)}\}$ of $H_1(C_i, \bZ)$.
Then, $\{\alpha, \beta\}= \{\alpha^{(1)}, \alpha^{(2)},
\beta^{(1)}, \beta^{(2)}\}$ is a symplectic basis of $H_1(C,\bZ)$.
Put
\begin{eqnarray*}
\langle \Phi_1| &= &
{\iota_1^*}^{-1}(\langle\omega(\gX_1, \{\alpha^{(1)}, \beta^{(1)}\})|),\\
\langle \Phi_2| &= &
{\iota_2^*}^{-1}(\langle\omega(\gX_2, \{\alpha^{(2)}, \beta^{(2)}\})|), \\
 \langle \Phi | &=& \langle \Phi_1| \otimes \langle \Phi_2|.
\end{eqnarray*}
Finally define
\begin{equation}
\label{dfn6.7}
\langle \omega(\gX, \{\alpha,\beta\})|=\iota_{+,-}^*(\langle \Phi|).
\end{equation}
For a general $N$-pointed nodal curve $\gX=(C;Q_1, \ldots, Q_N;\xi_1,\ldots,\xi_N)$
we define the preferred element $\langle \omega(\gX, \{\alpha,\beta\})|$
in a similar way.
The following lemma plays an important role in characterizing the preferred
section.
\begin{Lem}
\label{lem6.6}
Under the same notation as above, let the  genera of  $C_1$ and $C_2$ be
$g_1$ and $g_2$, respectively. Let $1=m_1<m_2<\cdots<m_{g_1}\le 2g_1-1$
and $1<n_1<n_2<\cdots<n_{g_2}\le 2g_2-1$ be the Weierstrass gap values of $C_1$
at $Q_1$ and $C_2$ at $Q_2$, respectively. Put
$$
| u_1\rangle = \overline{e}_{m_{g_1} -1/2}\wedge \cdots \wedge
\overline{e}_{m_1 -1/2}\wedge |-1\rangle,
\quad  | u_2\rangle =\overline{e}_{n_{g_2} -1/2}\wedge \cdots \wedge
\overline{e}_{n_1 -1/2}\wedge |0\rangle.
$$
Then we have that
$$
\langle (\iota_{+,-}^*)^{-1}\omega(\gX, \{\alpha,\beta\})|u_1\otimes u_2 \rangle
= \left|
\begin{array}{ccc}
I_{m_{g_1}}^1(Q_1)&\cdots & I_{m_{g_1}}^{g_1}(Q_1) \\
\vdots&\ddots&\vdots \\
I_{m_{1}}^1(Q_1) & \cdots & I_{m_{1}}^{g_1}(Q_1)
\end{array}
\right| \cdot  \left|
\begin{array}{ccc}
I_{n_{g_2}}^{g_1+1}(Q_2)&\cdots & I_{n_{g_2}}^{g_1+g_2}(Q_2)\\
\vdots&\ddots&\vdots \\
I_{n_{1}}^{g_1+1}(Q_2) & \cdots & I_{n_{1}}^{g_1+g_2}(Q_2)
\end{array}
\right| \ne 0
$$
where $\{\omega_1, \ldots , \omega_{g_1}\}$ and
$\{\omega_{g_1+1}, \ldots , \omega_{g_1+g_2}\}$ are the normalized
bases of holomorphic one-forms on $C_1$and  $C_2$, respectively, and
$$
\big(\sum_{m=1}^\infty I_m^k(Q_1)\xi_1^{m-1}\big)d\xi_1, \quad 1\le k \le g_1
$$
is the Taylor expansion of $\omega_k$ at $Q_1$ and
$$
\big(\sum_{n=1}^\infty I_n^{g_1+i} (Q_2)\xi_2^{n-1}\big)d\xi_2, \quad 1\le i \le g_2.
$$
is the Taylor expansion of $\omega_{g_1+i}$ at $Q_2$.
\end{Lem}

{\it Proof}. \quad
Let $\omega_1$ be a meromorphic one-form on $C_1$ which is holomorphic
outside $P_+$ and $Q_1$ such that at $P_+$ and $Q_1$, $\omega_1$
has the Laurent expansions:
\begin{eqnarray*}
\omega_{1,+}&=& \frac{dz}{z} + \hbox{\rm holomorphic} \\
\omega_{1, Q_1}&=& -\frac{d\xi_1}{\xi_1} + \hbox{\rm holomorphic}.
\end{eqnarray*}
Similarly let $\omega_2$ be a meromorphic one-form on $C_2$ which is holomorphic
outside $P_-$ and $Q_2$ such that at $P_-$ and $Q_2$ $\omega_2$
has the Laurent expansions:
\begin{eqnarray*}
\omega_{2, -}&=& -\frac{dw}{w} + \hbox{\rm holomorphic} \\
\omega_{2,Q_2}&=& \frac{d\xi_2}{\xi_2} + \hbox{\rm holomorphic}.
\end{eqnarray*}
Then, we have that
$$
\psi[\omega_{1,+}]|0\rangle = | -1\rangle, \quad
-\psi[\omega_{1,-}]|0\rangle = | -1\rangle.
$$
and therefore, we have that
\begin{eqnarray*}
\langle (\iota_{+,-}^*)^{-1}\omega(\gX, \{\alpha,\beta\})|u_1\otimes u_2 \rangle &=&
\langle \Phi |0_{+,-} \otimes u_1\otimes u_2 \rangle  \\
&=& \langle \Phi |0 \otimes -1 \otimes u_1\otimes u_2 \rangle
- \langle \Phi |-1 \otimes 0 \otimes u_1\otimes u_2 \rangle \\
&=& \langle \Phi_1 | 0 \otimes  u_1 \rangle
\langle \Phi_2 | -1 \otimes  u_2 \rangle -
\langle \Phi_1 | -1 \otimes  u_1 \rangle
\langle \Phi_2 | 0 \otimes  u_2 \rangle\\
&=& \langle \omega(\gX_1, \{\alpha^{(1)},\beta^{(1)}\}) |  u_1 \rangle
\langle \Phi_2 |\psi[\omega_{1,+}]|0 \rangle \otimes  u_2 \rangle \\
&& +
\langle \Phi_1 | \psi[\omega_{1,-}]|0 \rangle  \otimes  u_1 \rangle
\langle  \omega(\gX_2, \{\alpha^{(2)},\beta^{(2)}\}) |  u_2 \rangle\\
&=&
\langle \omega(\gX_1, \{\alpha^{(1)},\beta^{(1)}\}) |  u_1 \rangle
 \langle  \omega(\gX_2, \{\alpha^{(2)},\beta^{(2)}\}) |
 \psi[\omega_{1,Q_2}]u_2 \rangle\\
&& +
\langle \omega(\gX_1, \{\alpha^{(1)},\beta^{(1)}\}) |
 \psi[\omega_{1,Q_1}]u_1 \rangle
 \langle  \omega(\gX_2, \{\alpha^{(2)},\beta^{(2)}\}) |
 u_2 \rangle
\end{eqnarray*}

Since
$$
\psi[\omega_{1,Q_2}]|u_2\rangle =
\overline{e}_{n_{g_2} -1/2}\wedge \cdots
\wedge\overline{e}_{n_1 -1/2}\wedge |-1\rangle + * \wedge |0\rangle
$$
and $\langle  \omega(\gX_2, \{\alpha^{(2)},\beta^{(2)}\}) | $ does not
contains the term $e^{-1/2}$, we have that
\begin{eqnarray*}
\langle  \omega(\gX_2, \{\alpha^{(2)},\beta^{(2)}\}) |
 \psi[\omega_{1,Q_2}]u_2 \rangle
 &= & \langle  \omega(\gX_2, \{\alpha^{(2)},\beta^{(2)}\}) |
 \overline{e}_{n_{g_2} -1/2}\wedge \cdots
\overline{e}_{n_1 -1/2}\wedge |-1\rangle \\
&=& \left|
\begin{array}{ccc}
I_{n_{g_2}}^{g_1+1}(Q_2)&\cdots & I_{n_{g_2}}^{g_1+g_2}(Q_2)\\
\vdots&\ddots&\vdots \\
I_{n_{1}}^{g_1+1}(Q_2) & \cdots & I_{n_{1}}^{g_1+g_2}(Q_2)
\end{array}
\right| \ne 0
 \end{eqnarray*}
Since by  the same reason $\langle \omega(\gX_1, \{\alpha^{(1)},\beta^{(1)}\}) |
 \psi[\omega_{1,Q_1}]u_1 \rangle
 \langle  \omega(\gX_2, \{\alpha^{(2)},\beta^{(2)}\}) |
 u_2 \rangle=0$ and
$$
\langle \omega(\gX_1, \{\alpha^{(1)},\beta^{(1)}\}) |  u_1 \rangle
= \left|
\begin{array}{ccc}
I_{m_{g_1}}^1(Q_1)&\cdots & I_{m_{g_1}}^{g_1}(Q_1) \\
\vdots&\ddots&\vdots \\
I_{m_{1}}^1(Q_1) & \cdots & I_{m_{1}}^{g_1}(Q_1)
\end{array}
\right|\ne 0
$$
we obtain the desired result. \QED

Let $\gF=(\pi : \cC\rightarrow \cB, s_1, \ldots, s_N, \xi_1, \ldots,
\xi_N)$ be a family of $N$-pointed smooth curves.
For any point $b \in \cB$ there exists an open neighbourhood $U_b$ such that
$\pi^{-1}(U_b)$ is topologically  trivial so that we can choose smoothly
varying symplectic bases
$$
\{\alpha_1(t), \ldots, \alpha_g(t), \beta_1(t),
\ldots, \beta_g(t)\}, \quad t \in \cB.
$$
Then we can define
$\langle \omega(\gX_t, \{\alpha(t), \beta(t)\})|$ where
$\gX_t=(\pi^{-1}(t), s_1(t), \xi_1)$.

\begin{Thm}
\label{thm6.6}
The section $\langle \omega(\gX_t, \{\alpha(t), \beta(t)\})|$ is a holomorphic section of
$\cVd_{\ab}(\gF)$ over $U_b$.
\end{Thm}

{\it Proof}. \quad Put
$\gF' = (\pi : \cC\rightarrow \cB, s_1, \xi_1)$. Then, by Theorem \ref{thm3.4}
it is easy to show that we get an isomorphism
$\cVd_{\ab}(\gF) \cong \cVd_{\ab}(\gF')$ by the propagation of vacua construction.
Hence it is enough to
show that $\langle \omega(\gX_t, \{\alpha(t), \beta(t)\})|$ depends holomorphically
on $t$. This follows, since $\omega_{s_1(t)}^{(n)}$ varies holomorphically in $t$. {}\QED

Finally, we analyze the preferred section for families of
deformations of nodal curves.

Let $C_0$ be a complete curve with only one ordinary double point $P$
such that $C_0 \setminus \{P\}$ is non-singular.
Let $Q_1$, $Q_2, \ldots, Q_N$ be distinct non-singular points on $C_0$.
Let $\nu : \widetilde{C}_0
\rightarrow C_0$ be the normalization of the singular curve. Put $\{P_+, P_-\}
= \nu^{-1}(P)$. Note that the normalization might or might not be
connected but assume that each component contains at least one $Q_j$.
Let $\gF = (\pi : \cC \rightarrow D;
\sigma_1, \ldots, \sigma_N; \xi_1, \ldots, \xi_N)$
be the family constructed from $(C_0,Q_1,\ldots,Q_N,\xi_1,\ldots \xi_N)$
as described right before Lemma
\ref{lem5.7}. Suppose we now have a continuous basis
$\{\alpha_i(t),\beta_i(t)\}$ of $H_1(\pi^{-1}(t),\bZ)$, $t\in (0,1)\subset D$,
such that we get a well defined limit as $t$ goes to zero, which
gives a symplectic basis, say
$\{\alpha_1(0),\ldots,
\alpha_{g-1}(0),\alpha_{g}(0),\beta_{1}(0),\ldots,\beta_{g-1}(0)\}$ of $H_1(C_0,\bZ)$ as
described above for nodal curves and $\beta_g(0) = 0$. Let $\gX_t=(\pi^{-1}(t), s_1(t), \ldots s_N(t),
 \xi_1, \ldots \xi_N)$.

\begin{Thm}
\label{thm6.7}
We have that
$$\langle \omega(\gX_0, \{\alpha(0), \beta(0)\})| =
\lim_{t\to 0} \langle \omega(\gX_t, \{\alpha(t), \beta(t)\})|.$$
\end{Thm}

{\it Proof}. \quad  Assume the curve $\widetilde{C}_0$ is connected. It is then enough to prove
the theorem when $N=1$. Note by [F, Proposition 3.7] that if
$\{\omega_1,\ldots, \omega_g\}$ is a normalized
basis of holomorphic one-forms on $C_0$, then a normalized  bases of holomorphic one-forms
for
the family is of the form:
\begin{eqnarray*}
\omega_i(x, t) &=& \omega_i(x) + \frac14t(\omega_i(P_+) -\omega_i(P_-))(\omega(x, P_+) -
\omega(x, P_-))+ O(t^2), \quad 1\le i \le g-1\\
\omega_{g}(x,t) &=&\omega_{P_+-P_-}(x) + t u_g(x) + O(t^2).
\end{eqnarray*}
Here $\{\omega_1(x), \ldots, \omega_{g-1}(x) \}$ is a  normalized basis of
holomorphic one-forms of $\widetilde{C}_0$, and
$\omega_i(P_+)$ is the number  $f_i(0)$ where in a neighbourhood of
$P_+$,
$\omega_i$ is expressed as $f_i(z)dz$. The number $\omega_i(P_-)$  is defined
similarly. Moreover,
$\omega(x,y)$ is the normalized bidifferential of the curve $\widetilde{C}_0$
and in a neighbourhood of $(x, P_+)$ if we express $\omega(x,y) =
f(x,z)dxdz$ then $\omega(x, P_+) $ is defined as $f(x,0)dx$. The one-form
$\omega(x, P_-)$ defined similarly.  The form
$\omega_{P_+-P_-}(x)$ is a meromorphic one-form of $\widetilde{C}_0$
which has a pole of order one at $P_+$ with residue 1,  pole of order one at $P_-$ with
residue $-1$ and holomorphic outside $P_+$ and $P_-$, and $u_g$ is a meromorphic one-form
on $\widetilde{C}_0$ which has only poles at $P_\pm$ of order three. Finally, the
expression $O(t^2)$ means that it is a holomorphic one-form on
$\widetilde{C}_0\setminus\{P_+, P_-\}$
and $\lim_{t \to 0}\frac{O(t^2)}{t^2}$ is a holomorphic one-form on
$\widetilde{C}_0\setminus\{P_+, P_-\}$.


Let
$$
\omega_i(x,t) = (\sum_{n=1}^\infty I_n^i(t) \xi^{n-1})d\xi
$$
for $i=1, \ldots g$.

The  meromorphic one-forms $\omega_{g+n}(t)$ on $C_t$, as defined by (\ref{betaone}) and
(\ref{omegaQ1}) has a
pole only at $Q_1=s_1(t)$ of order $n+1$ with Laurent expansion at $s_1(t)$ of the form
$$
\frac{d\xi_1}{\xi_1^{n+1} }+ \hbox{\rm holomorphic}
$$
and satisfies
$$
\int_{\alpha_i(t)}\omega_{g+n}(t) = - \frac{2\pi \sqrt{-1}I_n^i(t)}{n}
\mbox{, }\int_{\beta_i(t)} \omega_{g+n}(t) = 0, \quad 1\le i\le g, \quad n \ge
0.
$$

Then,
$$\langle \omega(\gX_t, \{\alpha(t), \beta(t)\})| =
\langle \cdots \wedge e(\omega_{m}(t) )\wedge \cdots \wedge
e(\omega_2(t))\wedge e(\omega_1(t)) |.
$$
Thus we conclude that
$$
\lim_{t \to 0}\langle \omega(\gX_t, \{\alpha(t), \beta(t)\})|=
\langle \omega(\gX_0, \{\alpha(0), \beta(0)\})|.
$$

Next, assume that
the curve $\widetilde{C}_0$ has two connected components $C_1$ and $C_2$,  and
that $Q_1$, $P_+ \in C_1$ and $Q_2$, $P_- \in C_2$. Moreover,  it is enough to consider
the case $N=2$.  Let $\{\omega_1, \ldots, \omega_{g_1}\}$ and $\{\omega_{g_1+1},
\ldots, \omega_{g_1+g_2}\}$ be normalized bases of holomorphic one-forms
of $C_1$  and $C_2$, respectively.  Then we can find a family
$\{\omega_1(t), \ldots, \omega_{g_1+g_2}(t)\}$
of normalized bases of the family $\pi: \cC^*=
\pi^{-1}(D \setminus \{0\}) \rightarrow
D\setminus \{0\}$ such that $\lim_{t \to 0}\omega_k(t) = \omega_k$ according
to Proposition 3.1 \cite{F}. Here
we regard $\omega_k$ as a holomorphic section of the dualizing sheaf
$\omega_{C_0}$ of  the nodal curve, by
extending it by zero to the other component.

Over $D\setminus \{0\}$
$\langle \omega(\gX_t, \{\alpha(t), \beta(t)\})|$ is a holomorphic section
 of the sheaf of ghost vacua.
Let $1=m_1<m_2<\cdots<m_{g_1}\le 2g_1-1$
and $1=n_1<n_2<\cdots<n_{g_2}\le 2g_2-1$ be the Weierstrass gap values of $C_1$
at $Q_1$ and $C_2$ at $Q_2$, respectively. Put
$$
| u_1\rangle = \overline{e}_{m_{g_1} -1/2}\wedge \cdots\wedge
\overline{e}_{m_1 -1/2}\wedge |-1\rangle,
\quad | u_2\rangle =\overline{e}_{n_{g_2} -1/2}\wedge \cdots\wedge
\overline{e}_{n_1 -1/2}\wedge |0\rangle.
$$
Choose sections $f_j$
of $\pi_* \cO_\cC(n_j s_2(D)+*s_1(D))$ over $D$ such that it has a Laurent
expansion
$$
f_{j, 2}= \frac{1}{\xi_2^{n_j}} + \hbox{\rm holomorphic}
$$
along $s_2(D)$.
Put also
$$
\gX_t' = (\pi^{-1}(t);  s_1(t);  \xi_1).
$$
Then we have that
\begin{eqnarray*}
\langle \omega(\gX_t, \{\alpha(t), \beta(t)\})|u_1\otimes u_2 \rangle
&=&
\langle \omega(\gX_t, \{\alpha(t), \beta(t)\})|u_1\otimes
\ovpsi[f_{g_2,2}]\ovpsi[f_{g_2-1,2}]\cdots \ovpsi[f_{1,2}]|0\rangle\\
&=& (-1)^{g_2-1}
\langle \omega(\gX_t, \{\alpha(t), \beta(t)\})|
\ovpsi[f_{1,2}]\ovpsi[f_{2,2}]\cdots \ovpsi[f_{g_2,2}]|u_1 \otimes 0\rangle\\
&=&  (-1)^{g_2-1}
\langle \omega(\gX_t', \{\alpha(t), \beta(t)\})|
\ovpsi[f_{1,2}]\ovpsi[f_{2,2}]\cdots \ovpsi[f_{g_2,2}]|u_1 \rangle\\
&=& (-1)^{g_2-1}
\left|
\begin{array}{ccc}
\res_{Q_1}(f_1\omega_1(t))& \cdots &
\res_{Q_1}(f_1\omega_{g_1+g_2}(t))\\
\vdots&\ddots&\vdots \\
\res_{Q_1}(f_{g_2}\omega_1(t))& \cdots &
\res_{Q_1}(f_{g_2}\omega_{g_1+g_2}(t))\\
I_{m_1}^{1}(Q_1, t) &\cdots & I_{m_1}^{g_1+g_2}(Q_1, t)\\
\vdots&\ddots&\vdots \\
I_{m_{g_1}}^{1}(Q_1, t) & \cdots & I_{m_{g_1}}^{g_1+g_2}(Q_1, t)
\end{array}\right|
\end{eqnarray*}
where the Taylor expansion of $\omega_k(t)$ along $s_1(D)$ and $s_2(D)$ are
written as
\begin{eqnarray*}
\omega_k(t) &=& \sum_{m=1}^\infty I_m^k(Q_1, t)\xi_1^{m-1} \\
\omega_k(t) &=& \sum_{n=1}^\infty I_n^k(Q_2, t)\xi_1^{n-1}
\end{eqnarray*}
Now by the residue theorem we have
$$
\res_{Q_1}(f_j\omega_k(t)) = - \res_{Q_2}(f_j\omega_k(t))
= - I_{n_j}^k(Q_2, t).
$$
Since for any positive integers $m$, $n$ we have
\begin{eqnarray*}
\lim_{t \to 0}I_m^k(Q_1, t) &=& I_{m}^k(Q_1), \quad 1 \le k \le g_1,\\
\lim_{t \to 0}I_n^{g_1+i}  (Q_1, t) &=&  0, \quad 1 \le i \le g_2, \\
\lim_{t \to 0}I_m^k (Q_2, t) &=&  0, \quad 1 \le k \le g_1,\\
\lim_{t \to 0}I_n^{g_1+i} (Q_2, t) &=&  I_n^{g_1+i}(Q_2), 1 \le i \le g_2,
\end{eqnarray*}
we conclude that
\begin{eqnarray*}
&&
\lim_{t\to 0} \langle \omega(\gX_t, \{\alpha(t), \beta(t)\})|u_1\otimes u_2
\rangle \\
&&=
(-1)^{g_2-1} \lim_{t\to 0}
\left|
\begin{array}{ccc}
-I_{n_1}^1(Q_2, t)& \cdots &
-I_{n_1}^{g_1+g_2}(Q_2, t)\\
\vdots&\ddots&\vdots \\
-I_{n_{g_2}}^1(Q_2, t)& \cdots &
-I_{n_{g_2}}^{g_1+g_2}(Q_2, t)\\
I_{m_1}^{1}(Q_1, t) &\cdots & I_{m_1}^{g_1+g_2}(Q_1, t)\\
\vdots&\ddots&\vdots \\
I_{m_{g_1}}^{1}(Q_1, t) & \cdots & I_{m_{g_1}}^{g_1+g_2}(Q_1, t)
\end{array}\right| \\
&&=-
\left|
\begin{array}{ccc}
I_{m_1}^1(Q_1)&\cdots & I_{m_1}^{g_1}(Q_1) \\
\vdots&\ddots&\vdots \\
I_{m_{g_1}}^1(Q_1) & \cdots & I_{m_{g_1}}^{g_1}(Q_1)
\end{array}
\right| \cdot  \left|
\begin{array}{ccc}
I_{n_1}^{g_1+1}(Q_2)&\cdots & I_{n_1}^{g_1+g_2}(Q_2)\\
\vdots&\ddots&\vdots \\
I_{n_{g_2}}^{g_1+1}(Q_2) & \cdots & I_{n_{g_2}}^{g_1+g_2}(Q_2)
\end{array}
\right| \\
&&= -\langle\omega(\gX_0,\{\alpha(0), \beta(0)\})|u_1\otimes u_2 \rangle \ne 0.
\end{eqnarray*}

Hence if $\lim_{t\to 0} \langle \omega(\gX_t, \{\alpha(t), \beta(t)\})|$
exists, then it is $- \langle \omega(\gX_0, \{\alpha(0), \beta(0)\})| $.

Now let us show that $\lim_{t\to 0} \langle \omega(\gX_t,
\{\alpha(t), \beta(t)\})|$ exists. Note that for $t \ne 0$,
$\langle \omega(\gX_t', \{\alpha(t), \beta(t)\})|$ is given by
$$
\langle \cdots \wedge e(\omega_{m}(t) )\wedge \cdots\wedge
e(\omega_2(t))\wedge e(\omega_1(t)) |,
$$
where $\{\omega_1(0),\ldots, \omega_g(0)\}$
is a normalized basis of holomorphic
one-forms of $C_0$ and $\omega_{g+i}(0)$ is a normalized
meromorphic one-form which has only a pole at $s_1(0)=Q_1$. Now for
$t \ne 0$ $\{\omega_1(0), \ldots, \omega_g(0)\}$ is defined by
using Theorem \ref{thm3.4}. This means that for $|v_1\otimes v_2
\rangle \in \cF_2$, the evaluation $\langle \omega(\gX_t,
\{\alpha(t), \beta(t)\})|v_1\otimes v_2 \rangle$ is reduced to the
evaluation of $\langle \omega(\gX_t', \{\alpha(t),
\beta(t)\})|v_1'\rangle$ so that we can use similar arguments
as above. For example let us calculate $\langle \omega(\gX_t,
\{\alpha(t), \beta(t)\})|v_1\otimes v_2 \rangle$ for
\begin{eqnarray*}
v_1&=& \oe_{k_1+1/2} \wedge |0\rangle, \quad 1 \le k_1 , \\ v_2&=&
\oe_{l_g-1/2} \wedge \cdots \wedge \oe_{l_g-1/2} \wedge
i(\oe_{-m-1/2})|-1\rangle, \quad 0 \le l_1<\cdots <l_g, \quad 0
\le m.
\end{eqnarray*}
Choose meromorphic functions $f_i \in H^0(D,
\pi_*\cO_\cC(*(s_1(D)+s_2(D)))$ which have the Laurent expansion
along $s_2(D)$: $$
 f_{i,2} = \frac{1}{\xi_2^{l_i}} + \hbox{\rm holomorphic}.
$$ Also choose $\tau \in H^0(D, \pi_*\omega_\cC(*(s_1(D)+s_2(D)))$
which has the Laurent expansion along $s_2(D)$: $$ \tau_2 =
\big(\frac{1}{\xi_2^m }+ \hbox{\rm holomorphic}\big)d\xi_2. $$
Then by a similar argument as above
$$ \langle
\omega(\gX_t, \{\alpha(t), \beta(t)\})|v_1\otimes v_2 \rangle =
\pm \langle \omega(\gX_t', \{\alpha(t), \beta(t)\})|
\psi[\tau_1]\ovpsi[f_{1,1}]\cdots \ovpsi[f_{g,1}]|-1\rangle. $$
By
a simple calculation we have that
\begin{eqnarray*}
\psi[\tau_1]\ovpsi[f_{1,1}]\cdots \ovpsi[f_{g,1}]|-1\rangle &=&
\sum_{i=1}^g (-1)^{i-1}\res_{s_1(t)}(f_i\tau)
\ovpsi[f_{1,1}]\cdots \ovpsi[f_{i-1,1}]| \ovpsi[f_{i+1,1}] \cdots
|-1\rangle\\ &&\ovpsi[f_{1,1}]\cdots
\ovpsi[f_{g,1}]\psi[\tau_1]|-1\rangle.
\end{eqnarray*}
Moreover, we have that $$ \psi[\tau_1]|-1\rangle = \sum_{n=1}^s
a_n \oe_{-3/2}\wedge \cdots \wedge \widehat{\oe_{-n-1/2}}\wedge
\oe_{-n_3/2}\wedge \cdots, $$ Then, as above
\small
\begin{eqnarray*}
&&\langle \omega(\gX_t', \{\alpha(t), \beta(t)\})|
\psi[\tau_1]\ovpsi[f_{1,1}]\cdots \ovpsi[f_{g,1}]|-1\rangle\\ &&=
\pm \sum_{i=1}^g (-1)^{i-1} \res_{s_1(t)}(f_i\tau) \left|
\begin{array}{ccc}
\res_{s_1(t)}(f_1\omega_1(t))& \cdots &
\res_{s_1(t)}(f_1\omega_g(t))\\ \vdots&\ddots&\vdots \\
\res_{s_1(t)}(f_g\omega_1(t))& \cdots &
\res_{s_1(t)}(f_g\omega_g(t))
\end{array}\right| \\
&& + \sum_{n=1}^s \left( \left|
\begin{array}{cccc}
\res_{s_1(t)}(f_1\omega_{g+n}(t))& \res_{s_1(t)}(f_1\omega_2(t))&
\cdots & \res_{s_1(t)}(f_1\omega_g(t))\\ \vdots&\ddots&\vdots \\
\res_{s_1(t)}(f_g\omega_{g+n}(t))&\res_{s_1(t)}(f_g\omega_2(t))&
\cdots & \res_{s_1(t)}(f_g\omega_g(t))
\end{array}\right|   \right. \\
&&+ \left|
\begin{array}{cccc}
\res_{s_1(t)}(f_1\omega_1(t))& \res_{s_1(t)}(f_1\omega_{g+n}(t))&
\cdots & \res_{s_1(t)}(f_1\omega_g(t))\\ \vdots&\ddots&\vdots \\
\res_{s_1(t)}(f_g\omega_1(t))&\res_{s_1(t)}(f_g\omega_{g+n}(t))&
\cdots & \res_{s_1(t)}(f_g\omega_g(t))
\end{array}\right|   \\
&& + \left. \left|
\begin{array}{cccc}
\res_{s_1(t)}(f_1\omega_1(t))& \res_{s_1(t)}(f_1\omega_2(t))&
\cdots & \res_{s_1(t)}(f_1\omega_{g+n}(t))\\ \vdots&\ddots&\vdots
\\ \res_{s_1(t)}(f_g\omega_1(t))&\res_{s_1(t)}(f_g\omega_2(t))&
\cdots & \res_{s_1(t)}(f_g\omega_{g+n}(t))
\end{array}\right|  \right)
\end{eqnarray*}
\normalsize
which can be rewritten in the following form by using
$\res_{s_1(t)}(f_i\omega_k) =-\res_{s_2(t)}(f_i\omega_k)$.
\small
\begin{eqnarray*}
&& = \pm \sum_{i=1}^g (-1)^{i-1} \res_{s_1(t)}(f_i\tau) \left|
\begin{array}{ccc}
-\res_{s_2(t)}(f_1\omega_1(t))& \cdots &
-\res_{s_2(t)}(f_1\omega_g(t))\\ \vdots&\ddots&\vdots \\
-\res_{s_2(t)}(f_{g_1}\omega_1(t))& \cdots &
-\res_{s_2(t)}(f_{g_1}\omega_g(t))\\
\res_{s_1(t)}(f_g\omega_{g_1+1}(t))& \cdots &
\res_{s_1(t)}(f_g\omega_{g_1+1}(t))\\ \vdots&\ddots&\vdots \\
\res_{s_1(t)}(f_g\omega_1(t))& \cdots &
\res_{s_1(t)}(f_g\omega_g(t))
\end{array}\right| \\
&& + \sum_{n=1}^s \left( \left|
\begin{array}{cccc}
-\res_{s_2(t)}(f_1\omega_{g+n}(t))&
-\res_{s_2(t)}(f_1\omega_2(t))& \cdots &
-\res_{s_2(t)}(f_1\omega_g(t))\\ \vdots&\ddots&\vdots \\
-\res_{s_2(t)}(f_{g_1}\omega_{g+n}(t))&-\res_{s_2(t)}(f_{g_1}\omega_2(t))&
\cdots & -\res_{s_2(t)}(f_{g_1}\omega_g(t))\\
\res_{s_1(t)}(f_{g_1+1}\omega_{g+n}(t))&
\res_{s_1(t)}(f_{g_1+1}\omega_2(t))& \cdots &
\res_{s_1(t)}(f_{g_1+1}\omega_g(t))\\ \vdots&\ddots&\vdots \\
\res_{s_1(t)}(f_g\omega_{g+n}(t))& \res_{s_1(t)}(f_g\omega_2(t))&
\cdots & \res_{s_1(t)}(f_g\omega_g(t))
\end{array}\right|   \right. \\
\end{eqnarray*}
\begin{eqnarray*}
&&+ \left|
\begin{array}{cccc}
-\res_{s_2(t)}(f_1\omega_1(t))&
-\res_{s_2(t)}(f_1\omega_{g+n}(t))& \cdots &
-\res_{s_2(t)}(f_1\omega_g(t))\\ \vdots&\ddots&\vdots \\
-\res_{s_2(t)}(f_{g_1}\omega_1(t))&-\res_{s_2(t)}(f_{g_1}\omega_{g+n}(t))&
\cdots & -\res_{s_2(t)}(f_{g_1}\omega_g(t)) \\
\res_{s_1(t)}(f_{g_1+1}\omega_1(t))&
\res_{s_1(t)}(f_{g_1+1}\omega_{g+n}(t))& \cdots &
\res_{s_1(t)}(f_{g_1+1}\omega_g(t)) \\ \vdots&\ddots&\vdots \\
\res_{s_1(t)}(f_g\omega_1(t))&\res_{s_1(t)}(f_g\omega_{g+n}(t))&
\cdots & \res_{s_1(t)}(f_g\omega_g(t))
\end{array}\right|   \\
&& + \left. \left|
\begin{array}{cccc}
\res_{s_1(t)}(f_1\omega_1(t))& \res_{s_1(t)}(f_1\omega_2(t))&
\cdots & \res_{s_1(t)}(f_1\omega_{g+n}(t))\\ \vdots&\ddots&\vdots
\\
\res_{s_1(t)}(f_{g_1}\omega_1(t))&\res_{s_1(t)}(f_{g_1}\omega_2(t))&
\cdots & \res_{s_1(t)}(f_{g_1}\omega_{g+n}(t))\\
\res_{s_1(t)}(f_{g_1+1}\omega_1(t))&\res_{s_1(t)}(f_{g_1+1}\omega_2(t))&
\cdots & \res_{s_1(t)}(f_{g_1+1}\omega_{g+n}(t))\\
\vdots&\ddots&\vdots \\
\res_{s_1(t)}(f_g\omega_1(t))&\res_{s_1(t)}(f_g\omega_2(t))&
\cdots & \res_{s_1(t)}(f_g\omega_{g+n}(t))
\end{array}\right|  \right) .
\end{eqnarray*}
\normalsize
Now it is easy to see that the last expression has a meaning when
$t$ goes to 0.

Hence $\lim_{t \to 0} \langle \omega(\gX_t, \{\alpha(t),
\beta(t)\})|v_1\otimes v_2 \rangle$ always exist. \QED


\begin{thebibliography}{KNTY}

\bibitem[AU2]{AU2} J. Andersen \& K. Ueno, {\em Geometric construction of
 modular functors from conformal field theory}, MPS-preprint 2003 -- 5. math.QA/0304135.

\bibitem[AU3]{AU3} J. E. Andersen \& K. Ueno, {\em Construction of the Reshetikhin-Turaev
TQFT from conformal field theory}, Preprint in preparation.

\bibitem[BPZ]{BPZ} A. A. Belavin,  A. M. Polyakov \& A. B. Zamolodchikov,
{\em Infinite conformal symmetry in two-dimensional quantum
field theory}, Nucl. Phys.{\bf B241} (1984), 338 -- 380.

\bibitem[F]{F} J. D. Fay, {\em Theta functions on Riemann surfaces}, Lecture Notes in Math.
{\bf 352} (1973), Springer-Verlag.

\bibitem[FP]{FP} L. D. Faddeev \& V. N. Popov, {\em Feynman diagrams for the Yang-Mills
fields}, Physics Letters, {\bf B25} (1967), 29 -- 30.

\bibitem[KNTY]{KNTY} N. Kawamoto,  Y. Namikawa,
A. Tsuchiya \& Y. Yamada,
{\em Geometric realization of conformal field theory
on Riemann surfaces},
 Commun. Math. Phys.,{\bf 116} (1988), 247 -- 308.

\bibitem[KSUU]{KSUU} Kuroda, Y. Shimizu, Uematsu \& K. Ueno, {\em Abelian
conformal field theory under degenerations of curves}, Preprint, 2002.

\bibitem[SA]{SA} M. Sato, {\em Soliton equations as dynamical systems on
an infinite dimensional Grassmann manifold}, Res. Inst. Math. Sci., Kyoto
Univ.-- Kokyuroku {\bf 439}, 30 -- 46.

\bibitem[TK]{TK} A. Tsuchiya and Y. Kanie,
{\em Vertex operators in conformal field theory on ${\bf P}^1$
and monodromy representations of braid group},
Advanced Studies in Pure Math.,
{\bf 16 } (1988),  297 -- 372, Erratum, {\em ibid}.,
{\bf 19 }(1989), 675 -- 682.

\bibitem[TUY]{TUY} A.Tsuchiya,
K. Ueno  \&  Y. Yamada
{\em Conformal field theory on universal
family of stable curves with gauge symmetries},
Adv. Stud. in Pure Math., {\bf 19} (1989), 459--566

\bibitem[U1]{U1} K. Ueno, {\em On conformal
field theory}, London Math. Soc. Lecture Notes Series
{\bf 208} (1995), 283 -- 345, Cambridge Univ. Press.

\bibitem[U2]{U2}K.  Ueno,  {\em Introduction to conformal field
theory with gauge symmetry}, Lecture notes in pure and applied
mathematics {\bf 184},  Geometry and Physics,   603 -- 745,
Marcel Dekker, 1996.

\end{thebibliography}
\end{document}